\documentclass{article}
\usepackage{arxiv}
\usepackage[utf8]{inputenc}
\usepackage[T1]{fontenc}

\usepackage[authoryear,longnamesfirst]{natbib}
\usepackage{amsmath,amsfonts,amssymb}
\usepackage{url,nicefrac,microtype,doi}
\newcommand{\defeq}{\stackrel{\scalebox{0.5}{\ensuremath{\mathsf{\,def}}}}{=}}
\usepackage{booktabs,makecell,threeparttable,array}
\usepackage{multirow} % for \multirow in tables

\usepackage{graphicx}  % \scalebox
\usepackage{placeins}  % keep result floats close to their discussion

\usepackage{amsthm}
\newtheorem{observation}{Observation}

\newtheorem{proposition}{Proposition}
\newtheorem{deploymentrule}[proposition]{Rule}
\usepackage{algorithm}
\usepackage{algorithmic} % provides \REQUIRE/\ENSURE/\STATE etc.

\usepackage{subcaption}

\DeclareMathAlphabet{\pazocal}{OMS}{zplm}{m}{n}

\usepackage{hyperref}
\usepackage{xurl}

\long\def\rev#1#2#3{\ignorespaces#3\ignorespaces\unskip}
\begin{document}
% Correct float-fraction parameter names (previous \floatpagepagefraction /
% \textpagefraction were typos of the real LaTeX kernel commands and had no
% effect, leaving the strict defaults in place and forcing large figures onto
% dedicated float pages, leaving the preceding text page mostly blank).
\renewcommand{\topfraction}{0.9}
\renewcommand{\bottomfraction}{0.8}
\renewcommand{\textfraction}{0.07}
\renewcommand{\floatpagefraction}{0.75}
% Allow more floats per page so figures are not deferred behind tables
% (defaults: topnumber=2, bottomnumber=1, totalnumber=3).
\setcounter{topnumber}{3}
\setcounter{bottomnumber}{2}
\setcounter{totalnumber}{5}

\renewcommand{\shorttitle}{Predict-then-Optimize Public Transport Line Redesign}

\title{Predict-then-Optimize Framework for Public Transport Line Redesign under Fluctuating Traffic Conditions}

\author{
Zihao GUO \\
Institut Polytechnique de Paris \\
Palaiseau 91120, France \\
\texttt{zihao-eric.guo@ip-paris.fr}
\And
Andrea ARALDO\thanks{Corresponding author: \texttt{andrea.araldo@telecom-sudparis.eu}} \\
Institut Polytechnique de Paris \\
Palaiseau 91120, France \\
\texttt{andrea.araldo@telecom-sudparis.eu}
\And
Fay\c{c}al TOUZOUT \\
National School of State Public Works \\
Marne-la-Vall\'{e}e 77454, France \\
\texttt{faycal.touzout@univ-eiffel.fr}
\And
Mounim EL-YACOUBI \\
Institut Polytechnique de Paris \\
Palaiseau 91120, France \\
\texttt{mounim.el\_yacoubi@telecom-sudparis.eu}
}

\date{}

\hypersetup{
  hidelinks,
  pdftitle={Predict-then-Optimize Framework for Public Transport Line Redesign under Fluctuating Traffic Conditions},
  pdfauthor={Zihao GUO, Andrea ARALDO, Fay\c{c}al TOUZOUT, Mounim EL-YACOUBI},
  pdfkeywords={Public Transport network redesign, Predict-then-optimize, Multi-objective optimization, Travel time uncertainty, Structural continuity}
}

\maketitle

% Here goes the abstract
\begin{abstract}
Public Transport (PT) lines are traditionally designed to optimize performance under nominal traffic conditions. In practice, however, operating conditions frequently deviate from the nominal ones, leading to substantial performance deterioration. Existing adaptation mechanisms typically rely on corrective and reactive interventions, such as stop-skipping, which prove insufficient under large or recurrent traffic fluctuations. In such contexts, incremental adjustments may be insufficient. This paper evaluates the potential of deeper structural \emph{redesigns} to preserve performance.
To this aim, this paper proposes a method to proactively redesign appropriate parts of PT networks when faced with high traffic fluctuations that would otherwise deteriorate operator- and user- performance. We adopt a predict-then-optimize paradigm in which PT lines are reconfigured based on traffic forecasts via an algorithm based on the Non-dominated Sorting Genetic Algorithm~III (NSGA-III).
For operational feasibility, to avoid excessive structural changes, we enforce high Jaccard edge overlap between the original and redesigned networks. To account for the effect of prediction inaccuracies on redesign quality, we construct a statistical model of the errors made by a well-established deep learning-based prediction model (Diffusion Convolutional Recurrent Neural Network) trained on real-world data.
Computational results on Mandl's benchmark network and on the large-scale Beijing network show that the controlled redesign of PT lines brings significant user-centric performance improvement and operational cost reductions in high traffic fluctuation scenarios, while limiting topological modifications. 
On the Beijing network under high variability, average travel time improves by up to 25.8\% while preserving over 85\% line overlap; unlike stop-skipping baselines, which break connectivity for a large share of OD pairs, our redesign preserves full OD connectivity.
This indicates the potential of transitioning from current static planning paradigms toward continuous and adaptive PT network design.
\end{abstract}

% Use if graphical abstract is present
%\begin{graphicalabstract}
%\includegraphics{}
%\end{graphicalabstract}

% Keywords
% Each keyword is separated by \sep
\keywords{Public Transport network redesign \and Predict-then-optimize \and Multi-objective optimization \and Travel time uncertainty \and Structural continuity}

% Main text
% =============================================
% 1.Introduction
% =============================================
\section{Introduction}

Public Transport (PT) networks are typically planned for nominal, scheduled conditions~\citep{ibarra2015planning}, while real-world operations are continuously subject to variability caused by congestion and other factors~\citep{buchel2020review, danaher2020minutes}, especially in mixed traffic environments. Such routine perturbations occur daily and can propagate across the network, progressively degrading service reliability~\citep{marra2020delay}. For example, in Xi’an, 74.34\% of bus runs were delayed~\citep{zhang2022travel}. 
Beyond operational impacts, insufficiently reliable public transport can also induce behavioral responses, including increased reliance on private vehicles~\citep{rodier2025evaluation}.

Existing approaches to cope with recurrent perturbations can be broadly classified as \emph{preventive}, \emph{reactive}, or \emph{proactive}, depending on the decision horizon and intervention level. 
Preventive strategies enhance robustness ex ante through structural measures, such as dedicated bus lanes that isolate buses from mixed traffic~\citep{russo2022dedicated}. 
However, limited street space and high infrastructure costs often make such measures infeasible in dense urban areas~\citep{li2009evaluation}, and preventive capacity tends to be inflexible and underutilized~\citep{ITF2024, ITF2025}.
Reactive strategies attempt, at the operational level, to restore service after perturbations have occurred, for instance via stop-skipping, holding, or detouring~\citep{peled2021quality, tang2018optimal, zheng2023multi}. 
These actions deteriorate under high uncertainty, inconvenience passengers~\citep[\S{}2]{GkiotsalitisCats2021}, and are insufficient to manage cascading perturbations when no alternative network structure is available~\citep{mahdavi2024resilience}.

A key limitation shared by both preventive and reactive strategies is that the PT network structure typically remains fixed, despite persistent changes in operating conditions. 
Network design decisions are revised infrequently, and traditional models assume deterministic or fixed travel times~\citep{Fonseca2020, BuechelCorman2020}, limiting their ability to address recurrent traffic fluctuations.
Proactive approaches aim to overcome this limitation by anticipating variability through short-term predictions and revising the line design before the next operating interval~\citep{Peled2019}. 
%
%\rev{zg}{}{It should be noted that proactive network redesign, as proposed in this paper, also introduces a form of user inconvenience for passengers accustomed to fixed route patterns; We minimize this by controlling structural similarity~$Z_3$ (see~\S{}\ref{sec:timeframe}).}

The \textbf{main contribution} of this paper is to assess the potential of predict-then-optimize proactive redesign of bus networks. To this aim, we propose a new method that adapts bus lines to time-varying stochastic traffic congestion. We solve a multi-objective problem to minimize operational cost and user travel time, and to limit structural network modifications, so as to ease operations and user routing. Adaptation is triggered only when forecasted conditions severely deviate from the nominal ones; moreover, such adaptation is anticipated (e.g., 1h in advance) to avoid disrupting journeys of current travelers.

%Results on the Mandl benchmark and Beijing networks show that notable improvements in user- and operator-centric performance can be achieved while limiting the extent of structural modifications. On the large-scale Beijing network, average travel time improves by up to 25.8\% and operational cost by up to 22.6\% at highly variable traffic (coefficient of variation of link travel time = 100\%) while route overlap stays above 85\%;

\rev{zg}{}{Results on the Mandl benchmark and Beijing networks show that notable improvements in user- and operator-centric performance can be achieved while limiting the extent of structural modifications. In the Beijing large scale network, beyond a one-off static re-optimization gain of 11.4\%, with no link travel time variability, the benefits grow with such a variability. When it is very high, average travel time improves by up to 25.8\% and operational cost by up to 22.6\% at with high traffic variability, while line overlap stays above 85\% (Table~\ref{tab:stop_skipping_comparison}). Compared with a stop-skipping baseline, which improves travel time by at most 19.7\% but only by disconnecting 28--39\% of OD pairs under high variability, our redesign attains a larger gain while keeping every OD pair reachable.} These findings highlight the potential of proactive bus-line redesign in future public transport systems, where adaptive operational strategies may become easier to deploy thanks to vehicle automation~\citep{FagnantKockelman2015} and to users’ growing familiarity with flexible, app-mediated mobility services, as illustrated by the diffusion of on-demand platforms.

The remainder of the paper is organized as follows. 
Section~\ref{sec:relatedWork} reviews related work, Section~\ref{sec:metho} formulates the problem, Section~\ref{sec:alg} presents the solution approach, Section~\ref{sec:evaMETHO} describes the evaluation method, Section~\ref{sec:numericalRESULT} reports the results, and Section~\ref{sec:conclu} concludes.

% =============================================
% 2. Related Work
% =============================================
\section{Related Work}\label{sec:relatedWork}

%This section reviews the literature on public transport optimization according to decision levels and adaptation horizons.
Previous work related to ours typically focuses on either \textit{network design} at strategic and tactical levels, where PT lines are planned from scratch, or \textit{reactive operational control}, which manages daily operations without altering the network structure.
\textit{Network redesign} occupies an intermediate position by adapting existing networks rather than starting anew.
Within this domain, \textit{offline redesign} relies on re-optimization for known future scenarios, whereas emerging \textit{proactive redesign} approaches anticipate short-term uncertainty.
However, these proactive methods often overlook the critical trade-off between responsiveness and structural stability, a gap that this paper specifically addresses.

% % ++++++++++++++++++++ 2.1
\subsection{Public Transport Network Design under Nominal Conditions}
\label{sec:preventive}

The PT Network Design Problem (TNDP) consists of designing line routes and their frequencies, balancing user- and operator-centric objectives, and has been studied since the foundational formulations of~\cite{ceder1986bus} through successive reviews~\citep{guihaire2008transit, kepaptsoglou2009transit, duran2022survey} and reference texts~\citep{desaulniers2007public, ceder2016public},\cite{gkiotsalitis2020public}.

Most of these approaches formulate public transport networks \emph{from scratch}, without reference to any pre-existing service configuration; although some formulations begin from a predefined set of candidate lines~\citep{cervantes2023analyzing}, the network topology is typically treated as a fresh design decision rather than an adaptation of an existing operational configuration. These methods are primarily intended for long-term strategic or tactical planning.
The TNDP is NP-hard~\citep{magnanti1984network}.

User-centric performance has been expressed in numerous ways. In addition to overall travel duration~\citep{john2014improved,heyken2019adaptive}, recent studies also examine network coverage~\citep{cadarso2017rapid}, route efficiency~\citep{suman2019improvement}, and demand fulfillment~\citep{barahimi2021multi}. Operator-centric objectives are focused on limiting operational and capital expenditures~\citep{ahern2022approximate,chen2018continuum,calabro2023adaptive}, or proxies for them, such as the overall length of line routes.
Balancing these objectives draws on two distinct methodological strands. \emph{Multi-level} formulations capture hierarchical decision-making among different agents: \cite{szeto2014transit} established a bilevel framework in which an upper-level planner minimizes passenger transfers and a lower-level model handles PT assignment; \cite{tian2021service} extended this to a trilevel structure with passenger routing strategies at an intermediate level. \emph{Multi-objective} formulations (as ours) instead seek Pareto-optimal trade-off solutions across conflicting criteria: \cite{cervantes2023analyzing} analyzed the trade-off between travel-time reduction and fare reduction. 
%Although both paradigms address PT network design, they respond to fundamentally different modeling questions and should not be conflated.

Note that all these methods design PT network based on \emph{nominal} (e.g., average, peak, etc.) \emph{traffic conditions}.
Therefore, as noted by~\citep{tan2024optimizing}, the resulting networks may be inefficient faced with real-time variations in traffic patterns, deviating from nominal conditions. 

% ++++++++++++++++++++ 2.2
\subsection{Reactive Operational Control under Disruptions and Uncertainty}
\label{sec:reactive}

In most real-world operations, the structure of the public transport network remains fixed as originally planned. 
When deviations from nominal conditions occur due to exogenous factors~\citep{nissen2020does, chen2009analyzing}, 
operators rely on short-term corrective actions to mitigate supply--demand mismatches.
Typical examples include stop-skipping~\citep{peled2021quality}, short-turning, 
speed changing, holding~\citep{tang2018optimal, manasra2019optimization, gkiotsalitis2020public, olvera2023holding}, 
detouring of interrupted lines~\citep{zheng2023multi}, and replacement buses.
Each control action has distinct disadvantages.
Namely, bus holding increases the inconvenience of on-board passengers at the holding stops, short-turning requires passengers to disembark and wait for a subsequent bus trip, rescheduling affects the timetables and the crew/vehicle schedules; speed control cannot be implemented when the required speeds are beyond the safety limits, and stop-skipping results in refused passenger boardings~\citep{gkiotsalitis2021stop}.
While some works integrate real-time data into operational decision-making, such as schedule adjustment for enhanced reliability~\citep{cao2019real}, these frameworks focus on timetable or speed control within a fixed route structure and do not explicitly consider structural reconfiguration of the network topology under time-varying conditions.
It should be noted that proactive network redesign, as proposed in this work, may generate analogous and in some respects stronger inconveniences on the operational side, affecting vehicle scheduling and driver rostering, in particular with human-driven PT vehicles; it is thus important to limit excessive structural instability (as introduced in~\S{}\ref{sec:timeframe}).

More importantly, when actual conditions deviate significantly from the nominal ones, such corrective actions may not be sufficient to restore performance~\citep{guo2024restoring}, since the overall structure is maintained as planned.
Strategies such as stop-skipping and short-turning have also appeared in the context of offline or preplanned redesign, to generate alternative route configurations or long-term service plans under hypothetical disruption scenarios. They are outside the scope of this paper, which instead focuses on proactive redesigning lines for the hour ahead based on short-term forecasts.
%In contrast, within the operational control context considered here, they represent short-term corrective actions applied dynamically during daily operations to recover service regularity without altering the network topology.

% % ++++++++++++++++++++ 2.3
\subsection{Offline and Preplanned Network Redesign under Known Scenarios}
\label{sec:offline}

Various approaches can be employed in the bus route adjustment problem to generate adaptive bus routes, including rerouting, short-turning, skip-stopping, and branching~\citep{cao2019autonomous, ma2023integrated, wu2021modular}. 
In contrast to the network design approaches described in \S{}\ref{sec:preventive}, these studies explicitly start from an \emph{existing} public transport network and allow service or structural modifications under a limited set of anticipated scenarios.
Previous research on bus route localized adjustments has often aimed at two main objectives: accommodating demand fluctuations and enhancing integration with subway networks.

\cite{guan2023bi} introduced a bilevel planning model to optimize customized bus routes offline based on a spatiotemporal state network. While effective for static planning, this method is impractical for real-time dispatching. Concretely, it relies on a computationally expensive nested algorithm and requires a full re-solve for any new event.
\cite{zhang2018optimal} investigated long-distance replacement bus services to improve coordination with urban rail networks. While effective under static demand assumptions, this approach is inherently offline and cannot be applied in real-time operational contexts.
When passenger flow control is implemented in a subway system, \cite{zhou2021integrated} optimized existing bus routes to transport passengers temporarily restricted from entry, in case of disruptions.
\cite{wang2022optimization} proposed a two-stage method to mitigate disruptions to bus routes under planned construction scenarios, relying on predefined work-zone impacts and static travel-time estimates.
\cite{wu2024multi} proposed a multi-objective ant colony algorithm to adjust existing bus routes and introduced a sequence-level similarity to prevent the redesigned network from being excessively different from the original. However, their approach remains static and cannot capture time-varying adaptations.
\cite{zheng2024alternative} proposed localized adaptive bus routes under anticipated disruptions, while \cite{zhang2024improving} addressed large-scale emergency disruptions assuming available spare capacity.
These assumptions do not hold under recurrent traffic fluctuations, where no redundant lines exist and adaptation must occur within the same constrained topology.

% % ++++++++++++++++++++ 2.4
\subsection{Network Redesign under Predictive Uncertainty and Structural Consistency}
\label{sec:proactive}

Traditional network design handles uncertainty primarily by optimizing for expected scenarios (stochastic optimization) or worst-case bounds (robust optimization)~\citep{bertsimas2018data, wang2020robust}; see~\cite{ge2022robustness} for a recent review of robustness and disturbances in public transport.
However, these approaches typically determine a single network configuration before operations begin, with limited adaptation to realized road conditions (i.e., without reconfiguration). While this ensures stability and ease of operations, it may become ineffective in preserving performance faced with high traffic fluctuations.
%
%{However, these approaches typically determine a single network configuration prior to realization, not because they are unable to generate multiple configurations, but because operational contexts require stable, predictable network structures whose restructuring must be justified from an organizational standpoint~\citep{wu2024multi}. As a result, they provide no adaptation to the realized road conditions.}
%
Another line of work studies topological metrics to evaluate network resilience~\citep{Massobrio2024, Pu2022}, e.g., they measure vulnerability after a perturbation occurs. These approaches remain purely diagnostic and descriptive, rather than being prescriptive (as the approach we present in this paper), in that they do not provide guidance on the optimization of the network design. 
%
%Even modern predict-then-optimize frameworks~\citep{elmachtoub2022smart, mandi2024decision} overlook structural stability. They leverage short-term forecasts to adjust operational decisions—such as routing or scheduling—but ignore the underlying physical topology. Even if structural redesign is permitted, these models maximize performance under predicted conditions without preventing drastic, chaotic changes between redesign stages.

%In short, existing literature fundamentally isolates uncertainty, network design, and structural stability, lacking any approach that dynamically updates a physical network while enforcing structural similarity.
%
%To bridge this gap, our framework completely departs from this norm. We do not treat structural consistency as an \emph{ex post} evaluation metric or a weak penalty. Instead, we embed graph-level similarity directly into the optimization model. By forcing the new topology to remain close to the original configuration, we strictly constrain the redesign space under predictive uncertainty.

% ++++++++++++++++++++ 2.5

\subsection{Vehicle Routing Problems and Demand-Responsive Transport (DRT)}

In Vehicle Routing Problems (VRPs) and DRT, vehicle routes are constructed on-demand and solution stability has been employed to keep revised plans close to a baseline \citep{touzout2021modelling}. However, this stability applies to routing or scheduling assignments rather than the structure of the PT line network itself. Moreover, DRT and other VRP-based systems do not exhibit a structured network of lines, and therefore fail to achieve sufficient passenger consolidation~\citep{Calabro2022}. Hybrid fixed/flexible services partly bridge this gap~(\cite{calabro2023adaptive}), but still keep the fixed lines inflexible and flexible routes unstructured. By contrast, we study a high-capacity mass-transit system organized around a structured network of fixed lines, where existing lines are reconfigured rather than routes being computed from scratch. Our approach is therefore conceptually distinct from VRP- or DRT-based formulations. A direct comparison with DRT is also inappropriate, since we address high demand levels typical of dense urban cores, for which DRT is generally unsuitable in terms of both computational scalability and service performance~\citep{Calabro2022}.

%----

%In Demand-Responsive Transport (DRT) vehicle routes are constructed on-demand based on Vehicle Routing Problems (VRPs). \emph{Solution stability} concepts have been discussed, to limit plan revisions~(\citep{touzout2021modelling}). However, pure DRT generally lacks the fixed-line structure of high-capacity mass transit and therefore typically does not reach the same passenger consolidation levels~\citep{Calabro2022}. Hybrid fixed/flexible services partly bridge this gap~(\cite{calabro2023adaptive}), but still keep the fixed lines inflexible and flexible routes unstructured.
%Our approach operates on a high-capacity structured line network, reconfiguring existing routes rather than computing itineraries from scratch; it is therefore conceptually distinct from VRP/DRT-based formulations.

% ++++++++++++++++++++ 2.6

\subsection{Research gaps}
%Our literature review reveals several gaps in existing approaches to public transport network adaptation (\rev{zg}{see Table~\ref{tab:innovation_matrix} in the Appendix, which summarizes related work}
Table~\ref{tab:innovation_matrix} summarizes related work.
Most existing studies treat PT network design as a deterministic optimization problem, ignoring the inherent uncertainty and variability in urban traffic conditions. Only a few studies consider any form of uncertainty, and few redesign the network in response to short-term forecasts while also modeling travel-time variability in a systematic way.
An exception is~\cite{guo2025data}, where, however, (i)~future travel times on road links were assumed to be perfectly known (which is unrealistic), so the impact of prediction errors was not considered, while we model the error of a well established prediction model, (ii)~only a single simplified objective was adopted, while we solve a multiobjective optimization, and (iii)~their approach is tested only on a simple hypothetical scenario, while we consider real-world dataset to fit the prediction error models.

\begin{table*}[htbp]
\centering
\begin{threeparttable}
\resizebox{\textwidth}{!}{%
\begin{tabular}{ccccc}
\toprule
\textbf{Study}
& \textbf{Solution Method}
& \makecell{\textbf{Online/Proactive} \\ \textbf{Redesign}}
& \makecell{\textbf{Similarity} \\ \textbf{Control}}
& \textbf{Decision Focus} \\
\midrule
\citep{magnanti1984network} & Exact methods & No & No & ND \\
\citep{john2014improved} & Multi-objective GA & No & No & ND \\
\citep{heyken2019adaptive} & Scaled network GA & No & No & ND \\
\citep{cadarso2017rapid} & Branch-and-cut & No & No & ND \\
\citep{ahern2022approximate} & Gradient-based & No & No & ND \\
\citep{tan2024optimizing} & Congestion analysis & No & No & FS \\
\citep{guan2023bi} & Custom bilevel & No & No & OR \\
\citep{zhang2018optimal} & Tailored GA & No & No & OR \\
\citep{wang2022optimization} & Two-stage heuristic & No & No & DR \\
\citep{wu2024multi} & Ant Colony & No & NW & OR \\
\citep{zheng2024alternative} & Heuristic & No & No & DR \\
\citep{cao2019autonomous} & Skip-stop heuristic & No & No & VC \\
\citep{ma2023integrated} & Integer programming & No & No & OR \\
\citep{tang2018optimal} & Optimal control & No & No & VC \\
\citep{manasra2019optimization} & Branch-and-bound & No & No & VC \\
\citep{gkiotsalitis2020public} & Evolutionary & No & No & VC \\
\citep{olvera2023holding} & Holding optimization & No & No & VC \\
\citep{zheng2023multi} & Coordination heuristic & No & No & VC \\
\citep{guo2025data} & Exact methods & Partial & Penalty & NR \\
\midrule
\textbf{Our Approach}
& \textbf{NSGA-III based optimization}
& \textbf{Yes}
& \textbf{Jaccard}
& \textbf{NR} \\
\bottomrule
\end{tabular}%
}

\vspace{2pt}
\parbox{\textwidth}{\footnotesize
\textit{Notes:} Online/Proactive Redesign -- \textbf{No}: no online or proactive redesign; \textbf{Partial}: limited redesign capability (e.g., sampling-based or scenario-specific updates); \textbf{Yes}: proactive redesign enabled (e.g., predict-then-optimize).
Similarity Control -- \textbf{No}: no explicit similarity mechanism; \textbf{NW}: Needleman--Wunsch-based similarity; \textbf{Penalty}: penalty-based similarity; \textbf{Jaccard}: Jaccard-based similarity.
Decision Focus -- \textbf{ND}: from-scratch network design; \textbf{FS}: frequency setting; \textbf{OR}: offline or preplanned redesign of an existing network; \textbf{DR}: perturbation-oriented redesign or response; \textbf{VC}: vehicle or operational control; \textbf{NR}: network redesign under uncertain travel times.
}

\caption{Comparative analysis of public transport optimization approaches across decision scope, online/proactive redesign capability, and similarity control mechanisms.}
\label{tab:innovation_matrix}
\end{threeparttable}
\end{table*}

%\rev{zg}{Our work introduces a novel methodology for redesigning routes in response to predicted time-varying road conditions through controlled modifications that balance passenger service quality and operational efficiency.}{Our work introduces a methodology for redesigning PT lines for the hour ahead under fluctuating traffic: it uses short-term forecasts, is triggered when needed (i.e., when traffic conditions are highly perturbed causing user- and operator-performance degradation), and balances passenger travel time against operational cost through controlled structural changes.}

To our knowledge, this work is the first to systematically study the potential of proactively redesigning PT lines (e.g., for the hour ahead) under fluctuating traffic. Redesign decisions are based on short-term forecasts, are triggered when needed (i.e., when traffic conditions are highly perturbed causing user- and operator-performance degradation), and balance passenger travel time against operational cost through controlled structural changes.

% =============================================
% 3.Methodology
% =============================================
\section{Model and Problem Formulation}
\label{sec:metho}
\begin{figure}
\centering
\includegraphics[width=0.95\textwidth]{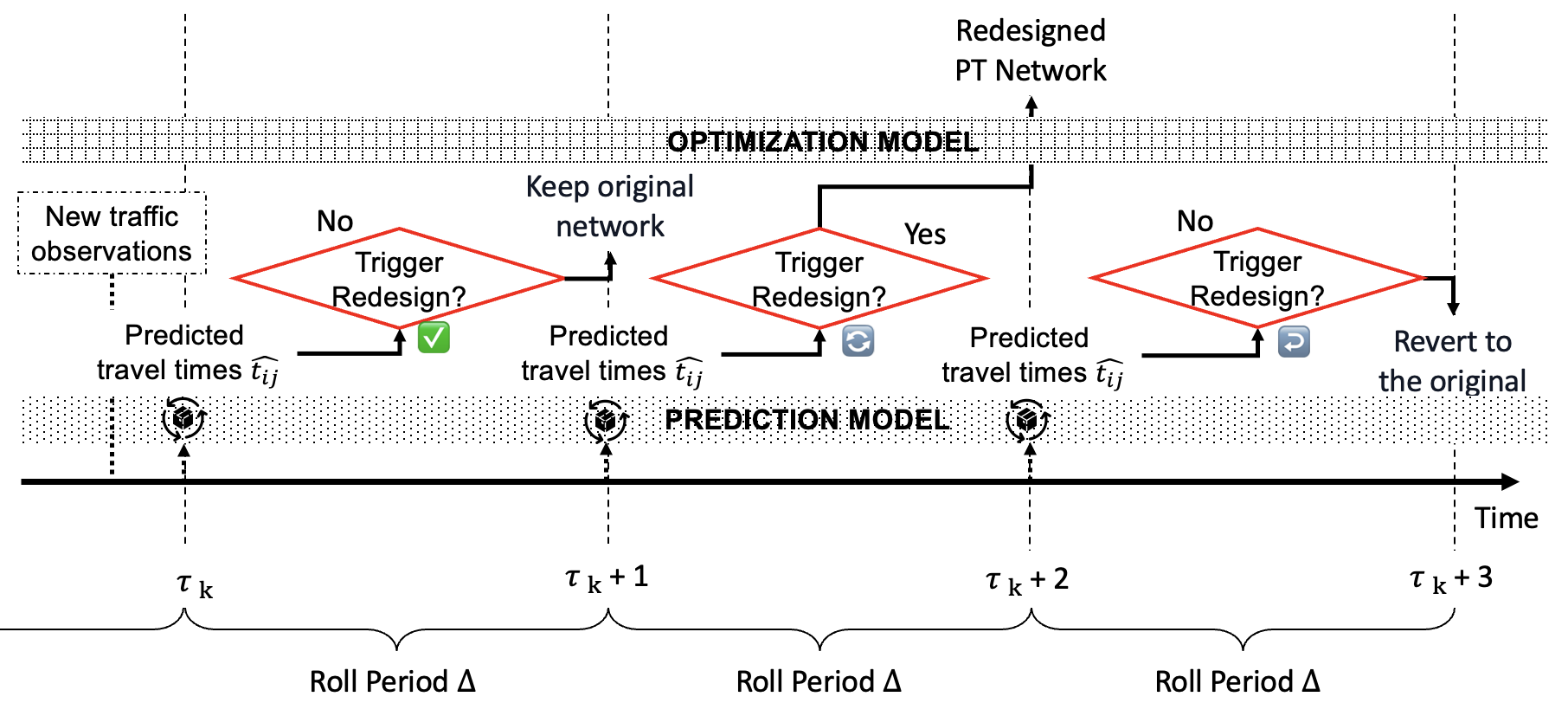}
\caption{\label{fig:rollFramework}Proactive Line Redesign Framework.}
\end{figure}

As illustrated in Fig.~\ref{fig:rollFramework}, our methodology follows a predict-then-optimize framework that explicitly integrates current observations into prediction. The prediction model is first trained offline on historical data to learn travel times on road links. During operation, it receives current observations from the latest time steps and outputs a short-term prediction of link travel times~$\hat{t}_{ij}$, which are then supplied to the optimization model for potential network redesign.
To ensure timely responsiveness, this procedure is repeated periodically, i.e., once per hour. At the start of each time slot, the new network is announced in advance (e.g., 1h), so passengers can plan a route that remains valid for their entire trip and is never changed mid-journey; we assume here that the longest trip is shorter than one time slot (1h). This also keeps the bus routes easier to manage. Further details are given in~\S\ref{sec:timeframe}.
The resulting network design is then evaluated against realized travel times $t_{ij}$ to assess the actual performance of our adaptive network under realistic conditions. The notation is summarized in Table~\ref{tab:notation}.

% ============================= FIGURE
% ============================= FIGUREend

% ++++++++++++++++++++ 3.1
\subsection{Graph Model of PT Network}\label{sec:graphMODEL}
% 物理基础设施（substrate） - Undirected
We represent the underlying road infrastructure as an undirected
substrate network $\pazocal{G}^{\text{sub}} = (\pazocal{V}^{\text{sub}}, \pazocal{E}^{\text{sub}})$, where $\pazocal{V}^{\text{sub}}$ are candidate bus stops and $\pazocal{E}^{\text{sub}} = \pazocal{V}^{\text{sub}} \times \pazocal{V}^{\text{sub}}$ the potential connections among them. Observe that a single $\ell \in \pazocal{E}^{\text{sub}}$ may represent a sequence of road links from a stop to another.
% 线路物理定义（service coverage） - Undirected
Let $\pazocal{L}$ denote the set of all PT lines in the network.
Each PT line $\ell \in \pazocal{L}$ is defined by a set of undirected substrate edges, representing the physical stop-to-stop connections served by the line.
% 数学编码 / 运行建模 - directed
To ensure path continuity and avoid disconnected or cyclic line configurations, the undirected edge set of each line is represented through a directed simple path graph $\pazocal{G}_\ell = (\pazocal{V}_\ell, \pazocal{E}_\ell)$ on the undirected substrate $\pazocal{G}^{\text{sub}}$,
where $\pazocal{V}_\ell \subseteq \pazocal{V}^{\text{sub}}$ is the set of stops served by line $\ell$, and $\pazocal{E}_\ell \subseteq \pazocal{E}^{\text{sub}}$ denotes the set of connections between consecutive stops along the path representation of line $\ell$.
\label{sentence:line-in-one-direction}
Vehicles operating on line $\ell$ are assumed to traverse this path in both directions; the reverse traversal is implicitly captured through symmetric travel times~\citep{schmid2014hybrid}.
The complete public transport network is constructed by aggregating all individual lines and is denoted by $\pazocal{G}_{\text{PT}} = (\pazocal{V}_{\text{PT}}, \pazocal{E}_{\text{PT}}, \pazocal{H})$.
In this formulation, $\pazocal{V}_{\text{PT}} = \bigcup_{\ell \in \pazocal{L}} \pazocal{V}_\ell$ and $\pazocal{E}_{\text{PT}} = \bigcup_{\ell \in \pazocal{L}} \pazocal{E}_\ell$, while $\pazocal{H} = \{h_\ell \mid \ell \in \pazocal{L}\}$ represents the set of line headways, with each $h_\ell$ defined as the time interval between two consecutive vehicles~\citep{vuchic1981timed}.
Accordingly, $\pazocal{G}_{\text{PT}}$ can be viewed as a line-labeled directed multigraph for modeling purposes, while operationally each line provides bidirectional service along its physical path.

For each PT line $\ell \in \pazocal{L}$, headway $h^\ell$ can be calculated by $h^\ell = \frac{T^\ell}{n^\ell}$, where $T^\ell$ represents the round trip time along the entire line, and $n^\ell$ denotes the number of vehicles operating on that line~(Eq.~4.4 of~\citep{vuchic1981timed}).  
Let $t^\ell$ denote the one-way travel time from the origin terminal to the destination terminal of line $\ell$, the round-trip cycle time is $T^\ell = 2t^\ell$.
Accordingly, the headway of line $\ell$ is given by
$h^\ell = \frac{2t^\ell}{n^\ell}$.
We use a fixed dwell time $d$, at each stop to account for the time a bus spends at a stop, including the time needed for passengers to board and alight, and the time lost for accelerating and decelerating.
While in practice dwell times vary with boarding and alighting counts, constant dwell time is standard at the network design level~\citep{HUYNH2022392, chen2022TRSC, schobel2022cheapest}, where demand is represented as an OD matrix rather than individual trip records and the optimization horizon (one hour) far exceeds individual dwell events. Variable dwell time is primarily relevant in operational-control models (bus bunching, holding) that operate at a finer temporal resolution.
Assuming random passenger arrivals, the average waiting time for line~$\ell$ is given by
\begin{equation}
    w^\ell = \frac{h^\ell}{2}\text{, where }h^\ell \text{denotes the service headway of line~}\ell.
    \label{eq:headwayTime}
\end{equation}
This assumption is standard in public transport modeling~\citep{parbo2018reducing, zhang2018service}. 
\cite{liu2020does} confirm that such an assumption is appropriate when the headway is relatively small (see their $\S$~3.3.1), and passengers simply go at the stop as soon as possible, without explicitly synchronizing at specific PT departure times; they also show that such a simple behavior is most time the most convenient for users, even better than a more complicated possible user strategies.
For a transfer at stop $i$ from line~$\ell$ to line~$\ell'$, the generalized transfer cost combines the expected waiting time
$w^{\ell'} = \frac{h^{\ell'}}{2}$
with a fixed 5~min transfer penalty, a common value in transit route system design~\citep{baaj1991ai}. Thus, each required transfer contributes $w^{\ell'}+5$~min, consistent with the waiting-time formulation in~\eqref{eq:headwayTime} while also accounting for transfer inconvenience.

% ++++++++++++++++++++ 3.1.1
\subsubsection{Virtual Node-Based Transfer Model}\label{sec:virtual_nodes}

\begin{figure}
   \centering
   \includegraphics[width=0.55\textwidth]{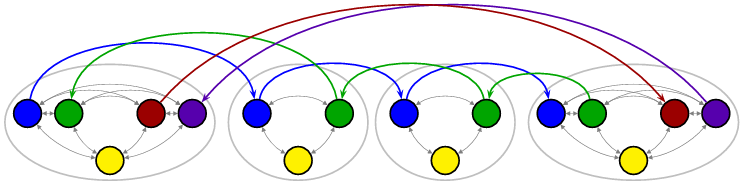}
   \caption{Example of a four-station PT network with asymmetric routes.}
   \label{fig:brt_example}
\end{figure}

To model passenger routing, we construct an auxiliary directed graph $\pazocal{G}^{\text{PAX}} = (\pazocal{V}^{\text{PAX}}, \pazocal{E}^{\text{PAX}})$, referred to as the virtual-node graph; the same object is denoted $\pazocal{G}'$ in the solution algorithm of \S\ref{sec:alg}.
The node set $\pazocal{V}^{\text{PAX}}$ comprises physical stop nodes $i \in \pazocal{V}^{\text{sub}}$ and line-specific virtual nodes $v_i^\ell$ for each line $\ell$ serving stop $i$. 
As illustrated in Fig.~\ref{fig:brt_example}, yellow nodes denote physical stops, while colored nodes represent line-specific virtual nodes; connections between colored nodes indicate potential passenger transfers.
The edge set $\pazocal{E}^{\text{PAX}}$ incorporates three distinct directed link types, each associated with a specific temporal cost. 
\textbf{(1) In-vehicle edges:} Connect virtual nodes $v_i^\ell$ and $v_j^\ell$ along the same line $\ell$. Their weights represent in-vehicle travel times.
\textbf{(2) Transfer edges:} Connect virtual nodes of different lines (e.g., $v_i^\ell$ to $v_i^{\ell'}$) at the same physical node. 
\textbf{(3) Access\& Egress edges:} Connect physical nodes to virtual nodes (e.g., $i$ to $v_i^\ell$). These represent initial boarding waiting times.
% \rev{zg}{}{\textbf{(4) Egress edges:} Connect each line-specific virtual node back to its physical stop (e.g., $v_i^\ell$ to $i$), allowing a passenger to alight from line $\ell$ at stop $i$. They carry zero cost, since alighting incurs no additional waiting, and they ensure that every OD path terminates at the physical destination node, so that shortest paths on $\pazocal{G}^{\text{PAX}}$ from $o_k$ to $d_k$ are well defined.}
\rev{zg}{The weights of (2) and (3) correspond to passenger waiting times, which are defined in~\eqref{eq:headwayTime} in $\S$~\ref{sec:graphMODEL}.}{The weight of each transfer edge in (2) equals the transfer waiting time plus the fixed 5~min transfer penalty defined above, whereas the weights of access/egress edges in (3) correspond to initial boarding waiting times; the waiting-time components are defined in~\eqref{eq:headwayTime} in $\S$~\ref{sec:graphMODEL}.}
Transfer edges are not created at origin or destination nodes.

% ========================= Figure
% ========================= FigureEND
% %
% To model transfers, for each physical node $i \in \pazocal{N}$ and each PT line $\ell \in \pazocal{L}$ serving it, a corresponding virtual node $v_{i}^\ell$ is generated. Virtual node identifier $v_i^\ell$ follows the encoding rule:
% $v_{i}^\ell = \ell \times 10^{s} + i$
% where $s$ is the minimum integer such that $10^s$ exceeds the maximum physical node identifier, ensuring uniqueness.

% % ++++++++++++++++++++ 3.3
\subsection{Stochastic Modeling of Travel Times}
\label{sec:Stochastic-Modeling-of-Travel-Times}
Let~$(\Omega,\pazocal F,\mathbb P)$ be a probability space. 
For the sake of realism, we model the aforementioned bus travel times~$t_{ij}$ from stop~$i$ to stop $j$, for each edge $(i,j)\in\pazocal{E}^{\text{sub}}$, at a given time of day, as random variables. For simplicity, we omit explicitly specifying the time of day in the notation. Let~$t_{ij}^\omega, \omega\in\Omega$ denote a realization of this variable. 
This is the travel time we would measure if we take a snapshot of the network at a certain point in time. Observe that, if we take another snapshot (at the same time of day, but in another day), we would observe another travel time~$t_{ij}^{\omega'}$, corresponding to another realization of variable~$t_{ij}$.
We assume that a prediction model provides predictions~$\hat t_{ij}$ for each edge.
If the prediction model were ``perfect'', for every sample~$\omega\in\Omega$, the realizations would satisfy $\hat t_{ij}^\omega = t_{ij}^\omega, \forall \omega\in\Omega,(i,j)\in\pazocal{E}^{\text{sub}}$. This is obviously impossible in reality, as any realistic prediction model is subject to some errors. Let the prediction error (which is another random variable) be:
\begin{equation}
\epsilon_{ij} \defeq \hat t_{ij}-t_{ij}\label{eq:link_time_error}
\end{equation}
%
%\rev{zg}{}{
%In the empirical residual-sampling procedure used in the experiments (\S\ref{sec:KDE}), we sample log-relative residuals $\eta_{ij}$ and construct predictions multiplicatively as $\hat t_{ij}=t_{ij}\exp(\eta_{ij})$. The additive error $\epsilon_{ij}=\hat t_{ij}-t_{ij}=t_{ij}(\exp(\eta_{ij})-1)$ may therefore be negative when the model underestimates actual travel time, while the predicted travel time remains positive by construction.
%}
% In practice, the PT planners consider baseline travel times $t_{ij}$ for design purposes (e.g., based on historical averages),
where $t_{ij}$ and $\hat t_{ij}$ denote the realized and predicted travel times, respectively. The sampling procedure for this error will be presented in \S{}\ref{sec:KDE}.
% These three quantities generally differ.
%
Let $e_{ij}^\ell$ indicate whether edge $(i,j)$ is used by line $\ell$. 
The time from one terminal to the other on line~$\ell$ is: 

\begin{align}
t^\ell
=
\underbrace{\sum_{i \in \pazocal{N}} \sum_{j \in \pazocal{N}} t_{ij} \, e_{ij}^\ell}_{\text{in-vehicle travel time}}
\;+\;
\underbrace{d \cdot \left(
\sum_{i \in \pazocal{N}} \sum_{j \in \pazocal{N}} e_{ij}^\ell - 1
\right)}_{\text{dwell time at stops}} .
\label{st:tl}
\end{align}
%
% Here, $d$ denotes a fixed dwell time per stop. Note that this~$t^\ell$ cannot be known in advance, because we cannot know in advance the true travel times~$t_{ij}$ on link~$(i, j)$. We can only calculate a ``nominal''~$t^{\ell,o}$ and a ``predicted''~$\hat t\ell$, obtained by replacing, in~\eqref{st:tl}, $t_{ij}$ with~$t_{ij}$ and $\hat t_{ij}$, respectively.
%
Here, $d$ denotes a fixed dwell time per stop. Since the actual link travel times $t_{ij}$ for edge $(i, j)$ are inherently uncertain, the true line travel time $t^\ell$ cannot be determined \textit{a priori}. We can only calculate a ``predicted''~$\hat t^\ell$, obtained by replacing $t_{ij}$ with~$\hat t_{ij}$ in~\eqref{st:tl}.
% --------------------------------------------------------------
% Notation Table
% --------------------------------------------------------------
\begin{table}
\centering
\renewcommand{\arraystretch}{0.95} % reduce row spacing
\setlength{\tabcolsep}{5pt} % adjust column spacing
\footnotesize % smaller font (you can also try \scriptsize if needed)
\begin{tabular}{p{0.35\textwidth} p{0.62\textwidth}}
\toprule
Symbol & Description \\
\midrule
\multicolumn{2}{l}{\textbf{Sets:}} \\
\midrule
$\pazocal{N}$ & Set of nodes, indexed by $i,j \in \pazocal{N}$ \\
$\pazocal{L}$ & Set of PT lines, indexed by $\ell \in \pazocal{L}$ \\
$\pazocal{K}$ & Set of origin–destination pairs, indexed by $k \in \pazocal{K}$ \\
$\pazocal{E}$ & Set of feasible edges in the network, indexed by $(i,j) \in \pazocal{E}$ \\
\midrule
\multicolumn{2}{l}{\textbf{Parameters:}} \\
\midrule
$q_k \in \mathbb{N},\; k \in \pazocal{K}$ & Demand for OD pair $k$ (passengers) - \S\ref{sec:obj}.
\\
% $\hat{t}_{ij}, t_{ij}^\text{nom} \in \mathbb{R}_{+},\; (i,j) \in \pazocal{E}$ & Predicted and baseline travel time on link $(i,j)$ \\
$t_{ij} \in \mathbb{R}_{+},\; (i,j) \in \pazocal{E}; \;\; t_{ij}^\omega, (i,j) \in \pazocal{E}, \omega\in\Omega$
& $t_{ij}$ is the true travel time on link $(i,j)$; this is a random variables; we denote its realizations as $t_{ij}^\omega$ - \S\ref{sec:Stochastic-Modeling-of-Travel-Times}. 
\\
$t_{ij}^0 \in \mathbb{R}_{+},\; (i,j) \in \pazocal{E}$ 
& Baseline link time (the ones referred to nominal conditions) - \S\ref{sec:trip-travel-times}. 
\\
$\hat{t}_{ij} \in \mathbb{R}_{+},\; (i,j) \in \pazocal{E}$ 
& Predicted travel time on link $(i,j)$ - \S\ref{sec:obj}.  
\\
$t_k^0 \in \mathbb{R}_{+},\; k \in \pazocal{K}$ 
& Baseline trip travel time: travel time for OD pair $k$ in the original network design $\pazocal G_\text{PT}^\text{org}$ under nominal (unperturbed) conditions - \S\ref{sec:trip-travel-times}.\\
$o_k, d_k \in \pazocal{N},\; k \in \pazocal{K}$ & Origin and destination node of OD pair $k$ - \S\ref{sec:trip-travel-times}. \\
$e_{ij}^{\ell, 0} \in \{0,1\},\; i,j \in \pazocal{N},\; \ell \in \pazocal{L}$ & Original network configuration (1 if link $(i,j)$ is used by line $\ell$) - \S\ref{sec:obj} \\
$M, d \in \mathbb{R}_{+}$ & Large constant, Fixed dwell time \\
\midrule
\multicolumn{2}{l}{\textbf{Upper-level decision variables (network design):}} \\
\midrule
\multicolumn{2}{l}{Binary Variables:} \\
\midrule
% 线路是否存在（建模层的 on/off 开关
$y^\ell,\; \ell \in \pazocal{L}$ & 1 if line $\ell$ is activated (operated) \\
% 线路的 stop-to-stop 序列（有向）
$e_{ij}^\ell,\; i,j \in \pazocal{N},\; \ell \in \pazocal{L}$ & 1 if edge $(i,j)$ is used by line $\ell$ in new network \\
% 存在时的唯一起点 / 终点
$a_i^\ell, b_i^\ell,\; i \in \pazocal{N},\; \ell \in \pazocal{L}$ & 1 if node $i$ is the start/end terminal of line $\ell$ \\
\midrule
\multicolumn{2}{l}{Discrete or Continuous Variables:} \\
\midrule
$n^\ell \in \mathbb{N},\; \ell \in \pazocal{L}$ & Number of vehicles on line $\ell$ (fleet size) 
\\

$u_i^\ell \in \mathbb{N},\; i \in \pazocal{N},\; \ell \in \pazocal{L}$ 
& Auxiliary variable for subtour elimination (MTZ): it denotes the order in which stop~$i$ appears in line $\ell$.
\\
$\hat{w}^\ell \in \mathbb{R}_{+},\; \ell \in \pazocal{L}$ & Average waiting time for line $\ell$ \\
$\hat{t}^\ell \in \mathbb{R}_{+},\; \ell \in \pazocal{L}$ & Travel time from one terminal to the other on line~$\ell$ \\
\midrule
\multicolumn{2}{l}{\textbf{Induced quantities (passenger assignment outputs):}} \\
\midrule
$p_{ij}^{k\ell} \in \{0,1\}$ & 1 if OD pair $k$ traverses arc $(i,j)$ using line $\ell$ \\
$\hat{t}_k \in \mathbb{R}_{+},\; k \in \pazocal{K}$ & \rev{zg}{Travel time of OD pair k}{Predicted trip travel time of OD pair~$k$ in the redesigned network~$\pazocal{G}_\text{PT}^\text{adpt}$ (see \S{}\ref{sec:trip-travel-times})} \\
\rev{zg}{}{$\hat{t}^0_k \in \mathbb{R}_{+},\; k \in \pazocal{K}$} & \rev{zg}{}{Predicted trip travel time of OD pair~$k$ on the \emph{original} network~$\pazocal{G}_\text{PT}^\text{org}$ under predicted link times~$\hat{t}_{ij}$, i.e., the predicted cost of not redesigning (see \S{}\ref{sec:trip-travel-times})} \\
\bottomrule
\end{tabular}
\caption{Notation for PT Network Re-Design Problem}
\label{tab:notation}
\end{table}
% --------------------------------------------------------------
% Notation TableEND
% --------------------------------------------------------------

% ++++++++++++++++++++ 3.4
\subsection{Trip Travel Times}
\label{sec:trip-travel-times}

Let $\pazocal G^{\text{sub}}$ and $\pazocal G_{\text{PT}}$ denote the substrate and public transport (PT) graphs, respectively (\S\ref{sec:graphMODEL}). Each OD pair $k \in \pazocal K$ is defined by an origin stop $o_k \in \pazocal V^{\text{sub}}$ and a destination stop $d_k \in \pazocal V^{\text{sub}}$.
We assume origins and destinations are both stops, since we take responsibility of user trips from the moment they enter PT to the moment they exit, as considered in recent graph-based public transport modeling studies~\citep{zhao2024origin}.

\begin{observation}
\label{obs:det_routing}
Given a fixed PT network design, including line structures $(e_{ij}^\ell,a_i^\ell,b_i^\ell)$ and fleet allocations $n^\ell$, we evaluate the service offered to each OD pair through the minimum generalized travel time attainable on the resulting PT network. This shortest-path evaluation is consistent with standard practice in transit network design and line-planning optimization, where candidate networks are often assessed by assigning OD demand to shortest paths or by ensuring that passengers can be served on shortest routes, rather than by embedding a full behavioral route-choice model~\citep{Mumford2019}. The computed shortest path should therefore be interpreted as the best travel-time option made available by the design, not as an operational constraint imposed on individual passengers. Users may still choose other feasible routes, but such behavioral heterogeneity is outside the scope of the present optimization model. With additive link costs, the OD-specific benchmark paths are induced by the design and can be computed via standard shortest-path algorithms.
\end{observation}

% \rev{aa}{This assignment is deterministic and \emph{uncongested}: link travel times are exogenous to passenger volumes, so the reported trip times are user-optimal times that abstract away vehicle capacity and in-vehicle crowding (a redesigned line is never penalized for the demand it attracts). This is a standard simplification at the network-design level~\citep{goerigk2017line, schobel2022cheapest, john2014improved}, and the travel-time improvements reported in \S\ref{sec:numericalRESULT} should be read as uncongested user-optimal gains; capacity- and crowding-aware assignment is left to future work (\S\ref{sec:conclu}).
% }
% {
Accordingly, when computing shortest paths we assume links are uncongested: 
link travel times are exogenous to passenger volumes, so the reported trip times are user-optimal times that abstract away vehicle capacity and in-vehicle crowding (a redesigned line is never penalized for the demand it attracts). This is a standard simplification at the network-design level~\citep{goerigk2017line, schobel2022cheapest, john2014improved}, and the travel-time improvements reported in \S\ref{sec:numericalRESULT} should be read as uncongested user-optimal gains; capacity- and crowding-aware assignment is left to future work (\S\ref{sec:conclu}).
% }

For any network configuration $\pazocal G_\text{PT}$ and any realization~$\omega\in\Omega$, each link $(i,j)$ can be traversed in time $t_{ij}^\omega$. Based on such values, we compute shortest paths for any OD pair $k\in\pazocal K$ and denote the trip travel time as $t_k^\omega$. \rev{zg}{Note that~$t_k^\omega$ includes (i) in-vehicle travel times, (ii)~dwell times, (iii)~initial waiting times, and (iv)~transfer waiting times, according to \S{}\ref{sec:virtual_nodes}.}{Note that~$t_k^\omega$ includes (i) in-vehicle travel times, (ii)~dwell times, (iii)~initial waiting times, and (iv)~transfer waiting times plus the fixed transfer penalties, according to \S{}\ref{sec:virtual_nodes}.} Let $t_{ij}$ be the random variable with realizations~$t_{ij}^\omega,\omega\in\Omega$. 
%
%, the travel time of OD pair $k \in \pazocal K$, denoted by $t_k$, is defined as the sum of the link times along the shortest path between $o_k$ and $d_k$ on the virtual-node graph constructed in \S\ref{sec:virtual_nodes}.  Edge weights in this graph explicitly account for: (i) in-vehicle travel times, (ii) dwell times, (iii) initial waiting times, and (iv) transfer waiting times. 
%
To formalize these paths, we define binary dependent variables $p_{ij}^{k\ell,\omega}$ indicating whether arc $(i,j)$ on line $\ell$ belongs to the shortest path for OD pair $k$ in realization $\omega\in\Omega$. These variables are uniquely determined by the network state and computed algorithmically (\S\ref{sec:alg}) rather than being direct decision variables.

Let~$\pazocal G_\text{PT}^\text{org}$ denote the \emph{original PT network}, i.e., ``master'' schedule of PT over a typical day of operation. Classically, this is the network that will be operated every day. We instead allow adaptation, which means that the PT line network at a certain time is $\pazocal G_{PT}$, which may or may not correspond to the original one (Fig.~\ref{fig:rollFramework}).
In our framework, three notions of OD travel time are defined:
(i)~the actual travel time ~$t_k$, described in the paragraph above; we consider its realizations when evaluating the performance of our method.
(ii) the predicted travel time $\hat t_k$, computed as the sum of the predicted link times $\hat t_{ij}$ along the shortest paths described in the paragraph above; we use these predicted values as input parameters of optimization~\eqref{eq:multi_objective}-\eqref{st:veh_budget};
(iii) the baseline travel time $t_k^0$, defined as the OD travel time on the original network design $\pazocal G_\text{PT}^\text{org}$ under baseline link times $t_{ij}^0,\forall (i,j\in\pazocal E)$.
\rev{zg}{}{For completeness, we additionally denote by $\hat t_k^0$ the \emph{predicted} trip travel time of OD pair~$k$ on the original network $\pazocal G_\text{PT}^\text{org}$ under the predicted link times $\hat t_{ij}$; this represents the predicted cost of not redesigning, and is distinguished from $\hat t_k$, which is the predicted trip travel time on the redesigned network $\pazocal G_\text{PT}^\text{adpt}$.}
% % =============================================
% % 3.4.Problem statement
% % =============================================
\subsection{Problem Statement}\label{sec:problemStatement}
This research tackles the \rev{aa}{}{proactive re-}optimization of public transport lines (e.g., one hour ahead), \rev{aa}{redesigning them only}{triggered} when \rev{aa}{the forecast-trigger condition in Rule~\ref{rule:forecast_trigger} is met to adapt to fluctuating and unknown road conditions}{traffic conditions heavily deviate from the nominal ones, in order to preserve user and operator performance.}.
The proposed optimization method has three complementary objectives: 
(i)~minimizing user cost by reducing trip travel time deterioration due to fluctuations of substrate link times, with respect to the nominal ones (resulting in changes of actual travel time $t_k$ experienced on OD pair~$k$, relative to the baseline travel time experienced in the original network, which we denote as $t_k^0$), 
(ii) lowering operator costs represented by total vehicle-hours traveled by the fleet of buses, and (iii) avoiding excessive changes in the PT network compared to the original to improve operational manageability and prevent user disorientation.

Our optimization problem addresses the following three decision problems:
(a) how to dynamically modify PT lines (line topology, fleet allocation, headways) to adapt to actual road conditions when they severely deviate from the expected ones?
(b) which route should we make available for users to travel among OD pairs? and
(c) where to locate the terminals of each line?

\subsubsection{Timeframe and Stability Considerations}\label{sec:timeframe}
Reconfiguring PT lines may affect users and operators in distinct ways. 

For passengers, changes in routes can reduce perceived regularity and predictability, a recognized determinant of transit network design quality~\citep{vanoort2009regularity}. 
However, large urban public-transport systems already impose frequent on-the-fly journey adaptation; such changes are therefore not a completely new burden introduced by our approach, which at least attempts to prevent the degradation of user performance. In New York City Subway, major incidents (delaying at least 50 trains) occurred more than 500 times a year \citep[p.~2, ``Major Incidents''; p.~3, Fig.~3]{OSC2025SubwayDelays}; moreover, 50\% of surveyed subway users reported having been forced by delays to change mode \citep[p.~9]{NYCComptroller2017HumanCost}. Users are therefore already accustomed to experiencing public transport as an intrinsically variable system.
As a consequence, navigation apps are increasingly used~\citep {SmartphoneUse2016} and more and more passengers rely on them~\citep{durand2024digital}. \citet[p 147]{ghahramani2016trends} report that 80\% of New York PT passengers use those apps twice or more in one day. \cite{Shalaby2024} show that ``people seek out transit information the most during times of uncertainty'' to find their best itinerary.
Building on this trend, we assume that, under our proposed proactive approach, smartphone navigation apps provide passengers with the best itinerary, in case traffic variability is so high that a network redesign is triggered, thereby enabling them to navigate the reconfigured line structure.
% As for the passengers, we assume that smartphone navigation apps are available that provide passengers with the best itinerary on the current PT network.\footnote{
% This is in line with real-world observation, where navigation apps are increasingly used regularly \citep{SmartphoneUse2016} more and more passengers rely on them \citep{durand2024digital}. \cite[p 147]{ghahramani2016trends} report that 80\% of New York PT passengers use those apps two or more times in one day. \cite{Shalaby2024} show that `` people seek out transit information the most during times of uncertainty'' to find their best itinerary.
% } 
%
% \rev{aa}{While smartphone adoption for transit guidance is high and growing~}{}\citep{SmartphoneUse2016, durand2024digital}\rev{aa}{, we acknowledge that not all passengers have or regularly use smartphones. Complementary information channels, such as stop-level displays and real-time driver announcements, remain important for equitable and robust service provision; smartphone navigation is therefore treated as a predominant rather than universal mechanism.}{}
In addition to the availability of smartphone apps, behavioral trends of PT users suggest that they could be willing to adapt their trips to reconfigured network structures if this would allow them to arrive at their destination more efficiently.
Indeed, empirical evidence suggests that passengers operate a trade-off between efficiency and regularity~\citep{devarasetty2012value}. Specifically, some travelers may prefer an alternative path that ensures on-time arrival over a habitual route that, under adverse traffic conditions on a given day, results in delay. In this sense, deviations from regular routes are not necessarily detrimental if they yield tangible gains in travel time reliability. 

For the operator, frequently updated itineraries may increase cognitive load on the drivers.
Nonetheless, in practice, several agencies operate services with ``route deviation'', in which vehicles may depart from established routes within a defined area~\citep{NationalRTAP_ADA_Toolkit_2024}.  We assume that such operating areas have been defined ex ante for liability reasons.
The operational feasibility of these assumptions—one-hour advance notice, predefined deviation zones, and app-based passenger guidance—is therefore context-dependent. They will be more realistic in settings where flexible routing will be organizationally and technologically supported (e.g., with forthcoming vehicle automation~\citep{Peled2019}), while in more rigid operating environments the results presented here should be interpreted as upper bounds on the performance gains that could be achieved as operational practices, vehicle automation, and passenger-information systems evolve.
%How demanding these assumptions are (an hour of notice, predefined deviation zones, and app-based guidance for passengers) depends on the setting, and this aspect allow to distinguish where the method is sufficiently realistic today, where it fits more naturally, and where instead the results presented here should be intended as bounds that can be attained only in the future. With automated buses there is no driver to retrain and the vehicles already reroute freely~\citep{Peled2019}, so redesigning the network each hour fits well. For emergencies and planned events, the redesign covers a single corridor and is prepared beforehand, so the hour of notice is planning time rather than a real-time deadline, and approved versions can be readied in advance. 
\rev{aa}{}{The ways in which the results of our method would be used in practice could vary. One could, for instance, assume that operators could first revise the redesign issued by the proposed method, instead of allowing the line changes in an automatic manner. Or the method could be used as a decision support system for operational service adjustment, subject to the final judgment of the operators.}

%Moreover, passengers are more and more inclined towards flexible mobility services~\citep{van2023marta}.
Although PT line reconfiguration is implementable, PT lines are more easily manageable if disruptive reconfigurations are limited. For this reason, we measure structural stability using the Jaccard index~\citep{jaccard1901etude}, defined for the edge sets of the original and redesigned networks as
\begin{equation}
J(\pazocal E_{\text{orig}},\pazocal E_{\text{new}})
=\frac{|\pazocal E_{\text{orig}}\cap \pazocal E_{\text{new}}|}
       {|\pazocal E_{\text{orig}}\cup \pazocal E_{\text{new}}|}\in[0,1].
\label{eq:jaccard}
\end{equation}
It provides a simple, system-wide measure of overlap and has been shown to be particularly effective in sparse networks~\citep{travieso2024analytical}, which reflect typical PT topologies. Unlike sequence-based metrics such as edit distance~\citep{cheng2025computing}, which focus on local route alignment and scale poorly in multi-line systems~\citep{jain2024graph}, the Jaccard index offers a tractable, network-level notion of similarity that directly matches our design objectives.
We assume access to a short-term traffic prediction module, as is often the case in urban operations~\citep{li2017diffusion, liao2018deep}, and focus on the subsequent optimization step. \rev{zg}{As shown in Fig.~\ref{fig:rollFramework}, we consider an hourly monitoring policy: if average relative deviation $(1/|\pazocal E^\text{sub}|)  \sum_{(i,j)\in\pazocal E^\text{sub}}|\hat t_{ij}-t_{ij}^0|/t_{ij}^0$ remain within a tolerance (e.g., $0.1$), the network is kept unchanged or reset to the original configuration to avoid accumulated drift; larger deviations trigger a redesign. This selective updating strategy preserves temporal consistency while facilitating response to fluctuating traffic conditions. Moreover, deciding reconfigurations (when needed) always 1 hour in advance, gives the operator and drivers the needed time to adapt. In the hourly implementation, this one-hour advance-notice window is used as an offline planning period before the next operating interval. Because structural continuity is carried by the objective $Z_3$, we report the full trade-off between structural overlap and passenger benefit along the Pareto front rather than committing to a single operating point; the monitoring tolerance and the redesign trigger introduced above are the only thresholds fixed a priori.}{As shown in Fig.~\ref{fig:rollFramework}, we use the following hourly trigger.

\begin{deploymentrule}[Forecast-triggered line redesign]
\label{rule:forecast_trigger}
At the start of a time slot, define
$
\bar{\Delta}
\defeq
\frac{1}{|\pazocal E^{\text{sub}}|}
\sum_{(i,j)\in\pazocal E^{\text{sub}}}
\frac{|\hat t_{ij}-t_{ij}^{0}|}{t_{ij}^{0}}.
$
For a monitoring tolerance $\delta_{\mathrm{trig}}\in[0,1]$, the announced network is kept if
$\bar{\Delta}\le\delta_{\mathrm{trig}}$. If $\bar{\Delta}>\delta_{\mathrm{trig}}$, the operator solves
\eqref{eq:multi_objective}--\eqref{st:veh_budget} using the predicted link times $\hat t_{ij}$, and the selected solution becomes
$\pazocal G_{\mathrm{PT}}^{\mathrm{adpt}}$ for the next one-hour slot, subject to operator approval.
\end{deploymentrule}
The rule determines \emph{whether} a redesign is attempted. Conditional on triggering, the extent of structural change is controlled by $Z_3$ in~\eqref{eq:multi_objective}.}
%
% 如果每次重新设计的线路–边集合
%% 都与原始网络保持足够相似，那么任意两次相邻重设计之间的差异也必然有限（即相似度不会太低）。
The Jaccard-based structural term also provides a simple stability guarantee. If each redesigned public-transport network $\pazocal{G}_{\mathrm{PT}}^{(k)}$ remains sufficiently similar to the original configuration $\pazocal{G}_{\mathrm{PT}}^{(0)}$, then the similarity between any two redesigns is bounded from below, regardless of how many time slots separate them.
\begin{proposition}[\rev{zg}{}{Horizon-independent pairwise stability}]
\label{prop:temporal_stability}
Let $\pazocal{E}^{(k)}$ denote the set of active edges in the redesigned network at step $k$, \rev{zg}{and analogously for the original network and for the next redesign. Then, for any consecutive redesigns k and k+1}{and $\pazocal{E}^{(0)}$ for the original network. Then, for any two redesign steps $k$ and $k+n$ ($n \ge 1$)}, it holds that\footnote{The proof is provided in App.~\ref{app:temporal_stability_proof}.}
\[
J\!\bigl(\pazocal{E}^{(k)},\pazocal{E}^{(k+n)}\bigr) \;\ge\; J\!\bigl(\pazocal{E}^{(k)},\pazocal{E}^{(0)}\bigr) + J\!\bigl(\pazocal{E}^{(k+n)},\pazocal{E}^{(0)}\bigr) - 1.
\]
\end{proposition}
%% (Horizon-independence discussion moved to App.~\ref{app:temporal_stability_proof} to keep the main text concise.)
%
Rule~\ref{rule:forecast_trigger} specifies a possible deployment mechanism, whereas Proposition~\ref{prop:temporal_stability} gives an analytical stability property induced by the Jaccard-based structural term. Together, they \emph{motivate} a rolling, multi-period deployment of the proposed method. The experiments in this paper, however, evaluate the redesign step itself: namely, the quality of a single forecast-conditioned redesign under varying travel-time variability. A full empirical evaluation of the rolling trigger over realized multi-hour traffic time series is left to future work.

% =============================================
% 3.4.1 obj + s.t.
% =============================================
% ++++++++++++++++++++ 3.4.1
\subsubsection{Objective Functions}
\label{sec:obj}
We formulate our public transport network adaptation problem as a three-objective optimization, as expressed in~\eqref{eq:multi_objective}. The first objective, $Z_1$, minimizes passenger travel time deterioration, where $q_k$ represents the travel demand flow. 
The use of $\max(\hat{t}_k - t_k^0, 0)$ guarantees that only travel-time increases relative to the nominal benchmark are penalized, thereby avoiding situations where improvements for certain OD pairs compensate for deterioration in others. 
The second objective $Z_2$ minimizes operational costs through total vehicle-hours.
% \rev{zg}{}{ Here, operational cost is operationalized as the total vehicle-hours traveled by the fleet, a common proxy at the network-design level~\citep{arbex2015trb}. \rev{zg}{}{Concretely, $2\hat{t}^\ell$ is the round-trip running time of line $\ell$ and $n^\ell$ the number of vehicles assigned to it, so $Z_2=\sum_\ell n^\ell\,2\hat{t}^\ell$ aggregates fleet size and route length into a single vehicle-time figure per service cycle; it is a design-level proxy rather than a literal vehicle-hour count over the one-hour window (which would equal the active fleet multiplied by the window length), and we report its effect through the relative Operational Cost Reduction (\S\ref{sec:performMETRICS}).} It captures the dominant fleet-operation effort but excludes crew wages, deadheading, layover, and depot costs, which depend on vehicle scheduling and driver rostering and lie outside the present scope. We retain the term ``operational cost'' throughout, with this vehicle-hour proxy understood.}
%
The third objective $Z_3$ minimizes network structural changes using the Jaccard index~($\S$~\ref{sec:timeframe}), calculated as the ratio between the number of shared edges (i.e., those present both in the original and in the re-designed network) and the total number of edges (i.e., the union of the edges in the two networks).
Structural continuity is handled as an \emph{objective} ($Z_3$).

%
% Throughout the paper, the objective $Z_3$, the Route Overlap metric, and Proposition~\ref{prop:temporal_stability} are computed on the same set of \emph{line-labelled} edges, i.e.\ triplets $(\min(i,j),\max(i,j),\ell)$, so that the same physical link served by different lines is counted as distinct (\S\ref{sec:evaluation}).} 

\begin{equation}
\label{eq:multi_objective}
\min_{\pazocal{G}_{\text{PT}}^{\text{adpt}}} \begin{pmatrix} Z_1 \\ Z_2 \\ Z_3 \end{pmatrix} = \begin{pmatrix}
\underbrace{\sum_{k \in \pazocal{K}} q_k \cdot \max(\hat{t}_k - t_k^0, 0)}_{\text{Passenger travel time deterioration}} \\[8ex]
\underbrace{\sum_{\ell \in \pazocal{L}} n^\ell \cdot 2\hat{t}^\ell}_{\text{Total vehicle-hours (operating cost)}} \\[8ex]
\underbrace{1 - \frac{\sum_{\ell \in \pazocal{L}} \sum_{(i,j) \in \pazocal{E}} \min(e_{ij}^\ell, e_{ij}^{\ell, 0})}{\sum_{\ell \in \pazocal{L}} \sum_{(i,j) \in \pazocal{E}} \max(e_{ij}^\ell, e_{ij}^{\ell, 0})}}_{\text{Line dissimilarity}}
\end{pmatrix}
\end{equation}
In~\eqref{eq:multi_objective}, $e_{ij}^{\ell, 0}$ and $e_{ij}^\ell$ denote whether line $\ell$ uses substrate link $(i, j)$ in the original and the redesigned network, respectively. Observe that $e_{ij}^{\ell, 0}$ are parameters of our problem, whereas $e_{ij}^\ell$ are decision variables.
Observe that when solving our optimization problem, the realized link travel times $t_{ij}$ are not available, we will thus optimize the decision variables based on prediction $\hat{t}_{ij}$ and then evaluate, ex-post in the results, the objective function evaluated using the realized link travel times $t_{ij}$.

\rev{zg}{}{\noindent\textbf{Variable domains.} Binary: $y^\ell, e_{ij}^\ell, a_i^\ell, b_i^\ell \in \{0,1\}$. Integer: $n^\ell \in \mathbb{Z}_{\geq 1}$, $u_i^\ell \in \{1,\ldots,|\pazocal{N}|\}$. Continuous: $\hat{w}^\ell, \hat{t}^\ell \in \mathbb{R}_{\geq 0}$. The induced OD travel time $\hat{t}_k$ is a function of the design variables and the predicted link times~$\hat{t}_{ij}$; it is computed algorithmically (see \S{}\ref{sec:evaluation}) and is not a free variable of the program.}

% ++++++++++++++++++++ 4.2
\subsubsection{Constraints}
% --------------------------------------------------------------
% Constraints
% --------------------------------------------------------------
We now report the constraints that, together with the objective~\eqref{eq:multi_objective},
formalize the upper-level network design problem. Passenger route variables are treated as induced (follower) variables and are discussed separately.

% ======================= 1.时间约束
\textbf{Performance-Related Constraints}~\eqref{eq:headwayTime} and~\eqref{st:tl}.

% ======================= 2.线路约束
\textbf{Line Structure Constraints}

\begin{align}
%-------------------- 1. Activation & Terminal Control --------------------
% Ensure exactly one origin and one destination if line is active
& \sum_{i \in \pazocal{N}} a_i^\ell = y^\ell, \qquad
  \sum_{i \in \pazocal{N}} b_i^\ell = y^\ell
&& \forall \ell \in \pazocal{L} \label{st:term_existence}\\
% Forbid coincident terminals (Origin != Destination)
& a_i^\ell + b_i^\ell \le 1
&& \forall i \in \pazocal{N},\ \forall \ell \in \pazocal{L} \label{st:term_distinct}\\
%-------------------- 2. Topology & Flow Conservation --------------------
% Edges can only be selected for active lines (Both directions checked against y)
& e_{ij}^\ell + e_{ji}^\ell \le y^\ell
&& \forall \{i,j\} \in \pazocal{E},\ \forall \ell \in \pazocal{L} \label{st:edge_activation}\\
% Flow conservation: Inflow + IsOrigin = Outflow + IsDestination
& \sum_{j:\{i,j\} \in \pazocal{E}} e_{ji}^\ell + a_i^\ell
=
\sum_{j:\{i,j\} \in \pazocal{E}} e_{ij}^\ell + b_i^\ell
&& \forall i \in \pazocal{N},\ \forall \ell \in \pazocal{L} \label{st:flow_balance}\\
% Simple Path Constraint: Max Out-degree <= 1 AND Max In-degree <= 1
& \sum_{j:\{i,j\} \in \pazocal{E}} e_{ij}^\ell \le 1, \qquad
  \sum_{j:\{i,j\} \in \pazocal{E}} e_{ji}^\ell \le 1
&& \forall i \in \pazocal{N},\ \forall \ell \in \pazocal{L} \label{st:simple_path}\\
%-------------------- 3. Subtour Elimination (MTZ) --------------------
% MTZ constraint applied to both potential directions of an undirected edge
& u_i^\ell - u_j^\ell + 1
\le |\pazocal{N}| \big(1-e_{ij}^\ell\big)
&& \forall \{i,j\} \in \pazocal{E},\ \forall \ell \in \pazocal{L} \label{st:mtz_fwd}
\\
% MTZ Bounds
& 1 \le u_i^\ell \le |\pazocal{N}|
&& \forall i \in \pazocal{N},\ \forall \ell \in \pazocal{L} \label{st:mtz_bound}\\
%-------------------- 4. Resource Allocation --------------------
& n^\ell \ge y^\ell, \qquad
  n^\ell \le \bar n \, y^\ell
&& \forall \ell \in \pazocal{L} \label{st:veh_bounds}\\
& \sum_{\ell \in \pazocal{L}} n^\ell \le V_{\max}
&& \label{st:veh_budget}
\end{align}
Constraints~\eqref{st:term_existence}--\eqref{st:veh_budget} establish the topological structure and operational limits of the designed network. 
Constraints~\eqref{st:term_existence}--\eqref{st:term_distinct} synchronize line activation with terminal selection, mandating that every active line ($y^\ell=1$) connects exactly one initial terminal to one final terminal, thereby preventing zero-length loops. 
Constraint~\eqref{st:edge_activation} ensures that if a line is not activated ($y^\ell=0$), there are no links associated to it. Otherwise, it ensures that the line is constructed in one direction (links $(i,j)$ and $(j,i)$ cannot be both selected) and then buses go back and forth along it (as specified in page~\pageref{sentence:line-in-one-direction}).
%
%To model directed routes on the undirected substrate graph~$\pazocal{G}^{\text{sub}}$, 
% constraint~\eqref{st:edge_activation} binds arc selection to line activation and explicitly enforces that, for each physical link, at most one traversal direction can be selected (i.e., $e_{ij}^\ell + e_{ji}^\ell \le 1$). 
%
%constraint~\eqref{st:edge_activation} links arc selection to activation; when  $y^\ell=0$ it forbids any traversal, and when $y^\ell=1$ it ensures that each physical link is used in at most one direction.
%
Bus route continuity and linearity are enforced by the flow conservation in~\eqref{st:flow_balance} combined with the degree Constraints in~\eqref{st:simple_path}; the latter restrict nodes to a maximum of one predecessor and one successor, strictly imposing a non-branching simple path topology. 
Acyclic integrity is further guaranteed by Constraints~\eqref{st:mtz_fwd}--\eqref{st:mtz_bound}, which adapt the Miller-Tucker-Zemlin (MTZ) formulation to eliminate disconnected subtours that might otherwise satisfy flow balance: if link $(i,j)$ is used in line~$\ell$, then node $j$ should appear after $i$ along the line (its order of appearance ($u_j$) must be higher than $i$).
Finally, resource allocation is governed by Constraints~\eqref{st:veh_bounds}--\eqref{st:veh_budget}, which assign vehicles only to active lines and enforce the global fleet budget.
%
% 起终点控制：激活线路恰有一个起点和终点，未激活线路没有
% 边选择与激活绑定：只有激活线路才允许选择边
% 路径结构约束：每条线路在结构上是一条无分叉有向路径
% 无环性：通过 MTZ 排除任何形式的有向回路
% 车辆配置约束：只有激活线路才分配车辆，且车辆数满足运营可行性
% 全局车队预算

% ======================= 3. Path feasibility constraints (induced variables)
\textbf{Induced Passenger Paths.}
Passenger routes are not decision variables of the network design problem.
Given a fixed network design $(\mathbf e,\mathbf n)$, passenger paths for each OD pair $k\in\pazocal{K}$ are deterministically induced by solving a shortest-path problem on the auxiliary passenger graph $\pazocal{G}^{\text{PAX}}$ defined in \S\ref{sec:virtual_nodes}.
Let $\hat t_k$ denote the resulting shortest-path travel time under predicted link costs $\hat t_{ij}$.
For bookkeeping and interpretation, we also define $p_{ij}^{k\ell}\in\{0,1\}$ as an induced indicator equal to one if the shortest path for OD pair $k$ traverses arc $(i,j)$ on line~$\ell$, and zero otherwise.
Accordingly, $\hat t_k$ can be expressed as the sum of the costs of the traversed arcs (including in-vehicle, dwell, access, and transfer components) on $\pazocal{G}^{\text{PAX}}$.
These induced quantities are computed algorithmically and are not part of the optimization decision set.
% %
% \rev{zg}{}{This design-evaluation decomposition is related to, but should not be conflated with, a bilevel (Stackelberg) program. Passengers do optimize their own travel time: given a fixed network, each OD demand follows a least-cost path, so a lower-level routing optimization is indeed present. What separates our setting from a Stackelberg game is that this lower-level response is non-strategic and decoupled from passenger flow. Because link travel times are exogenous to the routing decision (link costs do not depend on the passenger volumes assigned to them, i.e.\ there is no congestion interaction among travellers), the lower-level problem decomposes into independent shortest-path problems with a unique, deterministic solution that is computed exactly in polynomial time, rather than a flow-dependent user equilibrium requiring a fixed-point (Wardrop-type) solution. The lower-level optimum can therefore be substituted in closed form into the upper-level evaluation, collapsing the structure to a single-level combinatorial program over the design variables $(\mathbf e,\mathbf n)$, whose feasibility and boundedness follow from the finite design space and the vehicle budget constraint~\eqref{st:veh_budget}. This is also what avoids the nested re-solves characteristic of strategic bilevel models.}

% =============================================
% 4.Algorithm
% =============================================
% ============================= FIGURE
\begin{figure}
\centering
\includegraphics[width=0.85\textwidth]{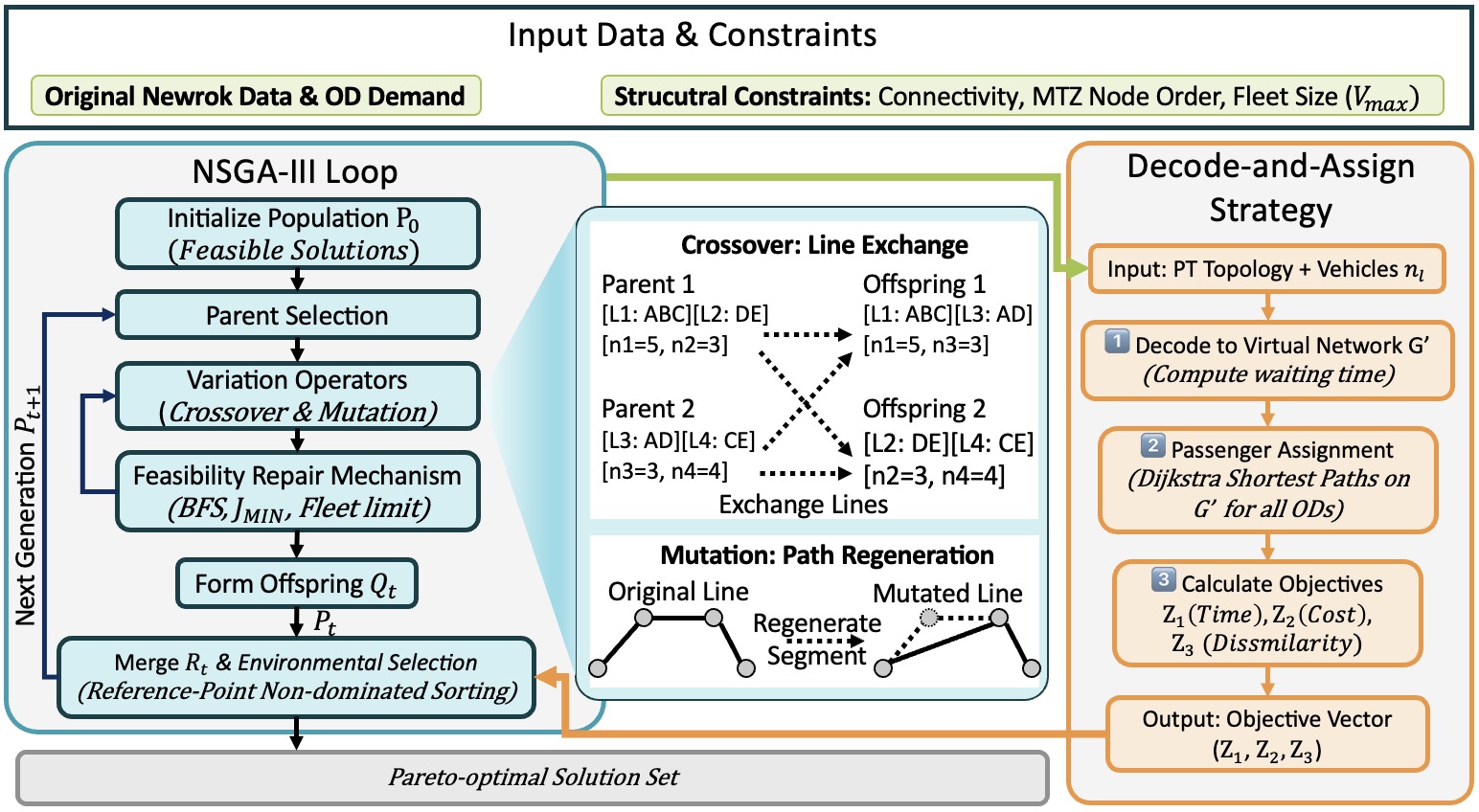}
\caption{\label{fig:nsgaFramework}Flowchart of the Proposed NSGA-III-based optimization method.}
\end{figure}
% ============================= FIGUREend

\section{Solution Algorithm}\label{sec:alg}
The optimization model presented in $\S$~\ref{sec:metho} is a mixed-integer formulation featuring discrete network design choices and nonlinear objective terms.
While the delay penalty in $Z_1$ is piecewise linear, the fractional Jaccard objective $Z_3$ introduces nonconvexity and global variable coupling.
Since an exact linearization of this ratio would require auxiliary constraints that severely compromise scalability, we retain the fractional structure.
Given the tri-objective nature of the problem and the nonconvex trade-offs, we adopt an NSGA-III framework (Algorithm~\ref{alg:NSGA3}) to approximate the Pareto front.
%
% As illustrated in Fig.~\ref{fig:nsgaFramework}, we decompose the decision process to ensure computational tractability.
% The \textbf{upper level} handles network configuration (topology $e_{ij}^\ell$ and fleet $n^\ell$), while the \textbf{lower level} treats passenger routing as a subproblem, where flow variables $p_{ij}^{k\ell}$ are deterministically derived via shortest-path calculations on the virtual graph $G'$ ($\S$\ref{sec:virtual_nodes}).
As illustrated in Fig.~\ref{fig:nsgaFramework}, we adopt a \textbf{design--evaluation decomposition} for computational tractability.
The evolutionary search operates on network design variables (line topologies $e_{ij}^\ell$ and fleet allocations $n^\ell$).
For any fixed design $(\mathbf e,\mathbf n)$, passenger routing is \emph{not} treated as an independent optimization model; instead, it is a deterministic \textbf{decode-and-assign mapping} induced by Observation~\ref{obs:det_routing}.
Specifically, OD travel times and the corresponding assignment outputs are computed by shortest-path calculations on the virtual graph $G'$ (\S\ref{sec:virtual_nodes}), which provides the objective values $(Z_1,Z_2,Z_3)$ used by NSGA-III.
Accordingly, the proposed workflow is structured into three coordinated modules:
(i) \textbf{Solution Encoding and Initialization} ($\S$~\ref{sec:encoding}) constructs a population of structurally feasible designs;
(ii) the \textbf{NSGA-III Evolutionary Loop} ($\S$~\ref{sec:search}) drives the feasible solution search, employing reference-point-based selection and variation \& repair operators to evolve the topology; and
(iii) the \textbf{Decode-and-Assign Strategy} (\S~\ref{sec:evaluation}) performs deterministic passenger assignment to compute the objective vector $(Z_1, Z_2, Z_3)$.

% --------------------------------------------------------- 
% --------------------------------------------------------- 
% ===================== algNSGA
\begin{algorithm}[htbp]
\small
\caption{NSGA-III based Adaptive PT Network Optimization}
\label{alg:NSGA3}
\begin{algorithmic}[1]
    \renewcommand{\algorithmicrequire}{\textbf{Input:}}
    \renewcommand{\algorithmicensure}{\textbf{Output:}}
    
    \REQUIRE $\pazocal{G}=(\pazocal{N}, \pazocal{E})$, $\pazocal{K}$, $\pazocal{L}$, $V_{\max}$, $q_k$, $\hat{t}_{ij}$, $\textit{PopSize}$, $\textit{MaxIter}$, $N_{\textit{ref}}$, $P_c$, $P_m$, $J_{\min}$
    \ENSURE Approximate Pareto set $\pazocal{P}^*$
    
    \STATE Initialize a population $\pazocal{P}$ of size $\textit{PopSize}$ via \S\ref{sec:encoding} and \S\ref{sec:init_population} (using parameter $J_\text{min}$) \label{line:initialization}
    \STATE Evaluate $(Z_1, Z_2, Z_3)$ for all $p \in \pazocal{P}$ via decode-and-assign (\S\ref{sec:evaluation}) \label{line:eval_init}
    
    \FOR{$iter = 1$ \TO $\textit{MaxIter}$}
        \STATE Select parents from $\pazocal{P}$ \label{line:parent_selection}
        \STATE Generate offspring $\pazocal{Q}$ by crossover (\S\ref{sec:crossover}) with $P_c$ and mutation (\S\ref{sec:mutation}) with $P_m$ \label{line:variation}
        
        % \STATE Repair infeasible lines using BFS path completion (Alg.~\ref{alg:bfs_repair})
        \STATE Apply feasibility-oriented repair to offspring when line-level structural infeasibility (e.g., disconnection or discontinuity) is detected (Alg.~\ref{alg:bfs_repair})
        \STATE Evaluate $(Z_1, Z_2, Z_3)$ for all $q \in \pazocal{Q}$ via decode-and-assign (\S\ref{sec:evaluation}) \label{line:evaluation}
        
        \STATE $\pazocal{R} \gets \pazocal{P} \cup \pazocal{Q}$
        \STATE Perform non-dominated sorting on $\pazocal{R}$ \label{line:nd_sort}
        \STATE Apply NSGA-III selection with $N_{\textit{ref}}$ reference points to obtain next $\pazocal{P}$ (\S\ref{sec:env_selection}) \label{line:env_select}
    \ENDFOR
    
    \STATE \textbf{return} the non-dominated solutions in $\pazocal{P}$ as $\pazocal{P}^*$
\end{algorithmic}
\end{algorithm}
% ===================== algNSGAend
% --------------------------------------------------------- 
% --------------------------------------------------------- 
% ------------------------------------------------------
% ============================ 4.1
% ------------------------------------------------------
\subsection{Solution Encoding}\label{sec:encoding}

% ********************* 线路组成
In Alg.~\ref{alg:NSGA3}, each PT line $\ell$ is represented as an ordered simple path
$p_\ell = (i_1, i_2, \dots, i_{n_\ell})$ with $i_j \in \pazocal{N}$ and no repeated nodes, which enforces acyclicity by construction according to constraints~\eqref{st:term_existence}--\eqref{st:veh_budget}.
PT lines may either be constructed from \textit{scratch} or obtained by \textit{reconstructing bus routes from the existing network}.
Starting from a current node $i$, the algorithm maintains a visited-node set $\pazocal{V}_{\text{visited}}^\ell$ and defines the set of feasible candidates as
$
\pazocal{N}_{\text{unvisited}}(i)
=
\{\, j \in \pazocal{N}(i) : j \notin \pazocal{V}_{\text{visited}}^\ell \,\},
$
which is used whenever a new path segment needs to be generated.
For lines constructed from scratch, the path is iteratively extended by selecting the next node $j \in \pazocal{N}_{\text{unvisited}}(i)$ according to the travel-time weighted probability
\begin{equation}
\label{eq:node-addition-probability}
\mathbb{P}(j \mid i)
=
\frac{1/t_{ij}}{\sum_{k\in\pazocal{N}_{\text{unvisited}}(i)} 1/t_{ik}},
\end{equation}
which favors temporally efficient connections while preserving stochastic exploration.
When a line is derived from the existing network, local modifications are applied by replacing a selected route segment with an alternative feasible path that preserves the original endpoints.

% ********************* 车辆分配
Simultaneously, vehicle allocation is managed by assigning integer variables $n^\ell$ to each line. 
During initialization (Alg.~\ref{alg:NSGA3}, line~\ref{line:initialization}),
each line is pre-allocated one vehicle to guarantee basic service availability.
The remaining $V_{\max}-|\pazocal{L}|$ vehicles are then distributed across
lines proportionally to their initial one-way travel times $t^\ell$ (computed based on nominal link travel times $t_{ij}$).
% and encoded as part of the initial chromosome representation,
% with integer rounding ensuring
% $\sum_{\ell\in\pazocal{L}} n^\ell = V_{\max}$.
% %
% Subsequently, fleet allocation is jointly evolved with line topologies through
% crossover and mutation operators, and waiting times $w^\ell$
% (defined in Eq.~\eqref{eq:headwayTime}) are updated whenever $n^\ell$ is modified.

% ------------------------------------------------------
% ============================ 4.2
% ------------------------------------------------------
\subsection{Evolutionary Search Procedure}
\label{sec:search}

Each individual encodes a complete PT network design, comprising a set of line path representations and the associated vehicle allocations.
Starting from an initial population, the algorithm evolves solutions in a generational manner.
At each generation, offspring are produced through selection, variation, and evaluation, followed by the NSGA-III reference-point-based \textit{environmental selection}, which prioritizes non-dominated solutions while maintaining diversity along the Pareto front.
The overall procedure is summarized in Alg.~\ref{alg:NSGA3}.

% -------------------- 4.2.1
\subsubsection{Population Initialization}
\label{sec:init_population}

During population initialization (Alg.~\ref{alg:NSGA3}, line~\ref{line:initialization}),
each individual is generated using the line construction rules described in $\S$~\ref{sec:encoding}.
PT lines are either constructed from scratch or obtained by locally perturbing routes reconstructed from the existing network.
Given a set of lines, vehicle allocations are initialized by assigning one vehicle to each active line to guarantee basic service availability.
The remaining fleet is then distributed across lines according to the chromosome encoding scheme described in $\S$~\ref{sec:encoding}, ensuring $\sum_{\ell \in \pazocal{L}} n^\ell = V_{\max}$.

To expedite initialization, candidate solutions are screened via a \textbf{connectivity check} on the undirected graph induced by
$\pazocal E_{\text{PT}} = \bigcup_{\ell\in\pazocal L} \pazocal E_\ell$,
i.e., the union of all active line edges in $\pazocal G_{\text{PT}}$.
By employing a Union--Find data structure, we efficiently identify connected components, and verify that, for every OD pair $k\in\pazocal K$, the origin $o_k$ and destination $d_k$
belong to the same connected component (which implies that it is possible to go from one to the other).
% This screening step avoids the computational overhead of full shortest-path evaluations during population initialization.
%
The procedure is repeated until $\textit{PopSize}$ individuals (i.e., PT network structures) are generated that satisfy:
(i)~structural validity (lines are simple-paths, with no branches nor cycles),
(ii)~connectivity-based serviceability, i.e., all OD pairs $k\in\pazocal K$ are connected within $\pazocal G_{\text{PT}}$,
and (iii)~similarity, i.e., individual~$\pazocal G_\text{PT}$ is not extremely different from original PT network~$\pazocal G_\text{PT}^\text{org}$, i.e., the Jaccard similarity index \eqref{eq:jaccard} must be above threshold $J_{\min}$.
These criteria enforce a necessary (but not sufficient) condition for OD reachability; full feasibility and network performance are subsequently evaluated using the decode-and-assign objective.

% -------------------- 4.2.2
\subsubsection{Crossover Operation}\label{sec:crossover}

With probability $P_c$, a crossover operator is applied to two parent solutions, denoted as parents $A$ and $B$.
Each solution consists of a fixed set of PT lines, where each line is treated as a basic inheritance unit characterized jointly by its bus route topology and vehicle allocation.
The crossover exchanges complete lines between parents as the primary recombination mechanism; any subsequent modifications are limited to feasibility restoration.

Line inheritance is governed by a binary selection mask $\boldsymbol{\xi}=(\xi_\ell)_{\ell\in\pazocal{L}}$, where each component is independently drawn as $\xi_\ell \sim \mathrm{Bernoulli}(0.5)$,
so that both parents are treated symmetrically and structural diversity is promoted at the line level.
Given $\boldsymbol{\xi}$, the offspring design is defined by block-wise inheritance:
\begin{align}
(e_{ij}^\ell)^{\text{off}}
&=
\xi_\ell (e_{ij}^\ell)^{A} + (1-\xi_\ell)(e_{ij}^\ell)^{B},
&& \forall \ell\in\pazocal{L},\ (i,j)\in\pazocal{E},
\label{eq:crossover_edges}\\
(n^\ell)^{\text{off}}
&=
\xi_\ell (n^\ell)^{A} + (1-\xi_\ell)(n^\ell)^{B},
&& \forall \ell\in\pazocal{L}.
\label{eq:crossover_fleet}
\end{align}
Whenever a line topology is inherited from a parent, its vehicle allocation is inherited jointly, thereby preserving the coupling between bus route geometry and service frequency.
% No interpolation across different line structures is performed.

% If any inherited line becomes structurally infeasible (e.g., disconnected or discontinuous), a localized repair mechanism is triggered (Alg.~\ref{alg:bfs_repair}).
If any inherited line becomes structurally infeasible (e.g., disconnected or discontinuous), a feasibility-oriented restoration operator is applied to recover a valid simple-path representation (Alg.~\ref{alg:bfs_repair}).
% Repairs are applied only to the affected lines and only when infeasibility is detected, by reconnecting disrupted segments using deterministic or randomized Breadth-First Search (BFS) on the substrate graph.
% This repair step is strictly feasibility-oriented and does not optimize route cost or passenger travel time.
%
Finally, fleet constraints ($n^\ell \ge 1$ for all active lines and $\sum_{\ell\in\pazocal{L}} n^\ell = V_{\max}$) are enforced via minimal $\pm 1$ adjustments on a small subset of lines.
This correction is a feasibility projection and is not guided by objective improvement.
% -------------------- 4.2.3
\subsubsection{Mutation Operation}\label{sec:mutation}

With probability $P_m$, a mutation perturbs a single line $\ell$ via local structural modifications, including node insertion, node deletion, or segment reconnection.
These operations introduce stochastic variations in line topology while preserving the overall chromosome structure.
If the mutation results in a disconnected or discontinuous line,
feasibility is restored by reconnecting the affected segments (Alg.~\ref{alg:bfs_repair}).
This repair step is feasibility-oriented and does not explicitly consider objective values.
Finally, global fleet conservation ($\sum_{\ell\in\pazocal{L}} n^\ell = V_{\max}$)
is enforced in the same manner as described in \S~\ref{sec:crossover}.

% -------------------- 4.2.4
\subsubsection{Environmental Selection}\label{sec:env_selection}
After offspring generation, parent and offspring populations are merged and ranked via non-dominated sorting.
Environmental selection follows the standard NSGA-III mechanism, employing a set of reference points $N_\text{ref}$ to maintain population diversity~\citep{deb2014nsga3}.\footnote{The implementation details are provided in App.~\ref{A.Evolution}.}
% ------------------------------------------------------
% ============================ 4.3
% ------------------------------------------------------
\subsection{Objective Evaluation via Decode-and-Assign}
\label{sec:evaluation}
% [4.3 内容]
% 关键部分！
% Step 1: Construct Virtual Graph from Decision Variables
% Step 2: Passenger Assignment (Dijkstra) -> "Induced Flows"
% Step 3: Compute Z1, Z2, Z3
The evaluation phase (Alg.~\ref{alg:NSGA3}, line~\ref{line:evaluation}) maps each chromosome into the objective space through a systematic decode-and-assign procedure. Following Observation~\ref{obs:det_routing}, we treat passenger routing as a deterministic response to the network configuration, thereby decoupling routing behavior from the upper-level network design decisions. 
For each individual, the algorithm first constructs a comprehensive virtual-node graph $\pazocal G'$ that integrates in-vehicle travel times and waiting times $w^\ell$, the latter being derived from the allocated vehicle units $n^\ell$ ($\S$~\ref{sec:virtual_nodes}).
To evaluate the travel-time deterioration $Z_1$, the algorithm solves a shortest-path problem (Dijkstra) for every OD pair $(o_k, d_k) \in \pazocal K$ on $\pazocal G'$. The resulting shortest path yields the predicted OD travel time $\hat t_k$. Feasibility is enforced by construction: $\pazocal G'$ includes only arcs corresponding to active line segments in the decoded network design, and OD paths are initialized through access/egress connections between physical stops and their associated line-specific virtual nodes (\S\ref{sec:virtual_nodes}).
This shortest-path evaluation is performed for each candidate network, consistently reflecting the current decoded topology and headways.
Simultaneously, the operational cost $Z_2$ and network dissimilarity $Z_3$ are computed. Specifically, $Z_3$ is quantified using a line-specific Jaccard similarity, where each network element is defined as a triplet $(\min(i,j), \max(i,j), \ell)$. This representation ensures that identical physical links served by different lines are treated as distinct entities, accurately capturing the topological shifts in the service network.

\section{Evaluation Method}\label{sec:evaMETHO}
In order to test if our \rev{zg}{real-time}{forecast-triggered} PT re-design method would perform well in a realistic situation with realistic uncertain fluctuations of link times, we need to tackle the following two challenges: (i)~link travel times are inherently stochastic, and (ii)~re-design decisions are based on travel time predictions, which might deviate from the true ones. To account for~(i), we generate stochastic travel times for all links, with different levels of variability ($\S$~\ref{sec:log_normal}). \rev{zg}{To account for~(ii), we then generate, in $\S$~\ref{sec:KDE}, prediction errors $\epsilon_{ij}, \forall (i,j) \in \pazocal{E}^{\text{sub}}$, statistically representative of the level of inaccuracy of well established prediction models (we consider Diffusion Convolutional Recurrent Neural Network (DCRNN) - see \S{}\ref{sec:KDE}).}{To account for~(ii), we generate empirical log-relative prediction errors $\eta_{ij}, \forall (i,j) \in \pazocal{E}^{\text{sub}}$, statistically representative of the inaccuracy of a well-established prediction model (we use the Diffusion Convolutional Recurrent Neural Network (DCRNN); see \S{}\ref{sec:KDE} and App.~\ref{A3}).}

% ------------------------------------------------------
% ++++++++++++++++++++ 5.2 讲解cv
% ------------------------------------------------------
\subsection{Model of Travel Time}
\label{sec:log_normal}

To account for the variability in travel times, we model them as random variables, according to $\S$~\ref{sec:Stochastic-Modeling-of-Travel-Times}. The shape we choose for this random variable is a log-normal distribution, which effectively captures the skewed nature of travel time data~\citep{guessous2014estimating} and accounts for extreme delay scenarios that can occur in reality, although not often.
On each link $(i,j)$, we generate travel time realization~$t_{ij}^\omega$ by performing Monte-Carlo sampling from random variable $t_{ij} \sim \text{log-normal}(\mu_{ij}, \sigma_{ij}^2)$ such that the median of $t_{ij}$ corresponds to the nominal time reported in the Mandl's dataset, and the coefficient of variation ($cv$) matches a predetermined value: we test $cv$ = 0.05, 0.25, 0.5, 1.0.
The parameterization used in this sampling is stated explicitly below.

\begin{proposition}[Log-normal parameterization for prescribed median and CV]
\label{prop:lognormal}
Let $t_{ij}$ be a lognormal random variable representing the travel time on link $(i,j)\in\pazocal{E}$ with median $m_{ij}$ and coefficient of variation $cv$. The parameters $(\mu_{ij},\sigma_{ij})$ of the underlying normal distribution satisfy
\[
    \mu_{ij} = \ln(m_{ij}), 
    \qquad
    \sigma_{ij} = \sqrt{\ln(1 + cv^{2})}.
\]
\end{proposition}

To test the performance of our optimization when faced with different levels of hypothetical traffic variability, for each $cv$, we generate a set of samples $\Omega_{cv}$, where each sample $\omega \in \Omega_{cv}$ corresponds to a collection of travel times $\{t_{ij}^{cv,\omega} \mid (i,j) \in \pazocal{E}^{\text{sub}}\}$ across the network. The proof of Proposition~\ref{prop:lognormal} is provided in App.~\ref{Appendix:log_normal}.
% ------------------------------------------------------
% ++++++++++++++++++++ 5.3 讲解noise
% ------------------------------------------------------
\subsection{Travel Time Prediction \rev{zg}{with KDE-Based Error Modelling}{via Empirical Log-Relative Error Sampling}}
\label{sec:KDE}

%\rev{zg}{Recall that our re-design method is based on predicted travel times (the input parameters of the optimization problem in $\S$~\ref{sec:problemStatement} are predicted link travel times~$\hat t_{ij}$ rather than the true values~$t_{ij}$, which are unknown in practice). We therefore need to evaluate whether the resulting decisions remain satisfactory despite the inevitable prediction errors~$\epsilon_{ij} \defeq \hat t_{ij}-t_{ij}$.
%}{
Our optimization uses predicted link travel times~$\hat t_{ij}$, while actual link travel times~$t_{ij}$ are unknown at decision time. We therefore evaluate whether the generated line plans remain effective even in presence of inevitable prediction errors.
Compliant with most traffic forecast models, we predict speeds on links, rather than directly predicting travel times. We express forecast errors as dimensionless quantities, valid both for travel times and speeds. This is why, in the following, we will just refer to ``errors''. In particular, let the \emph{logarithmic prediction errors} be 
\begin{equation}
\eta_{ij}\defeq 
\log(\hat t_{ij})-\log(t_{ij})=
\log(\hat t_{ij}/t_{ij})
=
\log\left( (d_{ij}/\hat v_{ij})/(d_{ij}/v_{ij}) \right)
=
\log( v_{ij}/ \hat v_{ij})
=
\log(v_{ij})-\log(\hat v_{ij})
\end{equation}
where $v_{ij}$ and $\hat v_{ij}$ are the true and predicted speed on link~$ij$, respectively, and $d_{ij}$ is the length of link~$ij$.

We can express the corresponding additive error (as defined in~\eqref{eq:link_time_error}) as $\epsilon_{ij}=\hat t_{ij}-t_{ij}=t_{ij}(\exp(\eta_{ij})-1)$.
%, but the sampling model is defined on $\eta_{ij}$ rather than on absolute errors.
%
%Let the \emph{empirical log-relative prediction residuals} be $\eta_{ij}\defeq \log(\hat t_{ij})-\log(t_{ij})=\log(\hat t_{ij}/t_{ij})$. We can express the corresponding additive error (as defined in~\eqref{eq:link_time_error}) as $\epsilon_{ij}=\hat t_{ij}-t_{ij}=t_{ij}(\exp(\eta_{ij})-1)$, but the sampling model is defined on $\eta_{ij}$ rather than on absolute errors.
%}

%\rev{zg}{Recall from $\S$~\ref{sec:Stochastic-Modeling-of-Travel-Times} that prediction errors are modeled as random variables. We generate them through a three-stage procedure:}{We model prediction uncertainty through random empirical residuals. We generate them through a three-stage procedure:}
% \todoi{aa: either always residual or always error}
We generate samples of prediction errors through a three-stage procedure:
\begin{itemize}
\item In \emph{Stage 1} (see $\S$~\ref{sec:errorStage1}), we train a prediction model to forecast traffic speeds.
\item In \emph{Stage 2} (see $\S$~\ref{sec:errorStage1}), we use the trained model to produce one-hour-ahead speed forecasts and collect empirical logarithmic errors $\eta_{ij}$.
%log-relative travel-time residuals.
\item \rev{zg}{
In \emph{Stage 3} (see $\S$~\ref{sec:errorStage3}), we estimate the error probability distribution using the collected error samples via Kernel Density Estimation (KDE).
}{
In \emph{Stage 3} (see $\S$~\ref{sec:errorStage3}), we sample with replacement from the set of empirical logarithmic errors and use those samples to construct noisy travel-time forecasts.}
\end{itemize}

\rev{zg}{As for the prediction model in Stage 1, we adopted DCRNN as a well-established spatio-temporal model; the goal was to obtain realistic prediction errors rather than compare forecasting models, and the adoption of DCRNN was presented as a conservative evaluation.}{We use the Diffusion Convolutional Recurrent Neural Network (DCRNN)~\citep{li2017diffusion} as a representative spatio-temporal predictor to generate empirical one-hour-ahead prediction errors. The objective is not to benchmark forecasting models, but to expose the line-plan generation method to realistic forecast errors. The dataset-specific training and error-extraction procedure is summarized below and detailed in App.~\ref{A2}.}
\rev{zg}{}{This prediction error-based construction is useful because downstream decisions may be sensitive not only to the average magnitude of forecast errors, but also to their distributional shape \citep{peled2021quality}.}

% =============================== 5.3.1

\subsubsection{Stage 1 and 2: Prediction model and empirical prediction error extraction}\label{sec:errorStage1}

We train a Diffusion Convolutional Recurrent Neural Network (DCRNN). 
Due to the absence of historical traffic data for Mandl’s Swiss network, the DCRNN is trained on the METR-LA dataset~\citep{li2017diffusion}, which contains traffic speed measurements from 207 sensors in Los Angeles County collected at 5-minute intervals over four months. Although METR-LA is a different network from Mandl's, we apply the error patterns derived via DCRNN in the METR-LA to our Mandl's network, as a practical way to reproduce the realistic errors of predictive models. This empirical setup does not affect the methodological validity of the presented approach, which focuses on optimization under uncertainty rather than on predicting models' accuracy.
For the real-world Beijing case study, instead, we utilize the Q-Traffic dataset~\citep{liao2018deep}, which provides traffic speed and road network data for the 6th Ring Road in Beijing, China. 

\rev{zg}{Post-training, we collect the 1-hour-ahead speed predictions in each link $(i,j)\in\pazocal E$ and the corresponding errors $\epsilon_{ij}$.}{Post-training, we collect the 1-hour-ahead speed predictions in each link $(i,j)\in\pazocal E$ and collect the resulting errors, as detailed in in App.~\ref{A2}.}

%derive prediction errors by evaluating 1-hour-ahead speed predictions at each sensor location.  These errors form the basis for modeling prediction uncertainty through additive link-level terms; specifically, for each link $(i,j) \in \pazocal{E}^{\text{sub}}$, the error $\epsilon_{ij}$ is defined as in~\eqref{eq:link_time_error}.  Since the prediction model operates on traffic speeds, we convert the predicted values into link travel times prior to optimization.  

% =============================== 5.3.2
\subsubsection{\rev{zg}{Stage 3: KDE modeling of prediction errors}{Stage 3: Empirical log-relative prediction error sampling}}
\label{sec:errorStage3}

% 人为设定一组 prediction noise distributions，然后用这些分布扰动真实 OD demand，生成 noisy demand predictions。

\rev{zg}{Following the previous framework, we modelled travel-time prediction errors as stochastic inputs to the network optimization problem. Earlier text contrasted pre-specified parametric noise distributions with a data-driven, non-parametric KDE approach.}{Following~\cite{peled2021quality}, we treat prediction errors as stochastic inputs to the network optimization problem. Rather than prescribing a parametric distribution, or fitting a KDE to absolute prediction errors, we collect a set of empirical prediction errors extracted from DCRNN test-set predictions. Sampling these errors directly preserves the trained predictor's empirical bias, skewness, and tail behavior without imposing a Gaussian, Weibull, or other smooth distributional shape. The error construction and sampling procedure are detailed in App.~\ref{A3}.}
\rev{zg}{In particular, instead of assuming a Gaussian or Weibull prediction error shape, we obtained prediction errors from an actual predicting model and fitted their distribution using KDE; the appendix detailed the implementation and reported that KDE performed best.}{The prediction errors are not recentred, so any mild empirical bias from the trained model is retained. Errors are not differentiated based on the links where they were calculated; consequently, we do not represent link-specific heterogeneity or spatio-temporal correlation.}

\rev{zg}{}{
To assess the sensitivity of our method to the magnitude of the prediction errors, in our optimization we set predicted travel times for a certain realization~$\omega\in\Omega$ as
%The additive prediction error is induced by scaling the sampled log-relative residual:
\begin{equation}
    \hat{t}_{ij}^{\omega}(\alpha)
    =
    t_{ij}^{\omega}\exp\!\big(\alpha\,\eta_{ij}^{\omega}\big),
    \label{eq:predTRAVEL_alpha}
\end{equation}
where $\eta_{ij}^{\omega}$ is sampled from the set of empirical logarithmic prediction errors. Equivalently, $\epsilon_{ij}^{\omega}=t_{ij}^{\omega}(\exp(\alpha\eta_{ij}^{\omega})-1)$. The parameter $\alpha$ scales the empirical forecast-error magnitude: $\alpha=1$ corresponds to the prediction error level extracted from the trained DCRNN, $\alpha=0$ gives the perfect-foresight benchmark, and $\alpha>1$ is used only in the forecast-sensitivity stress tests.}

% =============================== 5.3.2b  bias + correlation
% \rev{zg}{}{The residuals are not recentred, preserving any mild empirical bias from the trained model. They are sampled independently from a pooled link-level empirical distribution, so spatial heterogeneity and spatio-temporal correlation are not modelled; see App.~\ref{app:forecast_error_caveats} for details.}

% =============================== 5.3.3
\subsubsection{Simulating predicted travel times}

\rev{zg}{
For each simulation realization $\omega \in \Omega$ and link $(i,j) \in \pazocal{E}^{\text{sub}}$, we implement a two-step process: (1) sample actual travel time $t_{ij}^{\omega}$ from the log-normal distribution of $\S$~\ref{sec:log_normal}, and (2) sample a log-relative prediction error $\eta_{ij}^{\omega}$ via the empirical error sampling procedure detailed in App.~\ref{A3}. The resulting predicted travel time is:
}{
For each simulation realization $\omega \in \Omega$ and link $(i,j) \in \pazocal{E}^{\text{sub}}$, we implement a two-step process: (1) sample actual travel time $t_{ij}^{\omega}$ from the log-normal distribution of $\S$~\ref{sec:log_normal}, and (2) generate the additive prediction error $\epsilon_{ij}^{\omega}$ via the empirical log-relative error sampling procedure detailed in \S{}\ref{sec:errorStage3}. The resulting predicted travel time is:}
\begin{equation}
    \hat{t}_{ij}^{\omega}
    =
    \underbrace{t_{ij}^{\omega}}_{\text{$\S$~\ref{sec:log_normal}}}
    +
    \underbrace{\epsilon_{ij}^{\omega}}_{\text{$\S$~\ref{sec:KDE}}},
    \label{eq:predTRAVEL}
\end{equation}
\rev{zg}{}{In the main experiments we set $\alpha=1$ in Eq.~\eqref{eq:predTRAVEL_alpha}; alternative values are used only in the forecast-error sensitivity analysis~(App.~\ref{app:forecastValue}).}

% =============================================
% 6.Numerical Results
% =============================================

\section{Numerical Results}\label{sec:numericalRESULT}
\rev{zg}{We evaluate the proposed adaptive public transport (PT) redesign framework on both benchmark and real-world networks.}{We evaluate the proposed framework on Mandl's benchmark network (\S\ref{sec:mandl_RESULT1}) and a real-world Beijing network (\S\ref{sec:caseStudy}). In all experiments, the optimizer receives predicted link times~$\hat{t}_{ij}$ and produces a redesigned network; performance is then scored against the realized travel times~$t_{ij}^{\omega}$, so that reported figures reflect decision quality under imperfect forecasts.}
\rev{aa}{Each variability level~$cv$ represents the forecast-indicated stress of}{The following results concern} a \emph{single} decision slot; accordingly, the experiments assess the single-slot redesign rather than the multi-period rolling trigger of Rule.~\ref{rule:forecast_trigger}.

% ------------------------------------------------------ 
% % ++++++++++++++++++++ 6.0 - Considered scenarios 
% ------------------------------------------------------
\subsection{Experimental setup and metrics}\label{sec:setup}
Mandl’s Swiss network, initially used by~\cite{mandl1980evaluation}, consisting of 15 nodes, 21 links, and a symmetric demand of 15{,}570 passenger trips, is employed as a canonical benchmark for comparison with existing PT network design methods~\citep{liang2025novel, ahern2022approximate}. To further assess the method’s applicability under realistic operating conditions, we consider a real-world urban public transport network in a region in Beijing, China, recently analyzed by~\cite{wu2024multi}, which comprises 301 bus stops, 3,758 road sections, and 29 bus lines. Its origin--destination demand matrix is constructed from smart card data collected during peak hours by aggregating trips between boarding station $i$ and alighting station $j$. 

Key parameters are set as follows: a fixed dwell time of $d=1$ min per stop is assumed following~\cite{chen2022TRSC}.
For Mandl’s network, the maximum fleet size ($V_{\max}$) is capped at 40, following the service-frequency setting of~\cite{ahern2022approximate} to maintain comparable headways. 
For the Beijing network, $V_{\max}$ is set to 150 to accommodate a constant 12-minute headway across this large PT system. 

\subsubsection{Compared networks and baselines}
\label{sec:compared-PT-networks}
% \todoi{aa: Explain here how you pick a single design from the Pareto front.}

We wish to show the benefits of adapting PT network design \rev{zg}{in real-time}{for the hour ahead, and only when the forecast-trigger condition in Rule~\ref{rule:forecast_trigger} is met,} to varying road conditions via our method, \rev{zg}{against current network design approaches, which all keep the PT structure unchanged}{against selected static baselines, which keep the original line structure fixed}.

For \emph{Mandl’s Swiss network}, the current network design approaches are represented by \rev{zg}{two networks}{the following comparison networks}:
\begin{itemize}
\item A static integrated line-frequency design $\pazocal{G}_{\text{PT}}^{\text{ahe}}$, following~\citet{ahern2022approximate}, which jointly optimizes passenger travel time and operational cost.
\item A static state-of-the-art network design for trip travel-time minimization, $\pazocal{G}_{\text{PT}}^{\text{ver}}$, derived from the exact formulation of~\cite{vermeir2021exact}.
\item \rev{zg}{}{A frequency-only adjustment baseline, denoted by~$\pazocal{G}_{\text{PT}}^{\text{freq}}$, in which line structures are fixed and only service frequencies are adjusted.}
\end{itemize}

For the \emph{real-world Beijing network}, instead, the current design approaches are represented by
\begin{itemize}
\item \rev{zg}{The original bus line network, denoted by~$\pazocal{G}_{\text{PT}}^{\text{orig}}$, serving as the current operational standard against which the adaptive redesign is assessed.}{The original bus line network, denoted by~$\pazocal{G}_{\text{PT}}^{\text{orig}}$, serving as the current operational standard against which the adaptive redesign is assessed.}
\item \rev{zg}{}{A stop-skipping baseline, denoted by~$\pazocal{G}_{\text{PT}}^{\text{skip}}$, in which buses skip major congested \rev{zg}{stops}{links} while the underlying line structure is otherwise retained.}
\end{itemize}
\rev{zg}{Across all experiments, the baseline travel time for an OD pair $k \in \pazocal{K}$, denoted by $t_k^0$ (\S{}\ref{sec:trip-travel-times}), is computed as the shortest-path travel time between origin $o_k$ and destination $d_k$ in either network $\pazocal{G}_{\text{PT}}^{\text{ahe}}$, $\pazocal{G}_{\text{PT}}^{\text{ver}}$, or $\pazocal{G}_{\text{PT}}^{\text{orig}}$, depending on which one will be used as reference for comparison.}{Across all experiments, the baseline travel time for an OD pair $k \in \pazocal{K}$, denoted by $t_k^0$ (\S{}\ref{sec:trip-travel-times}), is computed as the shortest-path travel time between origin $o_k$ and destination $d_k$ in the relevant comparison network: $\pazocal{G}_{\text{PT}}^{\text{ahe}}$, $\pazocal{G}_{\text{PT}}^{\text{ver}}$, $\pazocal{G}_{\text{PT}}^{\text{orig}}$, $\pazocal{G}_{\text{PT}}^{\text{freq}}$, or $\pazocal{G}_{\text{PT}}^{\text{skip}}$.}
%
%strictly on the respective comparison PT network (i.e., $\pazocal{G}_{\text{PT}}^{\text{ahe}}$, $\pazocal{G}_{\text{PT}}^{\text{ver}}$, or $\pazocal{G}_{\text{PT}}^{\text{orig}}$).
\rev{aa}{}{Note that, overall, we use two dynamic benchmarks, i.e., (i)~frequency-only adjustment and (ii)~stop-skipping. We only keep in each scenario the one that performed the best (namely (i)~in Mandl and (ii)~in Beijing).}

To select representative compromise solutions under fluctuating travel time, we compute the solutions on the Pareto front via the proposed Alg.~\ref{alg:NSGA3}. 
Then, we compute the Euclidean distance from each of those solutions to an ideal point representing the median observed values in the Pareto front, across the three objectives. If not otherwise specified, the considered network design $\pazocal{G}_{\text{PT}}^{\text{adpt}}$ that we analyze in the following is the one with the objectives closest to that ideal point. 

All experiments are conducted on a server equipped with two Intel Xeon Gold 5220R CPUs (48 cores, 96 threads, up to 4.00 GHz) with a maximum computational time of 1h per instance, which we use here as an upper bound for our offline computational experiments within the hourly decision cycle (see Fig.~\ref{fig:rollFramework} and $\S$\ref{sec:timeframe}).\footnote{The values of the parameters of the proposed Alg.~\ref{alg:NSGA3}, as well as a thorough performance comparison against other multiobjective optimization methods, are extensively discussed in App.~\ref{sec:algorithmic-framework-and-evaluation}.}
%
% \footnote{
% The full implementation is publicly available at
% \rev{zg}{
\textit{(To respect double-blind review requirements, all problem instances and source code will be released in a public repository upon acceptance.)}
% }{\textit{Data and code availability: the problem instances and source code will be released in this public repository upon acceptance: \url{https://github.com/zihao-guo/disruption-model}.}}
% }

% ------------------------------------------------------
% ++++++++++++++++++++++++++++ 6.2. 旅行时间分析 (不同seed)
% ------------------------------------------------------
\subsubsection{Scenarios and evaluation protocol}
%To facilitate case-study interpretation, we identified the three most balanced solutions by measuring the Euclidean distance from each non-dominated solution to a central reference point. This reference point is defined as the component-wise median objective vector of all non-dominated solutions obtained across all algorithms and random seeds under a given uncertainty level. While not an optimal solution itself, this point serves as a measure of central tendency, enabling the selection of solutions that avoid extreme trade-offs between objectives.

We assess the effectiveness of our method under varying degrees of travel time uncertainty through four coefficient of variation levels: 
\begin{itemize}
    \item Minimal Variation (\textit{cv} = 0.05)
    \item Low Variation (\textit{cv} = 0.25)
    \item Moderate Variation (\textit{cv} = 0.50)
    \item High Variation (\textit{cv} = 1.00).
\end{itemize}
We evaluate the method across 30 realizations of link travel times and prediction errors, using 30 random seeds (each corresponding to a single~$\omega\in\Omega$ -  \S{}\ref{sec:Stochastic-Modeling-of-Travel-Times}), and quantify the resulting scenario-to-scenario uncertainty with standard deviations computed as described in Appendix~\ref{A.PerformanceMetrics}.
Following our predict-then-optimize framework, the optimization phase utilizes predicted travel times $\hat{t}_{ij}$ to design adaptive networks, while the evaluation phase assesses actual performance using realized travel times $t_{ij}^\omega,\omega\in\Omega$. This separation enables us to quantify the impact of the prediction error on the quality of network adaptation.

% -------------------------------------------- tablemanlCOMP
\begin{table}[tbp]
\centering
\footnotesize
\setlength{\tabcolsep}{0pt}
\renewcommand{\arraystretch}{1.05}

\begin{threeparttable}
\begin{tabular*}{\textwidth}{@{\extracolsep{\fill}}ccccccc@{}}
\toprule
\textbf{No. of lines} & \textbf{Metric}
& \textbf{CV=0}
& \textbf{CV=0.05}
& \textbf{CV=0.25}
& \textbf{CV=0.50}
& \textbf{CV=1.00} \\
\midrule

\multirow{3}{*}{\textbf{4}}
& 20\% improv. share (IS) (\%)
& $5.79\pm2.99$ & $8.25\pm4.65$ & $15.40\pm5.66$
& $24.63\pm7.91$ & $30.57\pm11.60$ \\
& Op. Cost Red. (\%)
& $10.20\pm5.31$ & $13.28\pm7.37$ & $16.86\pm6.74$
& $24.06\pm12.01$ & $39.20\pm12.94$ \\
& Line Overlap (\%)
& $84.34\pm5.44$ & $83.79\pm7.73$ & $78.83\pm6.65$
& $77.53\pm5.27$ & $77.26\pm5.19$ \\
\addlinespace[0.25em]
\midrule

\multirow{3}{*}{\textbf{6}}
& 20\% improv. share (IS) (\%)
& $8.96\pm4.47$ & $9.53\pm3.60$ & $16.31\pm5.10$
& $23.46\pm7.84$ & $31.50\pm10.46$ \\
& Op. Cost Red. (\%)
& $25.94\pm4.71$ & $28.94\pm5.40$ & $32.10\pm5.69$
& $36.75\pm9.55$ & $48.90\pm10.19$ \\
& Line Overlap (\%)
& $79.58\pm2.29$ & $81.20\pm2.80$ & $78.16\pm3.94$
& $79.63\pm4.47$ & $79.16\pm4.08$ \\
\addlinespace[0.25em]
\midrule

\multirow{3}{*}{\textbf{8}}
& 20\% improv. share (IS) (\%)
& $11.52\pm2.67$ & $13.99\pm3.90$ & $18.62\pm5.64$
& $26.84\pm8.76$ & $35.71\pm11.53$ \\
& Op. Cost Red. (\%)
& $26.97\pm4.45$ & $29.92\pm3.84$ & $31.51\pm5.27$
& $35.35\pm7.98$ & $45.35\pm9.23$ \\
& Line Overlap (\%)
& $83.58\pm4.07$ & $84.95\pm3.07$ & $82.51\pm4.03$
& $80.80\pm5.25$ & $81.14\pm4.23$ \\
\addlinespace[0.25em]
\midrule

\multirow{3}{*}{\textbf{10}}
& 20\% improv. share (IS) (\%)
& $8.83\pm3.15$ & $9.74\pm4.54$ & $14.27\pm5.61$
& $21.46\pm8.88$ & $31.21\pm11.56$ \\
& Op. Cost Red. (\%)
& $33.62\pm4.18$ & $34.65\pm4.93$ & $37.23\pm5.53$
& $40.62\pm8.71$ & $52.59\pm9.17$ \\
& Line Overlap (\%)
& $77.57\pm3.40$ & $78.14\pm2.76$ & $76.96\pm3.22$
& $78.77\pm4.09$ & $78.34\pm3.46$ \\
\addlinespace[0.25em]
\midrule

\multirow{3}{*}{\textbf{12}}
& 20\% improv. share (IS) (\%)
& $9.76\pm3.57$ & $10.34\pm3.19$ & $14.74\pm3.51$
& $21.44\pm6.96$ & $32.44\pm11.91$ \\
& Op. Cost Red. (\%)
& $35.94\pm4.95$ & $37.51\pm5.54$ & $39.56\pm5.18$
& $44.20\pm6.90$ & $54.27\pm9.00$ \\
& Line Overlap (\%)
& $80.69\pm2.98$ & $80.95\pm3.40$ & $80.61\pm2.78$
& $81.46\pm3.27$ & $80.77\pm4.11$ \\

\bottomrule
\end{tabular*}

\begin{tablenotes}
\scriptsize
\item No. of lines: number of bus lines; CV: coefficient of variation.
20\% improvement share (IS$^{20\%}$): fraction of OD pairs whose travel-
time improvement is at least 20\%;
Op. Cost Red. (OCR): Operational Cost Reduction (reduction in total
vehicle-hours; see \S\ref{sec:obj}).
The static route-frequency solutions for each number of bus line are taken from~\citet{ahern2022approximate}; entries report mean $\pm$ standard deviation over 30 realizations.
\end{tablenotes}

\end{threeparttable}
\caption{Performance comparison of the solutions from~\citet{ahern2022approximate} under different variability levels.}
\label{tab:routeNumberCV}
\end{table}
% -------------------------------------------- tablemanlCOMPEND

% =============================================
% 5.4 Performance Metrics
% =============================================
\subsubsection{Performance metrics}
\label{sec:performMETRICS}

To evaluate the adaptive network under uncertain travel times, we report three
application-specific performance indicators aligned with the objectives
in~\eqref{eq:multi_objective}:

\begin{itemize}
    \item \textbf{20\% improv. share (\%) (IS$_{20}$):} The percentage of OD pairs that achieve a travel time reduction of at least 20\% compared to the original PT network.
    \item \textbf{Operational Cost Reduction (OCR):} percentage reduction in
    total vehicle-hours operated relative to the original PT network (the operational-cost proxy defined in \S\ref{sec:obj}).
    \item \textbf{Line Overlap (LO):} Jaccard-based percentage overlap~\eqref{eq:jaccard} between the edge sets of the original and adaptive PT networks.
\end{itemize}

All three metrics depend on realized travel times and are therefore random
variables; in the numerical results, we report their empirical means (and
standard deviations) over 30 independent realizations.\footnote{Formal definitions are provided in App.~\ref{A.PerformanceMetrics}.}
\rev{zg}{}{We distinguish the \emph{objective} from the \emph{reporting metric}: the optimization objective $Z_1$ is demand-weighted by $q_k$ (\S\ref{sec:obj}), whereas mtric IS$_{20}$ counts OD pairs without demand weighting and therefore answers ``for how many OD relations does service improve,'' not ``how many passenger-trips benefit.'' The aggregate, demand-weighted passenger benefit is captured separately by the Avg.~Travel Time Improvement metric reported for the Beijing network (\S\ref{sec:caseStudy}, Table~\ref{tab:stop_skipping_comparison}), so that both the optimization and the reported benefit reflect the distributional, demand-driven question of which travellers gain from the redesign~\citep{fan2011bilevel, pereira2017equity}.}

% ------------------------------------------------------
% % ++++++++++++++++++++ 6.2.1 - 基于Mandl_旅行时间分析 (不同seed)
% ------------------------------------------------------

\subsection{Main performance results on Mandl's benchmark}\label{sec:mandl_RESULT1}

\rev{zg}{We organize the Mandl results in three steps. First, we report aggregate performance across different variability levels and bus-line counts. We then compare full route redesign with a frequency-only baseline and examine sensitivity to forecast information. Finally, we inspect trip-time distributions to show how the aggregate gains are distributed across OD pairs.}{We organize the Mandl results in four parts. First, we report aggregate performance across variability levels and route quantities. Second, we inspect trip travel-time distributions to show how aggregate gains are distributed across OD pairs. Third, we compare full route redesign with a frequency-only baseline. Finally, we examine sensitivity to forecast information.}

\subsubsection{Aggregate Performance across Variability Levels}

We use predicted travel times $\hat{t}_{ij}$ in Eq.~\eqref{eq:predTRAVEL} as input to the optimization problem, and we evaluate the resulting adaptive network against realized link travel times and prediction errors across multiple scenarios.
Table~\ref{tab:routeNumberCV} summarizes average performance across 30 independent realizations.
It shows pronounced OD-pair travel-time improvements under high variability ($cv=1.0$), while improvements remain significant even at lower variability levels.
% \footnote{
As expected, when variability is $\textit{cv}=0$, the central reference point (described in the previous paragraph) corresponds exactly to~$\pazocal G_\text{PT}^\text{ahe}$. 
% } 
As travel-time variability increases, both 20\% Improvement Share (IS$_{20}$) and Operational Cost Reduction (OCR) rise systematically, reflecting the increasing value of adaptive redesign under highly fluctuating conditions.
At low variability ($cv=0.05$), a sizable $IS_{20}$ is still observed. This is because extreme link travel times are still possible in reality, even at low $cv$ (the possibility of such extreme values is preserved by our choice of a Log-normal distribution~\S{}\ref{sec:log_normal}. For a small-scale scenario like Mandl’s network, such extremes inflate original trip times $t_k^0, k\in\pazocal K$.
% 尽管是0.05，但还是有extrame值，导致有提升
Accordingly, by exploring alternative trade-offs within the same objective space, the proposed approach can reduce passenger travel time while maintaining or improving operational costs, even in near-deterministic settings.
Across all variability levels and redesign schemes, Line Overlap (LO) remains high, indicating that these improvements result from targeted adjustments rather than wholesale restructuring. 
%The relatively narrow confidence intervals of route overlap reported in Table~\ref{tab:routeNumberCV} further confirm that the observed trends persist across realizations and are not driven by outlier cases.

% \rev{zg}{}{These results also clarify why the compromise selection is non-trivial. The three objectives pull in different directions: reducing passenger travel-time deterioration ($Z_1$) can require more structural change (higher $Z_3$, observed as lower RO), while reducing vehicle-hours ($Z_2$) can shorten service and reduce passenger benefit. Thus $S1$--$S3$ should be read as balanced Pareto compromises rather than single-objective optima; a $(Z_1,Z_2)$ projection of the static-case front is shown in Fig.~\ref{fig:NSGA3SELF} (Appendix~\ref{sec:algorithmic-framework-and-evaluation}).}

% ------------------------- FIG OD旅行时间分布比较
% cas-sc figure/table options are key-value: use pos=..., not a bare [htbp]
\begin{figure}[htbp]
    \centering
    % 左侧图片
    \begin{subfigure}[b]{0.47\textwidth}
        \centering
        \includegraphics[width=\textwidth]{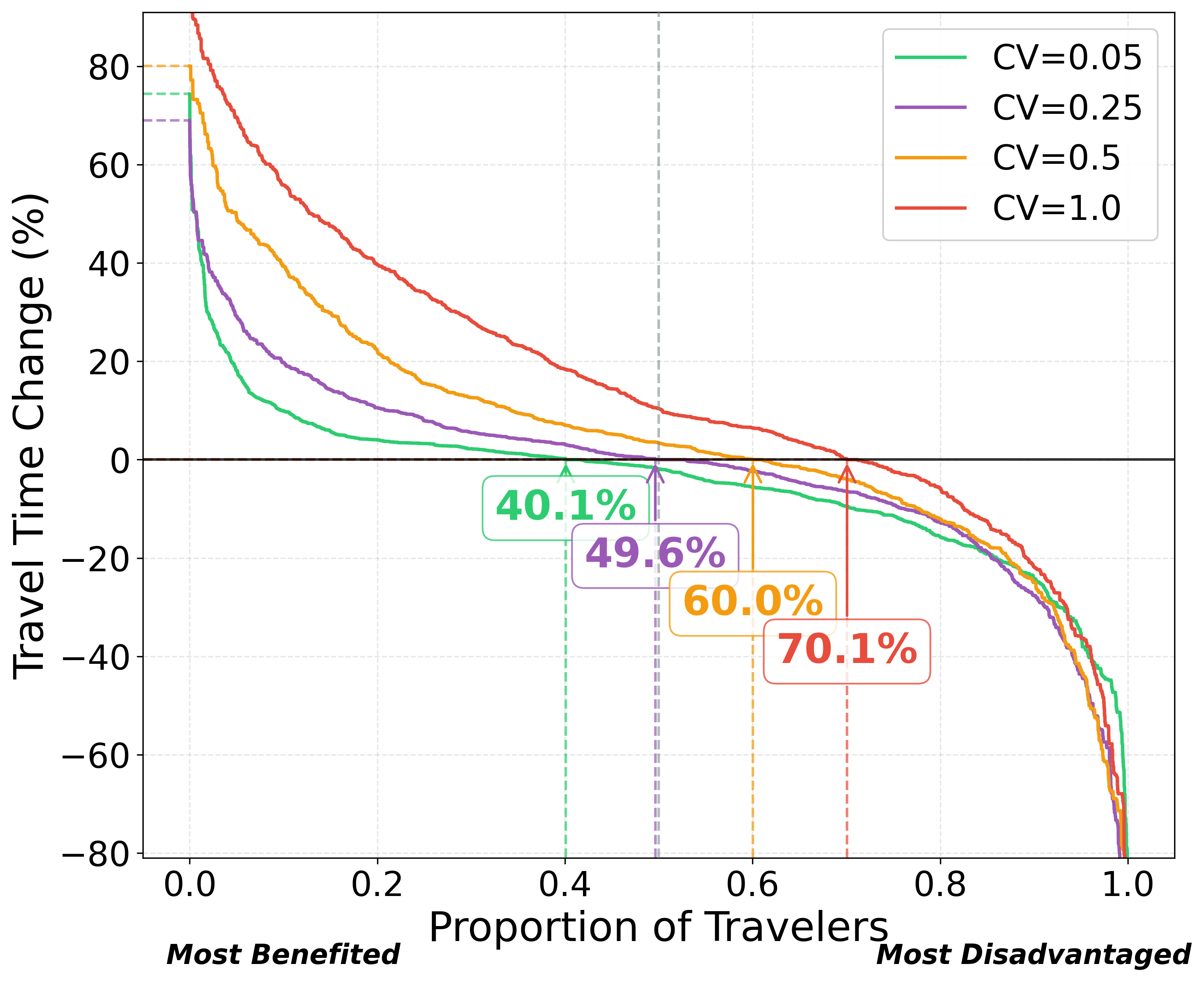}
        \caption{$\pazocal{G}_{\text{PT}}^{\text{adpt}}$ vs $\pazocal{G}_{\text{PT}}^{\text{ahe}}$~\citep{ahern2022approximate}}
        \label{fig:self_comparison}
    \end{subfigure}
    \hfill
    % 右侧图片
    \begin{subfigure}[b]{0.47\textwidth}
        \centering
        \includegraphics[width=\textwidth]{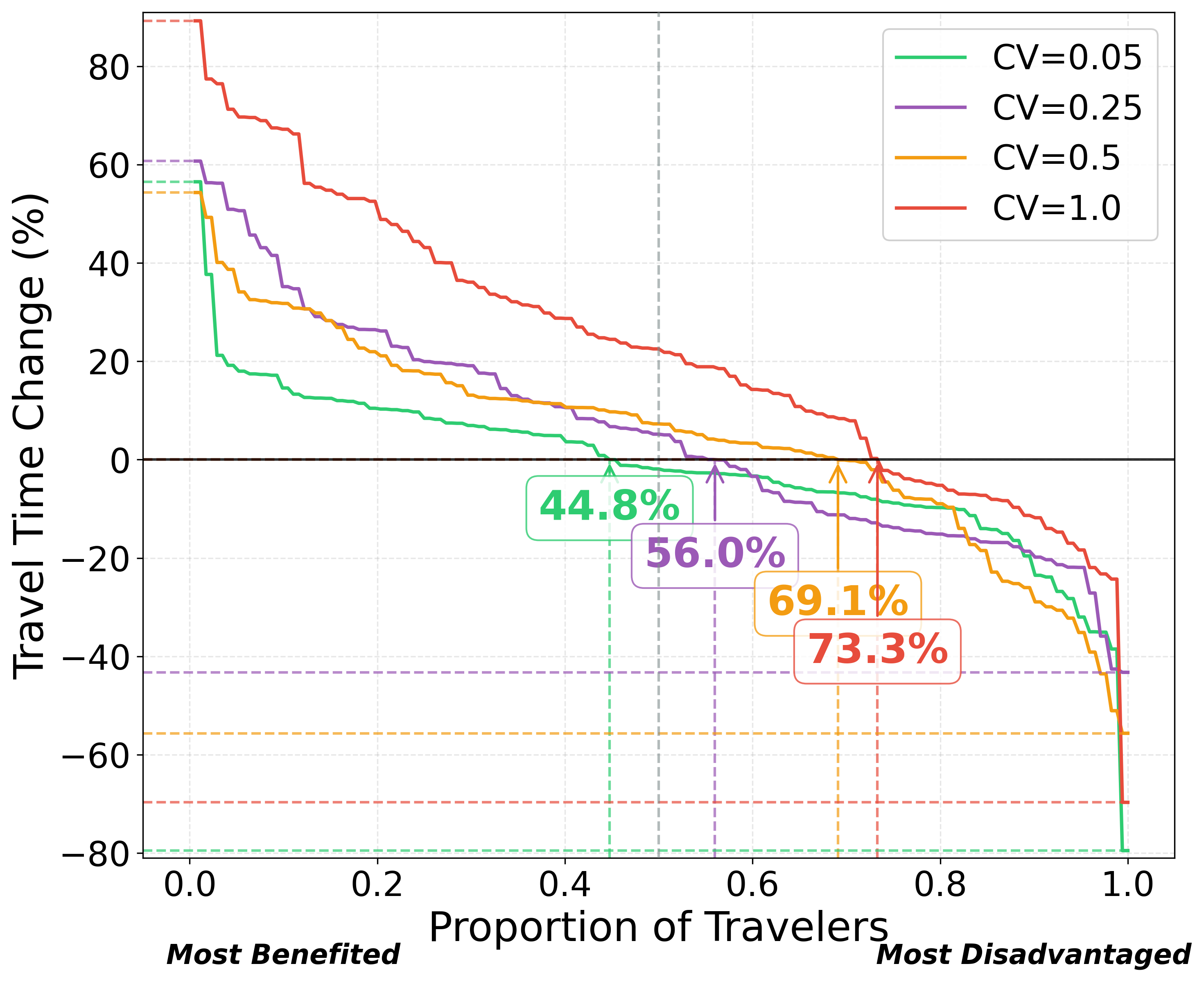}
        \caption{$\pazocal{G}_{\text{PT}}^{\text{adpt}}$ vs $\pazocal{G}_{\text{PT}}^{\text{ver}}$~\citep{vermeir2021exact}}
        \label{fig:ver_comparison}
    \end{subfigure}

    \caption{Trip Time Distribution Comparison Under Various Variability Conditions in a randomly chosen realization.}
    \label{fig:od_comparison}
\end{figure}
% ------------------------- FIG OD旅行时间分布比较end

\subsubsection{Trip Travel Time Distribution Analysis}

%In the following sections, we will show the performance of solution~$S1$.
\rev{zg}{To provide a conservative performance evaluation, we considered the ``worst'' realization, i.e., the seed that yielded the worst results in terms of the metrics of Table~\ref{tab:routeNumberCV}.}{For illustration, in this subsection we report the results for a single randomly chosen realization.}
Fig.~\ref{fig:od_comparison} illustrates trip time distributions across varying coefficients of variation (\textit{cv}). It highlights the benefits of our adaptive public transport planning method, as travel time uncertainty increases from minimal to high levels, against current network design approaches, represented by the original static network $\pazocal{G}_{\text{PT}}^{\text{ahe}}$ and Vermeir's network $\pazocal{G}_{\text{PT}}^{\text{ver}}$, in accordance to \S{}\ref{sec:compared-PT-networks}.
%
%As shown in Fig.~\ref{fig:self_comparison} and Fig.~\ref{fig:ver_comparison}, our adaptive approach demonstrates a clear pattern of increasing effectiveness as travel time uncertainty grows.
At low variability levels (\textit{cv}=0.05), less than 50\% of OD pairs improve their travel time. This occurs because these static networks are already well-optimized for stable conditions, leaving little room for enhancement.
Moreover, our optimization operates on predicted travel times that inevitably include errors. This is why, as illustrated in Fig.~\ref{fig:rollFramework}, we would not trigger any reconfiguration, faced with minor fluctuations in link travel times, such as those experienced with $\textit{cv}=0.05$.

However, as variability increases, the benefits of our approach ($\pazocal{G}_{\text{PT}}^{\text{adpt}}$) become pronounced, against the optimized static networks $\pazocal G_\text{PT}^\text{ahe}$ and $\pazocal G_\text{PT}^\text{ver}$.
For instance, under high uncertainty (\textit{cv}=1.0), in the considered realization, Fig.~\ref{fig:od_comparison} shows that about $73.3\%$ of travelers experience reduced travel time, with a reduction of at least $30\%$ for roughly $30\%$ of travelers.

\rev{zg}{This confirms the superiority of our proactive network redesign against current approaches, all keeping the network unchanged, in particular in highly fluctuating conditions.}{This suggests the potential benefit of our proactive network redesign over the selected static baselines, which keep the network unchanged, in particular in highly fluctuating conditions.}

Beyond this single-realization view, Table~\ref{tab:routeNumberCV} reports results averaged across 30 realizations: under high uncertainty (\textit{cv}=1.0), our adaptive network $\pazocal{G}_{\text{PT}}^{\text{adpt}}$ improves the travel time of about $30.6\%$ of OD pairs by at least $20\%$, relative to $\pazocal{G}_{\text{PT}}^{\text{ahe}}$. Further OD-level improvement-rate analysis aggregated across all 30 seeds is provided in App.~\ref{App:D3_ImprovementICDF}.

% -------------------------- TableFRE
\begin{table*}[htbp]
\centering
\small
\setlength{\tabcolsep}{0pt}
\renewcommand{\arraystretch}{1.12}

\begin{threeparttable}
\begin{tabular*}{\textwidth}{@{\extracolsep{\fill}}cccccc@{}}
\toprule
\textbf{Method} & \textbf{Metric}
& \textbf{CV=0.05}
& \textbf{CV=0.25}
& \textbf{CV=0.50}
& \textbf{CV=1.00} \\
\midrule

\multirow{3}{*}{\textbf{Frequency Only}}
& 20\% improv. share (\%)
& $3.96\pm2.94$ & $8.09\pm4.35$
& $12.94\pm6.92$ & $16.47\pm9.53$ \\
& Op. Cost Red. (\%)
& $12.85\pm4.66$ & $14.86\pm5.78$
& $19.27\pm9.04$ & $29.40\pm15.04$ \\
& Line Overlap (\%)
& $100.00\pm0.00$ & $100.00\pm0.00$
& $100.00\pm0.00$ & $100.00\pm0.00$ \\
\addlinespace[0.25em]
\midrule

\multirow{3}{*}{\textbf{Ours}}
& 20\% improv. share (\%)
& $8.25\pm4.65$ & $15.40\pm5.66$
& $24.63\pm7.91$ & $30.57\pm11.60$ \\
& Op. Cost Red. (\%)
& $13.28\pm7.37$ & $16.86\pm6.74$
& $24.06\pm12.01$ & $39.20\pm12.94$ \\
& Line Overlap (\%)
& $83.79\pm7.73$ & $78.83\pm6.65$
& $77.53\pm5.27$ & $77.26\pm5.19$ \\

\bottomrule
\end{tabular*}

\begin{tablenotes}
\footnotesize
\item FO: fixed-route frequency-only adjustment baseline $\pazocal{G}_{\text{PT}}^{\text{freq}}$; Ours: our full network redesign method.
CV: coefficient of variation.
20\% improvement share (IS$^{20\%}$-\S{}\ref{sec:performMETRICS}): fraction of OD pairs whose travel-time improvement is at least 20\%;
Op. Cost Red. (OCR): Operational Cost Reduction (reduction in total vehicle-hours; see \S\ref{sec:obj}).
The 4-line static line-frequency solution is taken from~\citet{ahern2022approximate}; entries report mean $\pm$ standard deviation over 30 realizations.
\end{tablenotes}

\end{threeparttable}
\caption{\rev{zg}{}{Performance comparison of full network redesign and fixed-line frequency adjustment for the 4-line Mandl instance.}}
\label{tab:frequencyBaselineCV}
\end{table*}
% -------------------------- TableFREEND

\subsubsection{Frequency-Only Baseline}

\begin{figure}[htbp]
\centering
\includegraphics[width=\textwidth]{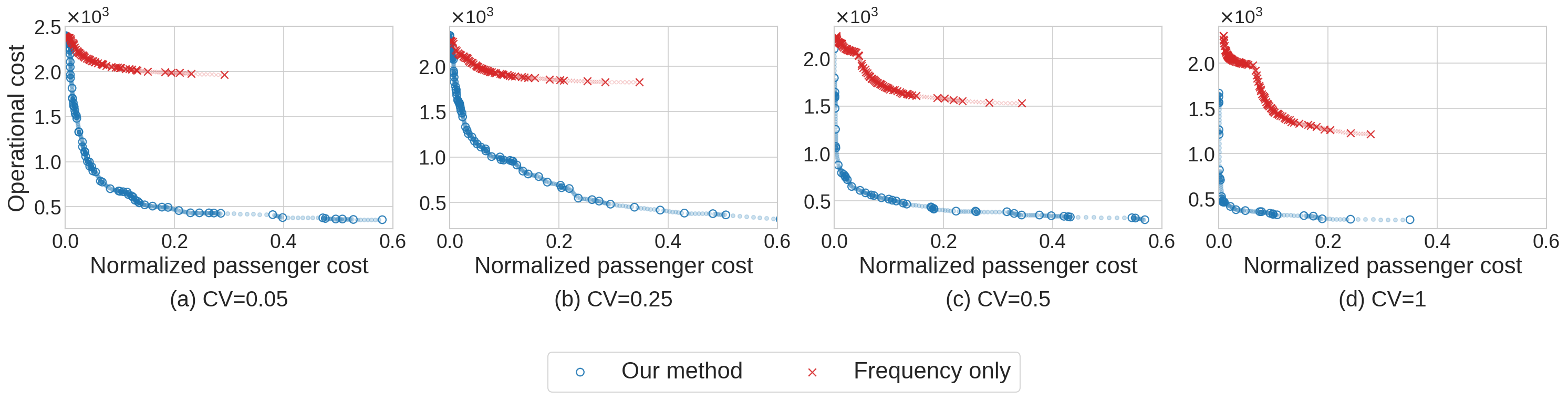}
\caption{\rev{zg}{}{Travel-time, cost ($Z_1$, $Z_2$) Pareto fronts for the 4-line Mandl instance under four travel-time variability levels.}}
\label{fig:objectiveZ1Z2Frontiers}
\end{figure}

We further evaluate whether adjusting frequencies alone can account for the observed gains, comparing our approach against $\pazocal{G}_{\text{PT}}^{\text{freq}}$, obtained by solving the same problem presented in~\S\ref{sec:obj}, by additionally imposing that all lines remain unchanged, and only allowing changing in line frequencies and vehicle allocation to line.
Fig.~\ref{fig:objectiveZ1Z2Frontiers} plots the 2D projection of the Pareto front, highlighting the trade-off between passenger travel-time $Z_1$ (which in the figure is normalized between~$0$ and~$1$) and operator cost in vehicle-hours ($Z_2$). Each point maps to a specific network design, with the lower-left corner representing the ideal direction for both minimization objectives. Meanwhile, the structural-continuity objective ($Z_3$) is reported as line overlap in Table~\ref{tab:frequencyBaselineCV}.
The frequency-only baseline $\pazocal{G}_{\text{PT}}^{\text{freq}}$ improves service under variability, but it captures substantially less of the OD-level improvement than full redesign; for example, at $cv=1$, its 20\%-improvement share is roughly half that of our method.
The results suggest topology change is beneficial when congestion patterns create larger spatial mismatches in the original network.
While frequency-only adjustment keeps the line topology fixed, the proposed line redesign can reroute service away from deteriorated links.

% -------------------------- Figure
% -------------------------- FigureEND

\subsubsection{Forecast-error Sensitivity}
We now test whether the proposed proactive redesign is driven by the intrinsic variability of link travel times rather than by forecast errors. We indeed want to rule out the hypothesis that the design adapts to the predictor, rather than to traffic conditions. To this aim, we hold all realized travel times~$\{t_{ij}^{\omega}\}$ fixed and vary only the error-scale parameter $\alpha$ in Eq.~\eqref{eq:predTRAVEL_alpha}. Recall that the error impacts the predicted travel times~$\hat t_{ij}$ (via~\eqref{eq:link_time_error}), which are the ones considered within the problem resolution~(as specified in~\S{}\ref{sec:obj}).
Here, $\alpha=0$ gives a perfect forecast, $\alpha=1$ reproduces the empirical prediction error, and $\alpha>1$ is evaluated for an amplified-error stress test.

For each setting, we score the resulting Pareto front by its hypervolume and compare it, scenario by scenario, to the fronts obtained by running the optimization in two cases: (i) using nominal link travel times ($t_{ij}^0$),  and (ii)~using perfect forecast ($\alpha=0$, which implies~$\hat t_{ij}=t_{ij}$). Detailed metrics and results are in App.~\ref{app:forecastValue} and indicate that, at low variability (e.g., low~$cv$), hypervolume is nearly the same in the two cases: the redesign already reaches an equivalent front on its own, so even a perfect forecast adds little to just using nominal values.
A real gap opens only once variability is large enough to create a genuine redesign opportunity. In this case, the Pareto front issued by the proposed optimization changes little when varying the forecast error factor~$\alpha$.
%there, the empirical forecast closes nearly all of the gap, and moderately amplifying further the forecast error factor~$\alpha$ does not erode this result.
The method's gains are therefore driven mainly by the redesign opportunities that variability creates, not by a dependence on high-precision forecasts.

\subsection{Large-scale validation on the Beijing network}\label{sec:caseStudy}

% ------------------------
\begin{figure}[htbp]
\centering
\begin{subfigure}[t]{0.48\textwidth}
    \centering
    \includegraphics[width=\linewidth]{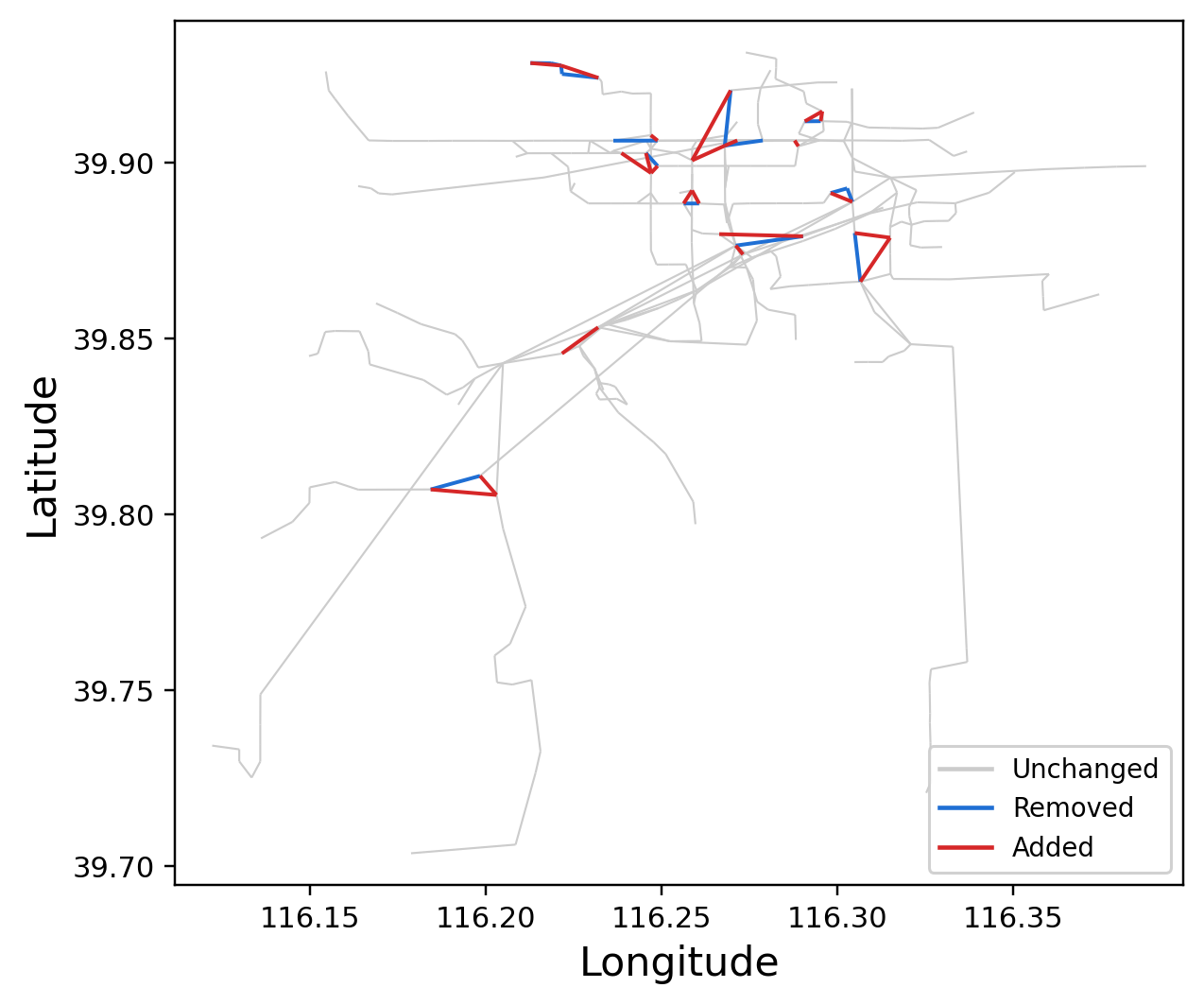}
    \caption{$cv = 0$}
    \label{fig:beijingCV0}
\end{subfigure}
\hfill
\begin{subfigure}[t]{0.48\textwidth}
    \centering
    \includegraphics[width=\linewidth]{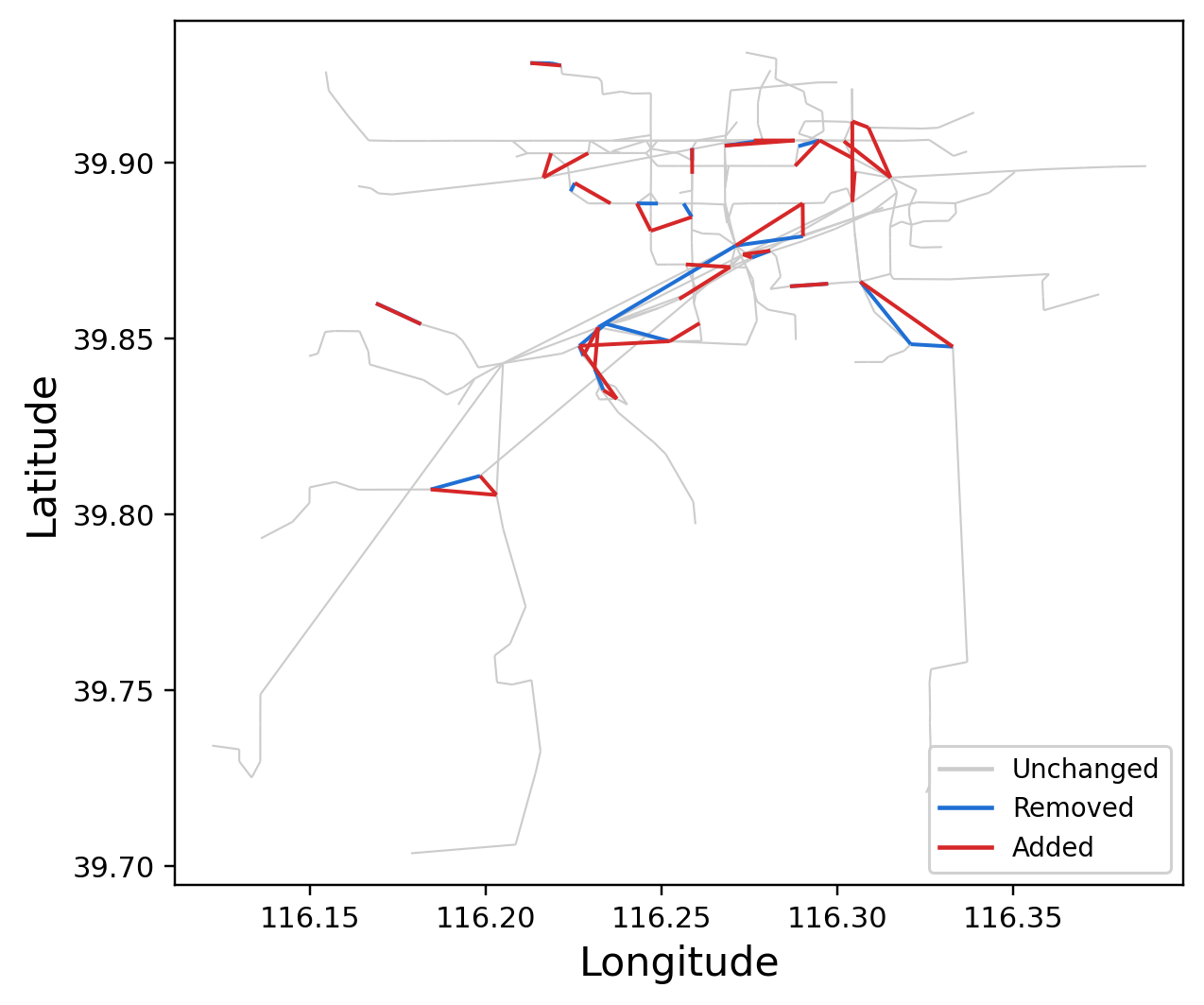}
    \caption{$cv = 1$}
    \label{fig:beijingCV1}
\end{subfigure}
\caption{\rev{zg}{Original and adjusted PT networks for the Beijing case study.}{Adapted Beijing PT networks under low ($cv=0$) and high ($cv=1$) travel-time variability. Additional variability levels are shown in App.~\ref{app:beijing_networks}.}}
\label{fig:beijingRouteNetworks}
\end{figure}

% =============================================
% [Deployment-boundary material integrated into \S\ref{sec:numericalRESULT} (framework output / scope of use) and \S\ref{sec:timeframe} (operational gates, near-term deployable forms, continuity calibration).]

%%% --------------------------------table
\begin{table}[htbp]
\centering
\scriptsize
\setlength{\tabcolsep}{3pt}
\renewcommand{\arraystretch}{1.12}
\begin{threeparttable}
    \begin{tabular}{c c c c c c c c}
    \toprule
    \makecell[c]{\textbf{Method}}
    & \makecell[c]{\textbf{Thr.}\\ \scriptsize Skip threshold}
    & \makecell[c]{\textbf{CV}\\ \scriptsize Coeff. var.}
    & \makecell[c]{\textbf{LO}\\ \scriptsize Line Overlap (\%)}
    & \makecell[c]{\textbf{OCR}\\ \scriptsize Op. cost red. (\%)}
    & \makecell[c]{\textbf{Avg.TTI}\\ \scriptsize Travel-time impr. (\%)}
    & \makecell[c]{\textbf{Unreach. OD}\\ \scriptsize Unreachable OD (\%)}
    & \makecell[c]{\textbf{Unserv. Dem.}\\ \scriptsize Unserved demand (\%)} \\
    \midrule

    \multirow{15}{*}{\makecell[c]{Stop-skip.}}
    & \multirow{5}{*}{\makecell[c]{25\%}}
        & 0.00 & 91.16 & 1.24  & 3.26  & 4.89  & 3.89  \\
    & & 0.05 & 91.71 & 1.15  & 2.73  & 4.29  & 3.46  \\
    & & 0.25 & 56.02 & 6.58  & 12.43 & 29.12 & 27.39 \\
    & & 0.50 & 47.58 & 6.95  & 14.08 & 35.28 & 32.99 \\
    & & 1.00 & 42.18 & 12.71 & 14.71 & 38.93 & 36.04 \\
    \cmidrule{2-8}

    & \multirow{5}{*}{\makecell[c]{50\%}}
        & 0.00 & 95.52 & 0.61  & $-$0.10 & 1.71  & 1.53  \\
    & & 0.05 & 95.24 & 0.63  & 1.94  & 1.71  & 1.53  \\
    & & 0.25 & 80.32 & 4.59  & 10.19 & 11.70 & 11.40 \\
    & & 0.50 & 56.13 & 8.69  & 19.42 & 28.63 & 27.56 \\
    & & 1.00 & 47.83 & 12.39 & 19.71 & 35.28 & 32.99 \\
    \cmidrule{2-8}

    & \multirow{5}{*}{\makecell[c]{100\%}}
        & 0.00 & 98.30 & 0.32  & 0.45  & 0.03  & 0.02  \\
    & & 0.05 & 98.30 & 0.35  & 0.49  & 0.03  & 0.02  \\
    & & 0.25 & 93.87 & 1.63  & 4.08  & 2.77  & 2.09  \\
    & & 0.50 & 78.57 & 6.68  & 17.77 & 12.23 & 11.61 \\
    & & 1.00 & 55.88 & 13.49 & 15.11 & 28.05 & 26.84 \\
    \midrule

    \multirow{5}{*}{\makecell[c]{\textbf{Ours}}}
    & \multirow{5}{*}{\makecell[c]{--}}
        & 0.00 & 90.08 & 0.50  & 11.40 & -- & -- \\
    & & 0.05 & 87.34 & 1.53  & 11.60 & -- & -- \\
    & & 0.25 & 86.51 & 3.85  & 16.99 & -- & -- \\
    & & 0.50 & 86.39 & 12.02 & 18.44 & -- & -- \\
    & & 1.00 & 85.42 & 22.56 & 25.78 & -- & -- \\
    \bottomrule
    \end{tabular}

    \begin{tablenotes}
    \scriptsize
    \item CV: coefficient of variation. Thr.: skip threshold used by the stop-skipping baseline.
    LO: Line Overlap. OCR: Operational Cost Reduction, measured as the reduction in total vehicle-hours; see \S\ref{sec:obj}.
    TTI: Travel-Time Improvement; see~\eqref{eq:avgTTI}.
    Unreach.\ OD and Unserv.\ Dem. report the fraction of OD pairs and demand, respectively, whose connectivity is broken by stop-skipping. These metrics are not applicable to our method, which removes no stops and preserves full OD connectivity.
    \end{tablenotes}
\end{threeparttable}
\caption{\rev{zg}{}{Comparison of our method and stop-skipping under varying travel-time variability on the Beijing network.}}
\label{tab:stop_skipping_comparison}
\end{table}

% \todoi{aa: In Table~\ref{tab:stop_skipping_comparison} add two lines of heading with a bit clearer writing}
%%% --------------------------------tableEND

We now move the analysis to the real Beijing network. 
The original and adjusted PT networks are shown in Fig.~\ref{fig:beijingRouteNetworks}. \rev{aa}{}{We color a link in red if there exists at least one line~$\ell$ that used that link in~$\pazocal G_\text{PT}^\text{orig}$ and that does not use that link in~$\pazocal G_\text{PT}^\text{adpt}$.
We color a link in blue if there exists at least one line~$\ell$ that did not use that link in~$\pazocal G_\text{PT}^\text{orig}$ and that does use that link in~$\pazocal G_\text{PT}^\text{adpt}$.
}

We compare our method against a \emph{stop-skipping} operational baseline on the Beijing network.
For each variability level, links whose predicted travel time exceeds the original by more than a preset threshold are classified as degraded.
The stop-skipping strategy then removes intermediate stops adjacent to degraded links whenever a direct shortcut is available, while leaving the overall line structure and vehicle allocation unchanged.
The threshold thus controls aggressiveness: a 25\% threshold skips stops whenever a link's predicted time grows by more than a quarter; a 100\% threshold acts only when it doubles.
This baseline represents a simple, local response to travel-time disruption, but structurally constrained: it cannot reroute vehicles or reassign fleet, and skipping stops may disconnect part of the demand at a stop from the network entirely, obliging users to walk to different stops.
Table~\ref{tab:stop_skipping_comparison} reports line overlap (LO), operational cost reduction (OCR), travel-time improvement (TTI), and the resulting fractions of disconnected OD pairs and disconnected demand (corresponding to users that would need to walk to other stops) for all threshold settings alongside our method.

In addition to the previous metrics, we also compute the following.
\begin{itemize}
    \item \textbf{Avg. Travel Time Improvement (Avg.TTI):} percentage reduction in aggregate travel time across all OD pairs $k$ relative to the original network:
    \begin{align}
    \label{eq:avgTTI}
        \mathrm{Avg. TTI} = \frac{\sum_{k \in \pazocal{K}} (t_k^0 - t_k^*)}{\sum_{k \in \pazocal{K}} t_k^0} \times 100\%, \text{ where } t_k^0, t_k^* \text{ are original and adaptive travel times.}
    \end{align}
\end{itemize}
The results in Table~\ref{tab:stop_skipping_comparison} show that the original network exhibits some inefficiencies: the solution of our optimization yields an 11.40\% Avg.\ TTI even at $cv = 0$. \rev{zg}{This reflects static re-optimisation of the incumbent network rather than prediction-driven adaptation. The incremental value of adaptation is the additional gain observed as $cv$ increases.}{This gain at $cv=0$ does not stem from prediction-driven adaptation: with no travel-time variability, there is nothing to forecast or to adapt to. Rather, it reflects a one-off re-optimisation of an incumbent network that, like most real-world bus systems, is not a mathematical optimum. Operational networks evolve incrementally over many years, with lines added or amended in response to historical demand and retained for service continuity, and are shaped by administrative, budgetary, and political constraints rather than periodic global re-optimisation, for which no optimization techniques are routinely used in practice~\citep{duran2022survey}; as travel demand and land use shift, the fixed topology gradually drifts away from the efficient configuration. Our optimizer recovers part of this accumulated static inefficiency, which explains the $11.40\%$ improvement even under nominal conditions. The genuine value of \emph{adaptation} is therefore the additional gain observed as $cv$ increases, beyond this one-time static correction.}
The benefit grows with variability: at higher $cv$, the algorithm produces designs increasingly different from the static original, achieving both travel-time improvements and line shortening, which in turn lowers operational cost.

\rev{zg}{}{Table~\ref{tab:stop_skipping_comparison} exposes a basic tension in stop-skipping: operational benefit comes at the cost of connectivity. Aggressive thresholds improve travel time on the reachable network but disconnect a large share of OD pairs as variability grows, whereas conservative thresholds keep connectivity yet act too rarely to help and still lose connectivity at high variability. No single threshold resolves this trade-off.

Our method achieves larger travel-time gains for two structural reasons. First, it redesigns entire line topologies rather than dropping individual stops, so vehicles can bypass congested corridors and serve demand more directly; this reaches reductions that stop-level changes cannot, since the bottleneck is usually the line's overall routing rather than a single link. Second, it jointly reallocates the fleet across lines, redeploying vehicles freed from shortened routes where they help most, while stop-skipping keeps the fleet fixed. Because no stops are removed, all OD pairs stay reachable and the reported TTI covers the full demand, making the gains comparable across variability levels without the connectivity issues of every stop-skipping variant.}

% =============================================
% 6.Conclusion
% =============================================
\FloatBarrier
\section{Conclusion}\label{sec:conclu}
This paper \rev{aa}{}{systematically studied the potential of proactive bus line redesign, targeted at limited parts of the network, in maintaining user and operator performance in the case of severe traffic fluctuations.}
\rev{aa}{introduced a novel framework for redesigning PT lines for the hour ahead from short-term forecasts, acting only when the forecast deviation from nominal conditions exceeds a preset tolerance, and keeping each redesign close to the network already in service.}{}
%real-time PT proactive network redesign that maintains high performance under uncertain traffic conditions
We formulated a multi-objective optimization problem balancing passenger travel time, operational costs, and network stability. Our predict-then-optimize approach relies on an NSGA-III-based optimization framework that proactively redesigns bus lines based on predicted line travel times. 

We evaluated our approach against the inevitable prediction errors of well-established traffic-predicting models. To this aim, we expose the optimizer to realistic forecast inaccuracies sampled from a well-established deep learning-based prediction model (DCRNN).
Experiments on Mandl's network and a real-world Beijing network demonstrate effective adaptation across varying levels of travel-time uncertainty. 
Results indicate that, faced with non-negligible link travel-time variability, the proposed approach \rev{zg}{yields substantial OD-level trip travel-time reductions, while achieving operational cost savings and preserving}{yields OD-level trip travel-time reductions, while achieving operational cost savings over the selected baselines and scenarios considered, and preserving} a high degree of structural similarity with the baseline network.
\rev{zg}{These findings highlight the importance of moving beyond static network design toward continuous, adaptive redesign in dynamic operating environments.}{Overall, the gains grow with travel-time variability and come with high line overlap, indicating that the benefit does not rely on large structural change.}
A natural extension is to enrich the forecast-error model with explicit spatial and temporal correlation, allowing network-wide common shocks to be represented.
A further direction is to
%relax the uncongested user-optimal assignment (Observation~\ref{obs:det_routing}) by 
incorporate in-vehicle crowding, for instance through a capacitated or congestion-dependent (frequency-based common-line) passenger assignment; this would let the redesign account for the load it induces on each line and report crowding-aware rather than purely uncongested travel-time gains.

%% (Limitations consolidated into the closing paragraph above; smartphone-availability caveat retained in \S\ref{sec:timeframe}.)

% Numbered list
% Use the style of numbering in square brackets.
% If nothing is used, default style will be taken.
%\begin{enumerate}[a)]
%\item 
%\item 
%\item 
%\end{enumerate}  

% Unnumbered list
%\begin{itemize}
%\item 
%\item 
%\item 
%\end{itemize}  

% Description list
%\begin{description}
%\item[]
%\item[] 
%\item[] 
%\end{description}  

% Uncomment and use as the case may be
%\begin{theorem} 
%\end{theorem}

% Uncomment and use as the case may be
%\begin{lemma} 
%\end{lemma}

%% The Appendices part is started with the command \appendix;
%% appendix sections are then done as normal sections
%% \appendix

% To print the credit authorship contribution details
% \printcredits

% =====================================================
% Declarations (anonymized for double-blind review).
% Author names, affiliations, funding and CRediT are kept
% off the anonymized manuscript and supplied on a separate
% Title Page via the submission system.
% =====================================================
% \section*{Declaration of competing interest}
% The authors declare that they have no known competing financial interests or personal relationships that could have appeared to influence the work reported in this paper.

\section*{CRediT author contribution statement}
\textbf{Zihao GUO:} Conceptualization, Data curation, Methodology, Software,
Validation, Visualization, Writing -- original draft.

\textbf{Andrea ARALDO:} Conceptualization, Methodology,
Project administration, Supervision, Validation, Writing -- original draft,
Writing -- review and editing, Funding acquisition.

\textbf{Fay\c{c}al TOUZOUT:} Conceptualization, Writing -- review and editing.

\textbf{Mounim EL-YACOUBI:} Methodology, Supervision, Writing -- review and editing.

\section*{Acknowledgements}
This work was supported by the French National Research Agency (ANR) under the MuTAS project (ANR-21-CE22-0025-01). The funder had no role in the study design; the collection, analysis, or interpretation of data; the writing of the manuscript; or the decision to submit the article for publication.

\section*{Declaration of competing interest}
The authors declare that they have no known competing financial interests or personal relationships that could have appeared to influence the work reported in this paper.

\section*{Data availability}
The Mandl benchmark instance is publicly available in the cited literature. The Beijing network data are derived from third-party sources cited in the paper. The problem instances and source code will be released in a public repository upon acceptance.

\section*{Declaration of generative AI and AI-assisted technologies in the writing process}
During the preparation of this work, ChatGPT (GPT-5.1, 2025 release; OpenAI, San Francisco, CA, USA) was used to refine the language of the manuscript and ensure compliance with formatting, grammar, and American English conventions. After using these tools, the authors reviewed and edited the manuscript as needed and take full responsibility for the final content of the publication.

%% Loading bibliography style file
%\bibliographystyle{model1-num-names}
\bibliographystyle{unsrtnat}

% Loading bibliography database
\bibliography{references}

% Biography
%\bio{}
% Here goes the biography details.
%\endbio

%\bio{pic1}
% Here goes the biography details.
%\endbio

% ===================================================
% APPENDIX
% ===================================================
\appendix

% --------------------------------------------------- 
% --------------------------------------------------- Appendix 1_Proof of Proposition
% ---------------------------------------------------

\section{Proof and Discussion of Proposition~\ref{prop:temporal_stability}}
\label{app:temporal_stability_proof}

\paragraph{Proof.}
Setting $A=\pazocal{E}^{(k)}$, $B=\pazocal{E}^{(0)}$, and $C=\pazocal{E}^{(k+n)}$.
The Jaccard distance between two sets $X$ and $Y$ is
\[
d_J(X,Y)=1-\frac{|X\cap Y|}{|X\cup Y|}=1-J(X,Y).
\]
Levandowsky and Winter~\citep{levandowsky1971distance} noted that this
dissimilarity measure satisfies the triangle inequality on finite sets,
so that $d_J$ may be regarded as a distance function. Hence,
\[
d_J(A,C)\le d_J(A,B)+d_J(B,C).
\]
Substituting $d_J(X,Y)=1-J(X,Y)$ into the above inequality yields
\[
1-J(A,C)
\le [\,1-J(A,B)\,]+[\,1-J(B,C)\,],
\]
which is equivalent to
\[
J(A,C)\ge J(A,B)+J(B,C)-1.
\]
\rev{zg}{, gives the desired result.}{which gives the desired result.}\hfill$\blacksquare$

\paragraph{Remark.}
This implies that if two consecutive designs exhibit a similarity of $\geq 90\%$ relative to the original (consistent with our numerical results), the similarity between these two consecutive designs is guaranteed to be at least 80\%.
%
% \rev{zg}{}{
% The crucial point is that, because the structural objective $Z_3$ is measured against the \emph{original} network $\pazocal{E}^{(0)}$ at \emph{every} step rather than against the preceding design, the right-hand side does not depend on the gap $n$. Variability therefore does not accumulate or propagate as the day progresses: the similarity between redesigns one hour apart and the similarity between redesigns twelve hours apart obey the \emph{same} lower bound. While the worst-case guarantee $2J_{\min}-1$ is loose (and vacuous at $J_{\min}=0.5$), the empirically observed similarity to the original network is much higher than $J_{\min}$: Tables~\ref{tab:routeNumberCV} and~\ref{tab:stop_skipping_comparison} report route overlaps of roughly $0.77$--$0.90$, which by the inequality above guarantees $J(\pazocal{E}^{(k)},\pazocal{E}^{(k+n)}) \ge 0.54$--$0.80$ between \emph{any} two redesigns of the day. The hourly monitoring policy (\S\ref{sec:timeframe}), which resets the network to $\pazocal{G}_\text{PT}^\text{org}$ whenever deviations fall within tolerance, provides a further safeguard against long-term drift.}
% --------------------------------------------------- 
% --------------------------------------------------- Appendix 2_log_normal
% --------------------------------------------------- 
\section{Stochastic Modelling of Travel Times and Prediction Errors}

% ================================ A.1
\subsection{Parameter Derivation for Log-Normal Distribution}\label{Appendix:log_normal}

\paragraph{Proof of Proposition~\ref{prop:lognormal}.}
For a lognormal random variable $X \sim \mathrm{LogNormal}(\mu,\sigma^{2})$,  
the median is $\mathrm{Med}(X)=\exp(\mu)$, so specifying the median $m$ yields
\[
\mu = \ln(m).
\]

The mean and variance of $X$ are
\[
\mathrm{E}[X] = \exp\!\left(\mu + \tfrac{1}{2}\sigma^{2}\right),
\qquad
\mathrm{Var}(X) = \bigl(\exp(\sigma^{2}) - 1\bigr)\exp\!\left(2\mu + \sigma^{2}\right).
\]
Hence the coefficient of variation is
\[
\begin{aligned}
cv
&= \frac{\sqrt{\mathrm{Var}(X)}}{\mathrm{E}[X]} \\
&= \frac{\sqrt{\bigl(\exp(\sigma^{2}) - 1\bigr)\exp(2\mu + \sigma^{2})}}
        {\exp(\mu + \tfrac{1}{2}\sigma^{2})} \\
&= \sqrt{\exp(\sigma^{2}) - 1}.
\end{aligned}
\]

Equating this expression with the prescribed value of $cv$ gives
\[
cv^{2} = \exp(\sigma^{2}) - 1
\qquad\Longrightarrow\qquad
\sigma^{2} = \ln\!\bigl(1 + cv^{2}\bigr)
\qquad\Longrightarrow\qquad
\sigma = \sqrt{\ln\!\bigl(1 + cv^{2}\bigr)}.
\]
\hfill$\blacksquare$

% ================================ A.2
\subsection{\rev{zg}{}{Traffic Speed Prediction and Error Estimation}}\label{A2}

To generate realistic stochastic prediction errors, we employ the Diffusion Convolutional Recurrent Neural Network (DCRNN) \citep{li2017diffusion} to model spatiotemporal traffic dependencies. Two benchmark datasets are utilized: \textbf{METR-LA} (4 months of highway data from 207 sensors) as a proxy for Mandl’s network, and \textbf{Q-Traffic} (61 days of Beijing 6th Ring Road segments) for the real-world case study. We split both datasets, chronologically,  as recommended by \cite{li2017diffusion} and \cite{lu2022graph}, respectively: We train each prediction model on the first part of the chronologically ordered
data and measure prediction errors on the last part, i.e., the test set.
Table~\ref{tab:prediction_error} reports the one-hour-ahead prediction errors
used in our experiments.

% ----------------- table_predERROR
\begin{table}[htbp!]
\centering
\caption{One-hour-ahead prediction error on benchmark datasets.}
\label{tab:prediction_error}
\begin{tabular}{llcc}
\hline
Dataset & Model & MAE & RMSE \\
\hline
\multirow{2}{*}{METR-LA}
& DCRNN & 3.590 & 7.607 \\
& MTGNN & 3.490 & 7.230 \\
\hline
\multirow{2}{*}{Q-Traffic}
& DCRNN & 3.932 & 6.030 \\
& MTGNN & 2.973 & 4.742 \\
\hline
\end{tabular}
\end{table}
% ----------------- table_predERRORend

The prediction models output traffic speeds. Since the downstream public-transport redesign model operates on travel times, the extracted prediction residuals must be expressed in travel-time space. If a link length $d$ is available, the realized and predicted travel times are $T=d/v$ and $\hat T=d/\hat v$, respectively. 
The corresponding log-relative travel-time residual is therefore $\eta=\log(\hat T)-\log(T)$. For a fixed link length, this is equivalently written as $\eta=\log(v)-\log(\hat v)$, since the length term cancels out algebraically.
% ~\footnote{A short derivation is provided in App.~\ref{app:speed_to_time_log_residual}.} 
%
This equivalence is particularly useful for \textbf{METR-LA}, where reliable segment lengths are not directly available for all sensor pairs. In this case, we compute the travel-time log-relative residual directly from the observed and predicted speeds. For \textbf{Q-Traffic}, where segment lengths are explicitly provided, the same residual can be obtained either after converting speeds into travel times or directly through the equivalent speed-based expression. The resulting residual pool is then used to generate noisy travel-time forecasts as described in App.~\ref{A3}.

% -----------------------------
% -----------------------------
% -----------------------------
% ================================ A.3
\subsection{\rev{zg}{}{Empirical Log-Relative Residual Sampling for Travel-Time Prediction}}\label{A3}

We pool all empirical prediction residuals on the test sets, across links and time steps, to construct a representative empirical error pool. Since the downstream experiments transfer prediction errors across links and networks with different travel-time scales, we sample log-relative residuals rather than absolute residuals. Specifically, for each realized travel time $T$ and its prediction $\hat T$, we define $\eta=\log(\hat T)-\log(T)=\log\left(\hat T/T\right)$.
This log-ratio formulation is appropriate for heteroscedastic prediction errors whose magnitude scales with the level of the predicted variable~\citep{Tofallis2015}; in our setting, this is consistent with traffic scenarios in which the coefficient of variation changes the scale of realized travel times. 
% 适合误差随变量尺度增长的 heteroscedastic 情况
%
Finally, independent samples $\eta_{ij}^{\omega}$ are drawn with replacement from this empirical pool for each link $(i,j)$ and scenario $\omega$. Given the realized travel time $t_{ij}^{\omega}$ generated under the prescribed traffic-variability level, the corresponding additive travel-time error is defined as $\epsilon_{ij}^{\omega}=t_{ij}^{\omega}\left(\exp!\left(\alpha\eta_{ij}^{\omega}\right)-1\right)$. The predicted link travel times~$\hat t_{ij}^{\omega}$, used to make decisions via solving~\eqref{eq:multi_objective}-\eqref{st:veh_budget}, are then computed by adding $\epsilon_{ij}^{\omega}$ to the actual travel time~$t_{ij}^{\omega}$, according to~\eqref{eq:predTRAVEL}. We do this in both Mandl and Beijing scenarios.

Moreover, because errors are sampled independently for each link and scenario, this first-order model abstracts away the spatial and temporal correlation of real forecast errors, such as corridor-level or incident-level common shocks. The reported results should therefore be read as an evaluation under sampled independent empirical residual scenarios, and correlation-aware forecast-error modeling remains an important extension.

\begin{equation}
\hat t_{ij}^{\omega} = t_{ij}^{\omega}
+
\epsilon_{ij}^{\omega}, \
\epsilon_{ij}^{\omega} =
t_{ij}^{\omega}
\left(
\exp!\left(\alpha\eta_{ij}^{\omega}\right)-1
\right)
\label{eq:predTRAVEL_appendix}
\end{equation}

% ----------------------

\paragraph{Diagnostic check on traffic variability.}
For the Mandl experiments, we further report the realized-to-nominal travel-time ratio
$R_{ij}^{\omega}=t_{ij}^{\omega}/t_{ij}^{0}$, where $t_{ij}^{0}$ is the nominal link travel time and $t_{ij}^{\omega}$ is the realized travel time under scenario $\omega$. Table~\ref{tab:travel_time_ratio} summarizes the distribution of $R$ under different traffic-variability levels. The results show that increasing \textit{cv} mainly amplifies the dispersion and upper tail of realized travel times. Small-\textit{cv} scenarios remain close to nominal conditions, whereas large-\textit{cv} scenarios generate severe congested realizations, as reflected by the increasing standard deviation and 95th percentile of $R$. In these diagnostics, the prediction-error scale is fixed at $\alpha=1.0$, corresponding to the empirical error magnitude extracted from the trained prediction model.

% -------------------- table_dig
\begin{table}[htbp!]
\centering
\caption{Distribution of realized-to-nominal travel-time ratios in the Mandl experiments.}
\label{tab:travel_time_ratio}
\begin{tabular}{lccc}
\hline
\textit{cv} & Mean($R$) & Std($R$) & P95($R$) \\
\hline
0.00 & 1.000 & 0.000 & 1.000 \\
0.05 & 0.988 & 0.047 & 1.076 \\
0.25 & 0.961 & 0.228 & 1.436 \\
0.50 & 0.973 & 0.456 & 2.001 \\
1.00 & 1.090 & 0.949 & 3.396 \\
\hline
\end{tabular}
\end{table}
\section{Algorithmic Framework and Evaluation}
\label{sec:algorithmic-framework-and-evaluation}

\subsection{Parameter values for the proposed Algorithm~\ref{alg:NSGA3}}
When computing adaptive network $\pazocal G_\text{PT}$ using the solution method described in $\S$~\ref{sec:alg}, we set the algorithmic parameters introduced therein as follows. 
The minimum path-overlap threshold is set to $J_{\text{min}}=0.5$ (\S{}\ref{sec:init_population}), the crossover and mutation probabilities to $P_c=0.85$ (\S{}\ref{sec:crossover}) and $P_m=0.35$ (\S{}\ref{sec:mutation}), and the number of reference points to $N_{\text{ref}}=5^3=125$ (\S{}\ref{sec:env_selection}), respectively.

% ------------------------------------------------E.1
\subsection{Performance Metrics for Multi-Objective Optimization}\label{app:metrics}
We evaluate algorithm performance using three standard metrics for multi-objective optimization: Hypervolume (HV), Inverted Generational Distance (IGD), and Generational Distance (GD). 
The HV metric quantifies the volume of the objective space dominated by the Pareto solutions:
\begin{equation}
% \text{HV}(X, P) = \bigcup_{x \in X} v(x, P)
\text{HV}(X, P) = \Lambda \left( \bigcup_{x \in X} v(x, P) \right)
\end{equation}
where $X$ is the approximation set, $P$ is the reference point, and $v(x, P)$ is the hypervolume dominated by solution $x$ relative to $P$. Higher HV values indicate better coverage of the objective space. 
To ensure a fair comparison across algorithms, we report the HV metric calculated relative to a unified reference point. 
This reference point is constructed by collecting all solutions generated by all algorithms for each objective and selecting the worst (largest) value observed for each objective. 
% The same reference point is then used for all methods.
% Under this uniform evaluation standard, NSGA-III attains the largest HV$^\ast$, indicating the most favorable Pareto approximation among the compared methods.
%
The IGD metric measures the closeness of the approximation set to the true Pareto front:
\begin{equation}
\text{IGD}(X, \pazocal{P}^*) = \frac{\sum_{x^* \in \pazocal{P}^*} d(x^*, X)}{|\pazocal{P}^*|}
\end{equation}
where $\pazocal{P}^*$ is a set of reference points on the true Pareto front, $d(x^*, X)$ is the minimum distance from $x^*$ to any point in $X$, and $|\pazocal{P}^*|$ is the cardinality of $\pazocal{P}^*$. Lower IGD values indicate better convergence to the true Pareto front.
The GD metric complements IGD by measuring how far the obtained solutions are from the true Pareto front:
\begin{equation}
\text{GD}(X, \pazocal{P}^*) = \frac{\sum_{x \in X} d(x, \pazocal{P}^*)}{|X|}
\end{equation}
where $d(x, \pazocal{P}^*)$ is the minimum distance from a solution $x$ to the true Pareto front $\pazocal{P}^*$, and $|X|$ is the cardinality of the approximation set. Lower GD values indicate better convergence of the approximation set to the true Pareto front.
\rev{zg}{}{As the true Pareto front $\pazocal{P}^*$ is unknown for the instances considered, we approximate it, following standard practice in multi-objective optimization, by the non-dominated subset of the union of all solutions produced by all compared algorithms across all runs at a given variability level. The HV reference point is constructed analogously from the worst objective values observed across all algorithms (as described above), so that HV, IGD, and GD are all evaluated against a common, algorithm-agnostic reference.}

% ------------------------------------------------------
% % ++++++++++++++++++++ 6.1.1 - Static Network Optimization
% ------------------------------------------------------
% ------------------------------------------------E.1
\subsection{Implementation Details of the NSGA-III Environmental Selection}\label{A.Evolution}

The algorithm adopts the reference-point-based environmental selection mechanism
of NSGA-III~\citep{deb2014nsga3} to maintain diversity across the objectives~\eqref{eq:multi_objective}.
Since all objectives are formulated as minimization problems, normalization is
performed on the combined population $R_t = P_t \cup Q_t$.
For each objective $j \in \{1,2,3\}$, objective-wise lower and upper bounds are
computed as
$z_j^{\min} = \min_{x \in R_t} Z_j(x)$ and
$z_j^{\max} = \max_{x \in R_t} Z_j(x)$,
and each objective is linearly scaled to $[0,1]$ accordingly.

A set of uniformly distributed reference points on the unit simplex is generated following the simplex-lattice design in~\cite{deb2014nsga3}.
% [Revision for aa: Which population? & No ellipses & Define Dominated]
As in standard multiobjective optimization, we say that solution $x$ dominates $y$ if $x$ is no worse than $y$ in all objectives and strictly better in at least one. We employ non-dominated sorting, which partitions the \textbf{combined population} $R_t = P_t \cup Q_t$ (comprising both parent and offspring individuals) into a sequence of Pareto fronts $F = \{F_1, F_2, \dots, F_l\}$.
% [Revision for aa: Exact list of partitions logic]
Here, $F_1$ contains the globally non-dominated solutions
% \footnote{A solution $x$ dominates $y$ if $x$ is no worse than $y$ in all objectives and strictly better in at least one.}
$F_2$ contains solutions dominated exclusively by members of $F_1$, and $F_l$ represents the \textit{critical front} where the cumulative count of solutions first meets or exceeds the population capacity $\textit{PopSize}$.
Each solution is associated with its nearest reference point by computing perpendicular distances from the solution to the reference lines in the normalized space. 
The environmental selection process (line~\ref{line:env_select} of Alg.~\ref{alg:NSGA3}) adds complete Pareto fronts $F_1$ through $F_{l-1}$ sequentially; for the last front $F_l$, slots are filled by prioritizing solutions associated with reference points that have fewer assigned solutions to ensure diversity.

% ===================== algRepair
\begin{algorithm}[htbp]
\small
\caption{Feasibility-Oriented Path Completion for Line Repair}
\label{alg:bfs_repair}
\begin{algorithmic}[1]
    \renewcommand{\algorithmicrequire}{\textbf{Input:}}
    \renewcommand{\algorithmicensure}{\textbf{Output:}}
    
    \REQUIRE Substrate graph $\pazocal{G}=(\pazocal{N}, \pazocal{E})$, endpoints $(u, v)$, visited node set $\pazocal V_{\mathrm{used}}^\ell$
    \ENSURE Connecting path $\pazocal{P}_{u \to v}$ or $\emptyset$
    
    \STATE $\mathit{visited} \gets \{u\} \cup \pazocal V_{\mathrm{used}}^\ell$, $\mathit{queue} \gets \langle u \rangle$, $\mathit{parent}(\cdot) \gets \emptyset$
    
    \WHILE{$\mathit{queue}$ is not empty}
        \STATE $i \gets \mathit{queue}.\text{pop\_front}()$
        \IF{$i = v$}
            \STATE \textbf{break}
        \ENDIF
        \FOR{each neighbor $j$ of $i$ in $\pazocal{G}$}
            \IF{$j \notin \mathit{visited}$}
                \STATE $\mathit{visited} \gets \mathit{visited} \cup \{j\}$
                \STATE $\mathit{parent}(j) \gets i$
                \STATE $\mathit{queue}.\text{push\_back}(j)$
            \ENDIF
        \ENDFOR
    \ENDWHILE
    
    \IF{$v \notin \mathit{visited}$}
        \RETURN $\emptyset$
    \ENDIF
    \STATE $\pazocal{P}_{u \to v} \gets \text{Backtrack from } v \text{ to } u \text{ using } \mathit{parent}(\cdot)$
    \RETURN $\pazocal{P}_{u \to v}$
\end{algorithmic}
\end{algorithm}
\subsection{Algorithm Performance Evaluation}

Algorithmic performance is evaluated via 
the Hypervolume (HV), the Inverted Generational Distance (IGD), and the Generational Distance (GD), which are standard performance metrics in multi-objective evolutionary optimization
as surveyed in~\citep{riquelme2015performance}.\footnote{Formal definitions and computational details are provided in App.~\ref{app:metrics}.}
The Hypervolume (HV) evaluates both convergence and diversity by measuring the dominated volume in the objective space, where higher values indicate better performance ($\uparrow$). In contrast, GD and IGD are distance-based indicators that assess convergence to a reference Pareto front, with lower values indicating better convergence ($\downarrow$).

To evaluate the effectiveness of NSGA-III for \rev{zg}{real-time}{forecast-triggered} PT redesign, we compare it with NSGA-II~\citep{deb2002fast} and hybrid multi-objective algorithms that combine non-dominated sorting with variable neighborhood descent, based on Mandl’s Swiss network.
Specifically, we consider two Non-dominated Sorting Variable Neighborhood Descent (NSVND) variants inspired by the integration of NSGA-II with variable neighborhood search strategies~\citep{hansen2001variable}.
The first variant (NSVND-II) follows the NSGA-II framework while replacing crossover and mutation with variable neighborhood descent, whereas the second variant (NSVND-III) further incorporates an adaptive neighborhood selection mechanism to enhance exploration of the Pareto frontier.

% ========================== tableHVIGD©
\begin{table}[htbp]
\centering
\footnotesize
\setlength{\tabcolsep}{4.5pt}
\renewcommand{\arraystretch}{1.15}
\begin{threeparttable}
\begin{tabular}{@{}l*{4}{c}@{}}
\toprule
\textbf{Performance Metric} & \textbf{NSGA-II} & \textbf{NSGA-III} & \textbf{NSVND-II} & \textbf{NSVND-III} \\
% \midrule
% Hypervolume (HV, higher is better) 
% & $3.71\times10^9$ 
% & $7.09\times10^9$ 
% & $7.28\times10^9$ 
% & \textbf{$7.36\times10^9$} \\
% HV Improvement (\%)\tnote{a} 
% & -- 
% & 91.01 
% & 96.23 
% & \textbf{98.38} \\
\midrule
Inverted Generational Distance (IGD, lower is better) 
& 347.64 
& \textbf{22.63} 
& 38.90 
& 34.80 \\
IGD Change\tnote{a} 
& -- 
& \textbf{-325.01} 
& -308.74 
& -312.84 \\
\midrule
Generational Distance (GD, lower is better) 
& 75.36 
& \textbf{17.17} 
& 24.40 
& 22.80 \\
GD Change\tnote{a} 
& -- 
& \textbf{-58.19} 
& -50.96 
& -52.56 \\
\midrule
Hypervolume (HV, higher is better) \tnote{b} 
& $30.52\times10^9$ 
& \textbf{$31.49\times10^9$} 
& $31.23\times10^9$ 
& $31.31\times10^9$ \\
HV Change\tnote{a} 
& -- 
& \textbf{0.97$\times10^9$} 
& 0.71$\times10^9$ 
& 0.79$\times10^9$ \\
\bottomrule
\end{tabular}

\begin{tablenotes}\footnotesize
\item[a] Improvement/Change are computed relative to NSGA-II (baseline).
\item[b] HV uses a single fixed reference point shared by all algorithms.
\end{tablenotes}
\end{threeparttable}
\caption{Multi-objective algorithm comparison under noise-free travel times ($cv=0$).}
\label{tab:algorithm_comparison_simple}
\end{table}
% ========================== tableHVIGDEND

Table~\ref{tab:algorithm_comparison_simple} presents an algorithm performance comparison for adaptive network optimization assuming no travel time deviation from the nominal ones (\textit{cv}=0).
% ~\footnote{
Although no adaptive redesign is required when $cv=0$, this setting provides a noise–free benchmark for comparing optimization algorithms.
% } 
% While NSVND-series algorithms demonstrate superior hypervolume metrics, they exhibit significantly inferior convergence as evidenced by their poor IGD and GD values. In contrast, 
NSGA-III maintains an optimal balance between solution diversity (high hypervolume) and convergence precision.
%
% To enable a fair cross-algorithm comparison, we additionally report the hypervolume computed with a common fixed reference point (HV$^\ast$).%
% \footnote{The common reference point is constructed by first collecting all solutions generated by all algorithms and then taking, for each objective, the worst (largest) value observed. This single reference point is shared by all algorithms, enabling their hypervolume values to be directly comparable.}
% Under this uniform evaluation standard, NSGA-III attains the largest HV$^\ast$, indicating the most favorable Pareto approximation among the compared methods.

\section{Network Performance Metrics}\label{A.PerformanceMetrics}

This appendix details the performance indicators used in $\S$~\ref{sec:performMETRICS}.

For each realization $\omega \in \Omega$ of stochastic travel times, we define the following performance metrics: \textbf{20\% Improvement Share (IS$_{20}$)}, \textbf{Operational Cost Reduction (OCR)}, and \textbf{Line Overlap (LO)}. Since these metrics vary with the realized scenario $\omega$, they are treated as random variables. In our numerical results, we report their empirical averages (and standard deviations) over 30 independent realizations.

\begin{itemize}
    \item \textbf{20\% improv. share (\%):} The proportion of OD pairs experiencing a travel time reduction of at least 20\% under realization $\omega$.
    \begin{align}
        \label{eq:metricIS20}
        \text{IS$_{20}$}(\omega) = \frac{\sum_{k \in \pazocal{K}} \mathbb{I} \left( \frac{t_k^0 - t_k^*(\omega)}{t_k^0} \geq 20\% \right)}{|\pazocal{K}|} \times 100\%, \text{ where } \mathbb{I}(\cdot) \text{ is the indicator function.}
    \end{align}
    Here $t_k^0$ and $t_k^*(\omega)$ denote the original and adaptive travel times for OD pair $k$, respectively. Since both depend on random travel-time realizations, $\text{IS$_{20}$}$ is a random variable.

    \item \textbf{Operational Cost Reduction (OCR)}
\begin{equation}
\label{eq:metricOCR}
\operatorname{OCR}(\omega)
=
\frac{
    \sum_{\ell \in \pazocal{L}} n^\ell t^\ell
    -
    \sum_{\ell \in \pazocal{L}} n^{\ell,*}(\omega) t^{\ell,*}(\omega)
}{
    \sum_{\ell \in \pazocal{L}} n^\ell t^\ell
}
\times 100\%.
\end{equation}
    Here $t^\ell$ and $t^{\ell,*}(\omega)$ are the realized travel times of line $\ell$ in the original and adaptive networks, respectively,
    and $n^\ell$ and $n^{\ell,*}$ are the corresponding vehicle allocations.
    Because $t^\ell$ and $t^{\ell,*}(\omega)$ vary with the realization, OCR is also a random variable.

    \item \textbf{Line Overlap (LO)}
    To quantify structural similarity between the original and adaptive networks, we use the Jaccard index:
\begin{equation}
\label{eq:metricRO}
\operatorname{LO}(\omega)
=
\frac{
    \left| \pazocal{E}_{\mathrm{PT}}
    \cap
    \pazocal{E}_{\mathrm{PT}}^{*}(\omega) \right|
}{
    \left| \pazocal{E}_{\mathrm{PT}}
    \cup
    \pazocal{E}_{\mathrm{PT}}^{*}(\omega) \right|
}
\times 100\%.
\end{equation}
    Here $\pazocal{E}_{\mathrm{PT}}$ and $\pazocal{E}_{\mathrm{PT}}^*(\omega)$ denote the realized edge sets of the original and adaptive networks, respectively.
    Since these sets depend on the travel-time realization, LO is likewise a random variable.

\end{itemize}

The numerical results therefore report the \emph{empirical means} $\pm$ \emph{standard deviations} of IS$_{20}$, OCR, and LO over the 30 realizations $\omega \in \Omega$, providing a stochastic assessment of adaptive-network performance.

% -------------- App.D.1
\subsection{Forecast-error Sensitivity Analysis}
\label{app:forecastValue}

This appendix provides the formal basis and the full results for the forecast-error sensitivity analysis summarized in the main text. App.~\ref{app:forecastValueMetrics} defines the two hypervolume-based summaries used to quantify \rev{aa}{}{the impact of the} forecast \rev{aa}{value}{error in the optimization output} (the \emph{forecast-value headroom} and the \emph{forecast value retained}). App.~\ref{app:forecast_sensitivity} then applies these metrics to the Mandl instance and reports the detailed results across variability levels and error scales~$\alpha$.

We deliberately grade the entire front rather than any single selected design: the optimizer returns a set of non-dominated trade-offs, and the specific compromise that gets deployed (e.g., the ideal-point solution of $\S$\ref{sec:compared-PT-networks}) is chosen by the operator after the front is produced, according to priorities the algorithm itself does not fix. Tying the sensitivity analysis to one such pre-selected point would conflate two different effects — forecast quality and the arbitrary choice of operating point — since the operator could simply pick a different point on the front to offset a degraded solution, masking a genuine impact of the forecast accuracy on the quality of the solutions (or manufacturing a spurious one). Hypervolume sidesteps this confound by scoring the whole menu of options the operator has to choose from in a single scalar: a drop in $HV$ reflects an actual shrinkage of what is available to the operator, not an artifact of which point happens to be selected.

\subsubsection{Forecast-value metrics}
\label{app:forecastValueMetrics}

We quantify whether the forecast helps the redesign, and by how much. Holding the realized \rev{aa}{scenario}{link travel times} $\{t_{ij}^{\omega}\}$ fixed, we vary \emph{only} the forecast link times $\hat{t}_{ij}^{\omega}(\alpha)$ handed to the optimizer (Eq.~\eqref{eq:predTRAVEL_alpha}) and score the returned Pareto front on the true realized times $t_{ij}^{\omega}$. Since only the error scale $\alpha$ changes, any difference is attributable to forecast quality alone. Each front is graded by its hypervolume $HV$, the volume it dominates in the space of the three performance metrics $(\mathrm{IS}_{20},\mathrm{OCR},\mathrm{LO})$ defined in \S\ref{sec:performMETRICS}, measured from the single fixed reference point of App.~\ref{app:metrics}; a larger $HV$ is better\rev{aa}{ and needs no weighting of the three metrics}{}. For each seed $s\in\{1,\dots,S\}$ (with $S=30$ random seeds) and each variability level $\textit{cv}$ we compute three fronts:
\begin{itemize}
    \item \emph{No forecast}: \rev{aa}{}{the optimizer uses the} nominal link times $t_{ij}^{0}$\rev{aa}{, giving the \emph{floor}}{We call the corresponding hypervolume} $HV^{\mathrm{NF}}_{s,\textit{cv}}$ \rev{aa}{(lower bound)}{\emph{floor}, since it corresponds to the lower bound}.

    \item \emph{Perfect foresight} ($\alpha=0$): the optimizer uses the true realized times $t_{ij}^{\omega}=\hat{t}_{ij}^{\omega}(0)$. We denote the corresponding hypervolume by $HV^{\mathrm{PF}}_{s,\textit{cv}}\equiv HV_{s,\textit{cv},0}$ and call it the \emph{ceiling}, since it represents the upper bound.

    \item \emph{Error-$\alpha$ forecast}: the forecast $\hat{t}_{ij}^{\omega}(\alpha)$ of Eq.~\eqref{eq:predTRAVEL_alpha} ($\alpha=1$ the empirical DCRNN level, $\alpha>1$ a stress test), giving the intermediate hypervolume $HV_{s,\textit{cv},\alpha}$.
\end{itemize}
Because a seed fixes both the realized scenario and the sampled errors $\eta$, each $HV_{s,\textit{cv},\alpha}$ is a single realization, and all reported 95\% confidence intervals reflect variation across the $S$ seeds. We measure the climb from floor to ceiling in two steps.

\paragraph{Step 1: is there anything to gain?} The \emph{forecast-value headroom} $H_{s,\textit{cv}}=\bigl(HV^{\mathrm{PF}}_{s,\textit{cv}}-HV^{\mathrm{NF}}_{s,\textit{cv}}\bigr)/HV^{\mathrm{PF}}_{s,\textit{cv}}$ is the floor-to-ceiling gap as a fraction of the ceiling: $H\approx0$ means even a perfect forecast barely beats nominal times, whereas a large $H$ signals real value to capture. We report the mean $\bar H_c$ across seeds with a 95\% confidence interval.

\paragraph{Step 2: how much do we actually capture?} The \emph{forecast value retained} $R_{s,\textit{cv},\alpha}=\bigl(HV_{s,\textit{cv},\alpha}-HV^{\mathrm{NF}}_{s,\textit{cv}}\bigr)/\bigl(HV^{\mathrm{PF}}_{s,\textit{cv}}-HV^{\mathrm{NF}}_{s,\textit{cv}}\bigr)$ places an error-$\alpha$ forecast on that same scale: $R=1$ matches perfect foresight, $R=0$ matches no forecast, and $R<0$ is worse than no forecast. Since this gap can be tiny for a single seed (making per-seed ratios erratic, unlike $H_{s,\textit{cv}}$, whose denominator $HV^{\mathrm{PF}}_{s,\textit{cv}}$ stays bounded away from zero), we aggregate it as a ratio of averages,
\begin{equation}
\bar R_{c,\alpha}=\frac{\tfrac1S\sum_s\big(HV_{s,\textit{cv},\alpha}-HV^{\mathrm{NF}}_{s,\textit{cv}}\big)}{\tfrac1S\sum_s\big(HV^{\mathrm{PF}}_{s,\textit{cv}}-HV^{\mathrm{NF}}_{s,\textit{cv}}\big)},
\end{equation}
with a 95\% confidence interval and $\bar R_{c,0}=1$ by construction. At the lowest variability ($cv=0.05$) the floor and ceiling are statistically indistinguishable (Fig.~\ref{fig:forecast_headroom}), so $R$ is undefined and we omit this level from Fig.~\ref{fig:forecast_retained}.

\subsubsection{\rev{zg}{}{Forecast-value results}}\label{app:forecast_sensitivity}

Fig.~\ref{fig:forecast_value_hv} reports the forecast-value metrics of App.~\ref{app:forecastValueMetrics} for the 4-line Mandl instance. By construction, these metrics separate two otherwise-conflated effects: the redesign opportunity created by travel-time variability (the \emph{headroom}), and the accuracy with which the forecast exploits it (the \emph{retained value}). Fig.~\ref{fig:forecast_value_hv}(a) shows that the headroom is negligible at $cv=0.05$ and grows with variability, reaching roughly half of the perfect-foresight hypervolume at $cv=1.0$: a forecast is valuable primarily when variability creates a genuine redesign opportunity, whereas in stable conditions even an oracle adds little because the nominal times $t_{ij}^{0}$ already describe the operating state. Fig.~\ref{fig:forecast_value_hv}(b) shows that the empirical forecast ($\alpha=1$) retains essentially all of this available value, and moderate amplification ($\alpha=1.5$) causes no visible collapse, consistent with the intended robustness claim. \textbf{The method's gains are therefore driven mainly by the redesign opportunities that variability creates, not by a dependence on high-precision forecasts.}

% -------------------------- Figure: Forecast-value headroom + retained
\begin{center}
\captionsetup{type=figure}
\centering
\begin{subfigure}[t]{0.49\textwidth}
\centering
\includegraphics[width=\textwidth]{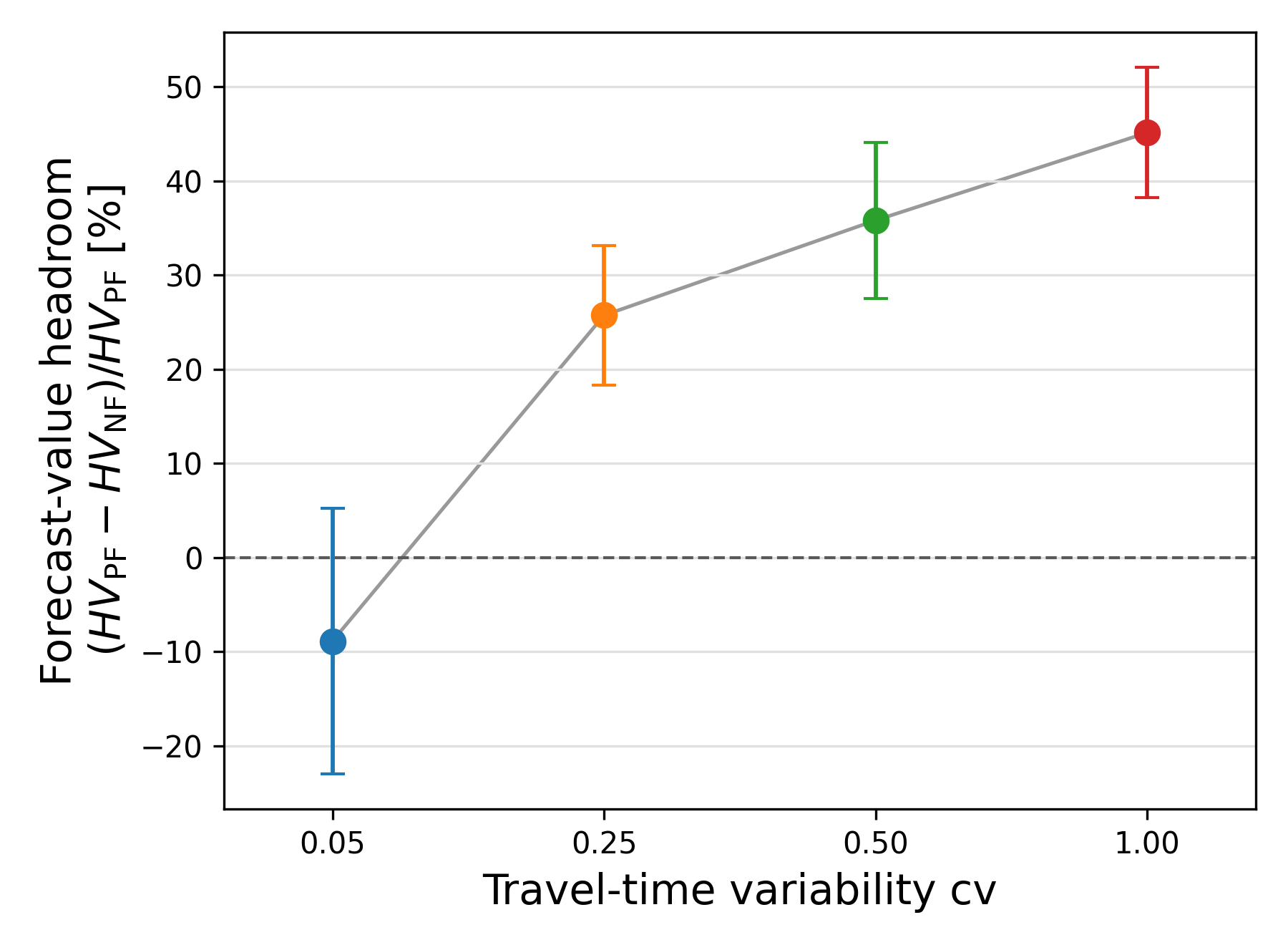}
\caption{\rev{zg}{}{\textbf{Forecast-value headroom $\bar H_c$ vs.\ $cv$} (mean over 30 seeds; 95\% CI bars). Headroom is the share of the perfect-foresight hypervolume attributable to using any forecast.}}
\label{fig:forecast_headroom}
\end{subfigure}
\hfill
\begin{subfigure}[t]{0.49\textwidth}
\centering
\includegraphics[width=\textwidth]{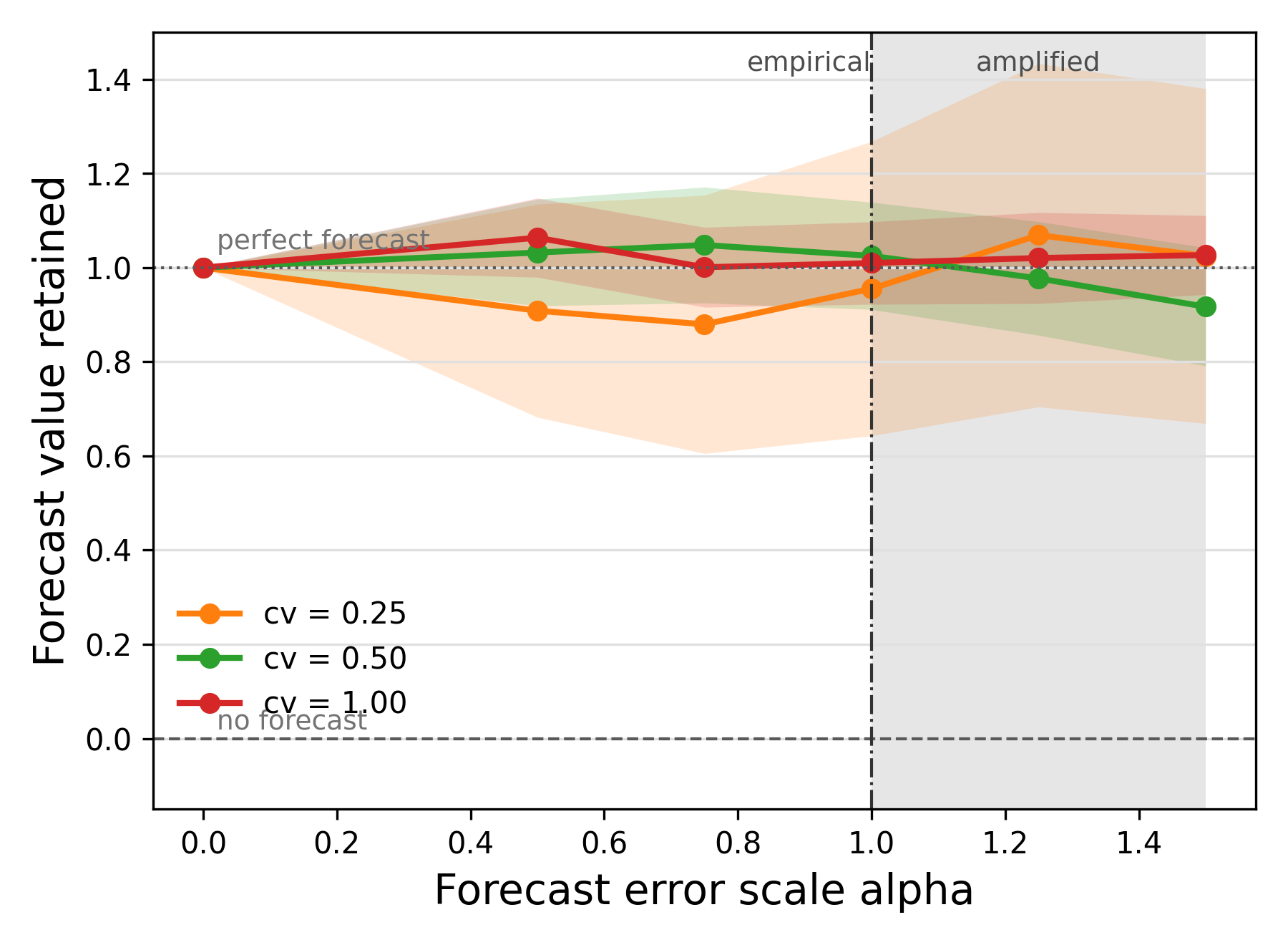}
\caption{\rev{zg}{}{\textbf{Forecast value retained $\bar R_{c,\alpha}$ vs.\ error scale $\alpha$} (ratio-of-means; shaded 95\% CI). $\alpha=1$ is the empirical DCRNN error; $\alpha>1$ is a stress test.}}
\label{fig:forecast_retained}
\end{subfigure}
\caption{\rev{zg}{}{Hypervolume-based forecast value on the Mandl redesign (mean over 30 seeds): (a) headroom available from any forecast as variability grows; (b) fraction of headroom retained as forecast error is amplified. The method retains essentially all available value at the empirical error level and degrades only under severe amplification.}}
\label{fig:forecast_value_hv}
\end{center}

\subsection{\rev{zg}{}{Exceedance Probability of OD-Level Improvement Rates by Bus-Line Quantity}}\label{App:D3_ImprovementICDF}

\rev{aa}{}{For each variability level \(\textit{cv}=0.05, 0.25, 0.5, 1\) and each random seed \(\omega \in \Omega\), we compute the empirical complementary cumulative distribution function (CCDF) of the percentage travel-time improvement over the set of OD pairs. Thus, for a threshold \(x\), the CCDF gives the fraction of OD pairs whose percentage improvement is larger than \(x\). The solid lines in Fig.~\ref{fig:appendix_improvement_exceedance_all_routes}  report the pointwise average of these empirical CCDFs over the 30 seeds.
The figure shows that proactive redesign becomes beneficial only under large travel-time variability, namely for ($cv \geq 0.5$). For lower variability levels, it worsens travel times for a larger fraction of OD pairs than it improves.
}

\rev{aa}{Fig.~\ref{fig:appendix_improvement_exceedance_all_routes} reports the trimmed-mean exceedance curves of OD-level improvement rates. For each seed $\omega$, the OD-level improvement rates are first sorted to form a per-seed quantile curve $Q_\omega(p)$, where $p$ denotes the exceedance probability. The extreme seed curves are then removed and the remaining per-seed quantile curves are averaged pointwise. Thus, the plotted line is a trimmed mean of per-seed quantile curves, rather than a curve computed from pooled OD pairs across all seeds. Across line quantities, the curves make it possible to compare both the intensity of improvements and the share of OD pairs receiving a given minimum percentage reduction under different variability levels.}{}

\begin{center}
\captionsetup{type=figure}
    \centering
    \begin{subfigure}[t]{0.47\textwidth}
        \centering
        \includegraphics[width=\textwidth]{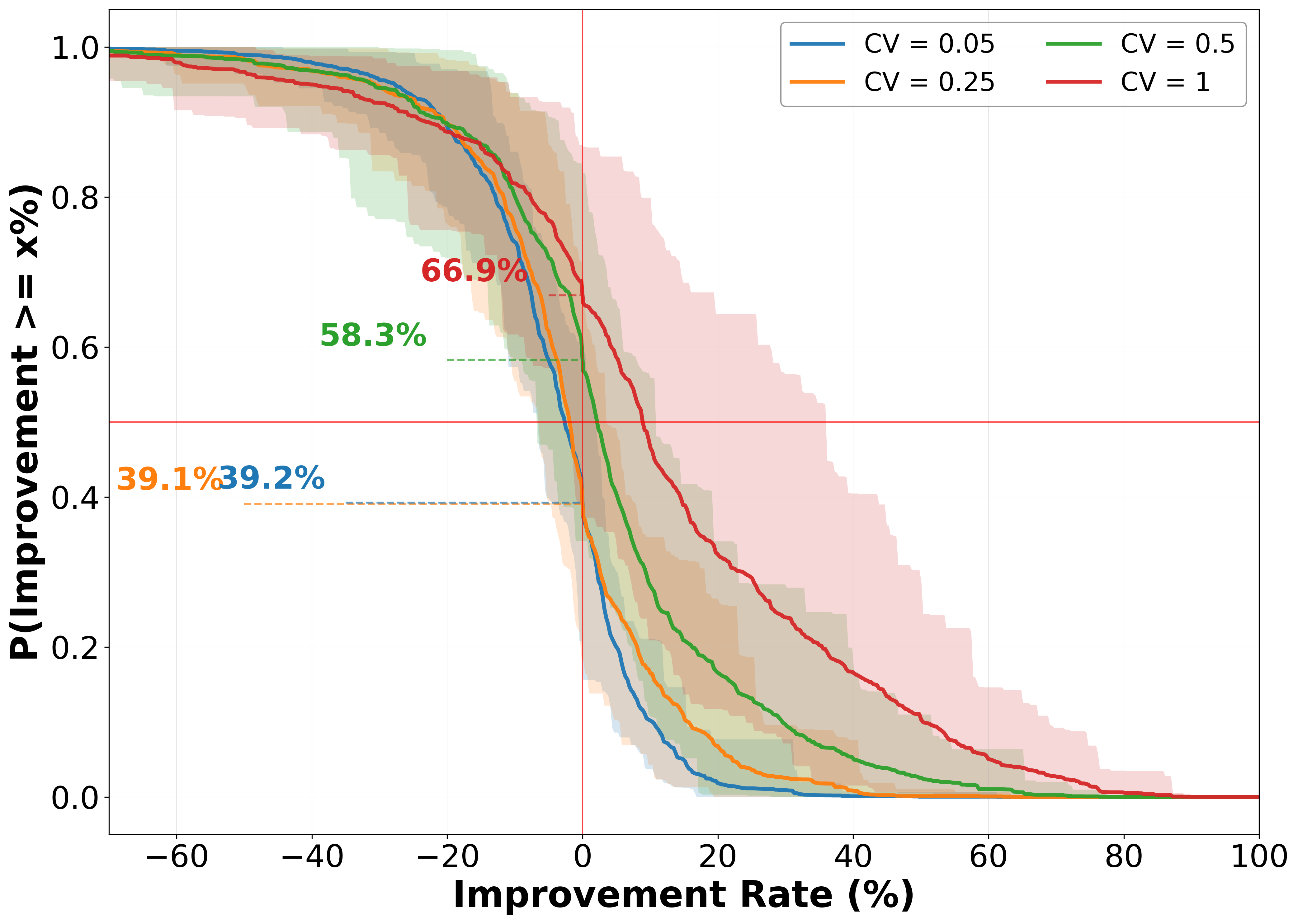}
        \caption{4 lines}
        \label{fig:appendix_improvement_exceedance_4routes}
    \end{subfigure}
    \hfill
    \begin{subfigure}[t]{0.47\textwidth}
        \centering
        \includegraphics[width=\textwidth]{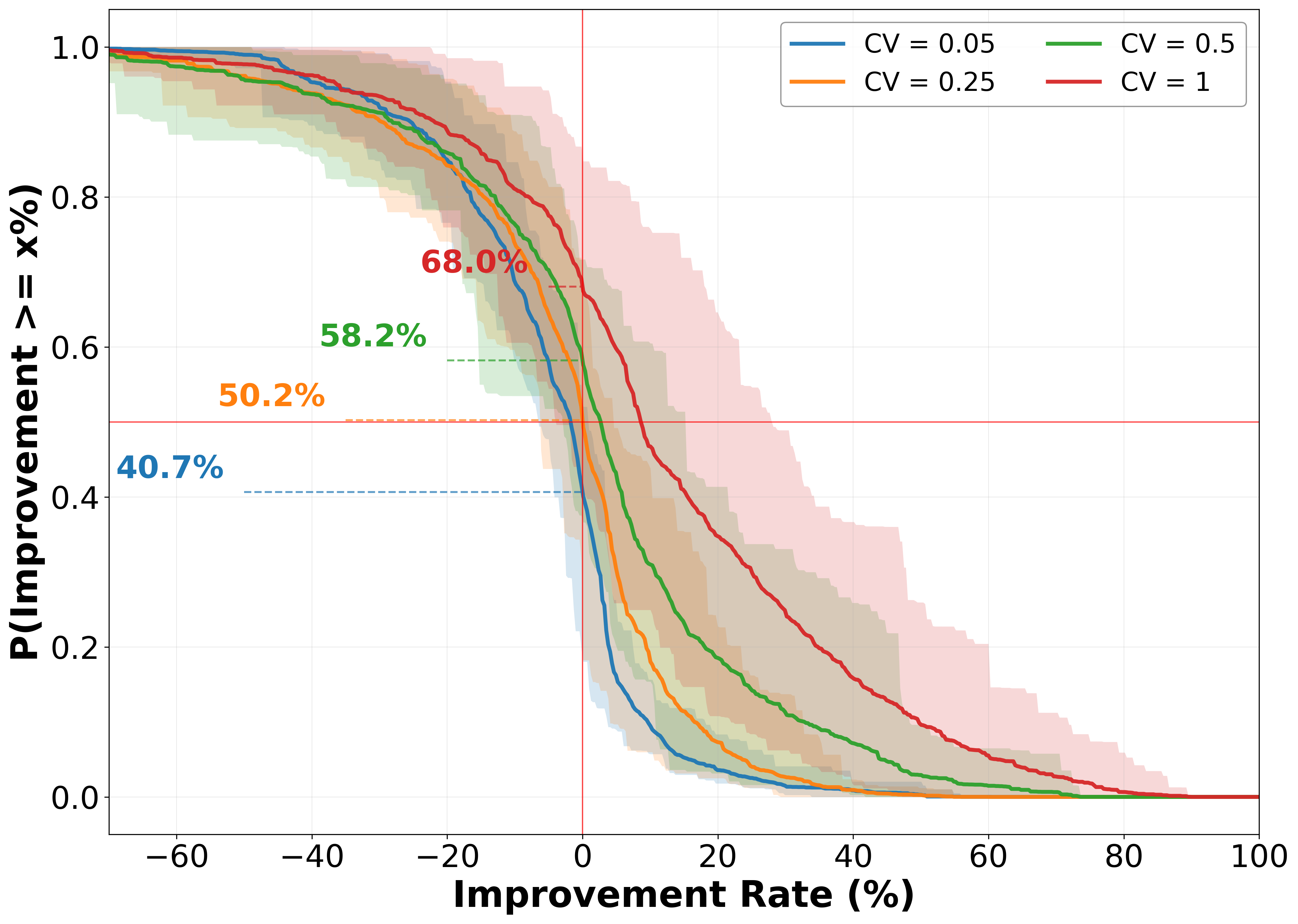}
        \caption{6 lines}
        \label{fig:appendix_improvement_exceedance_6routes}
    \end{subfigure}
    
    \medskip
    \begin{subfigure}[t]{0.47\textwidth}
        \centering
        \includegraphics[width=\textwidth]{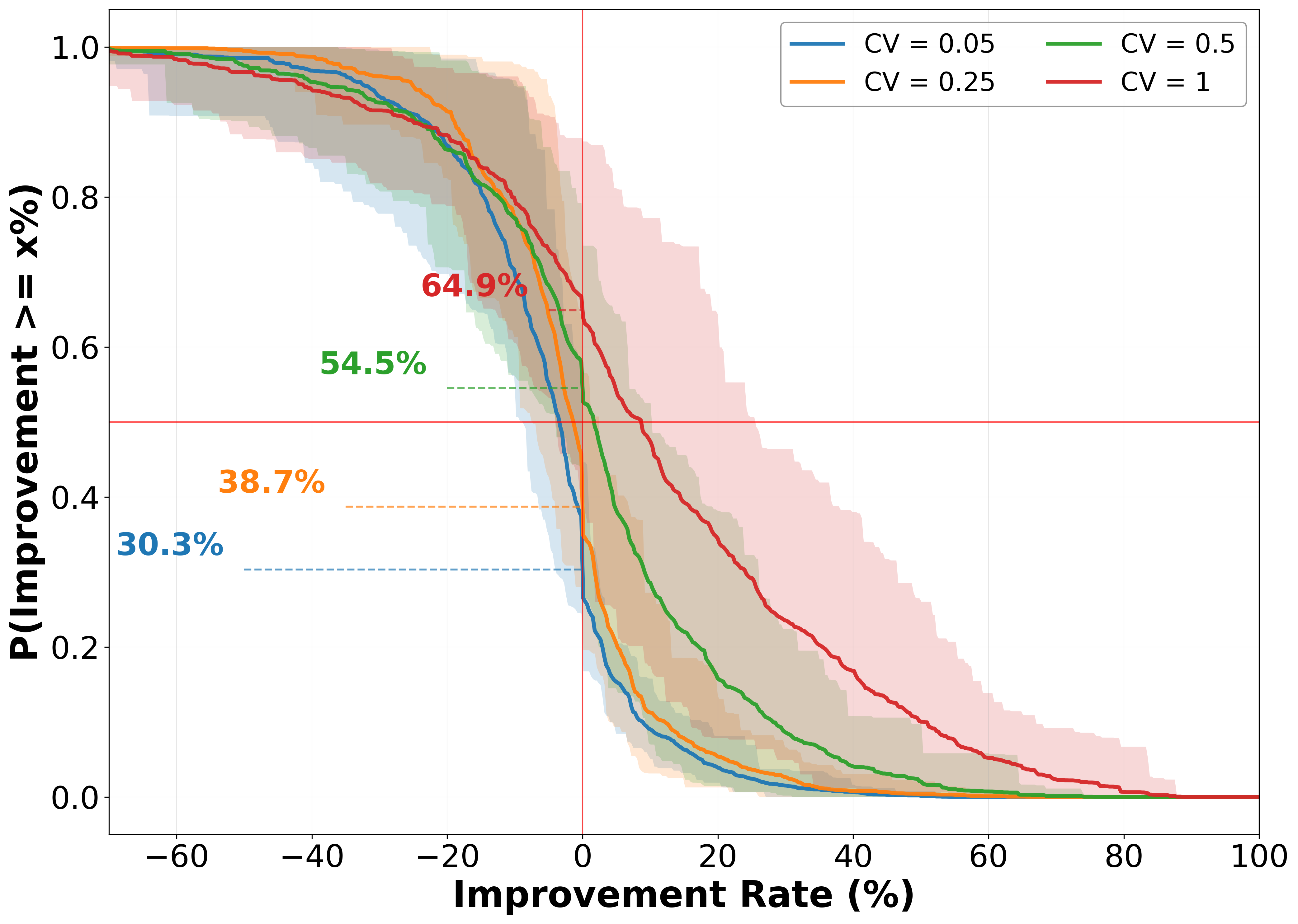}
        \caption{8 lines}
        \label{fig:appendix_improvement_exceedance_8routes}
    \end{subfigure}
    \hfill
    \begin{subfigure}[t]{0.47\textwidth}
        \centering
        \includegraphics[width=\textwidth]{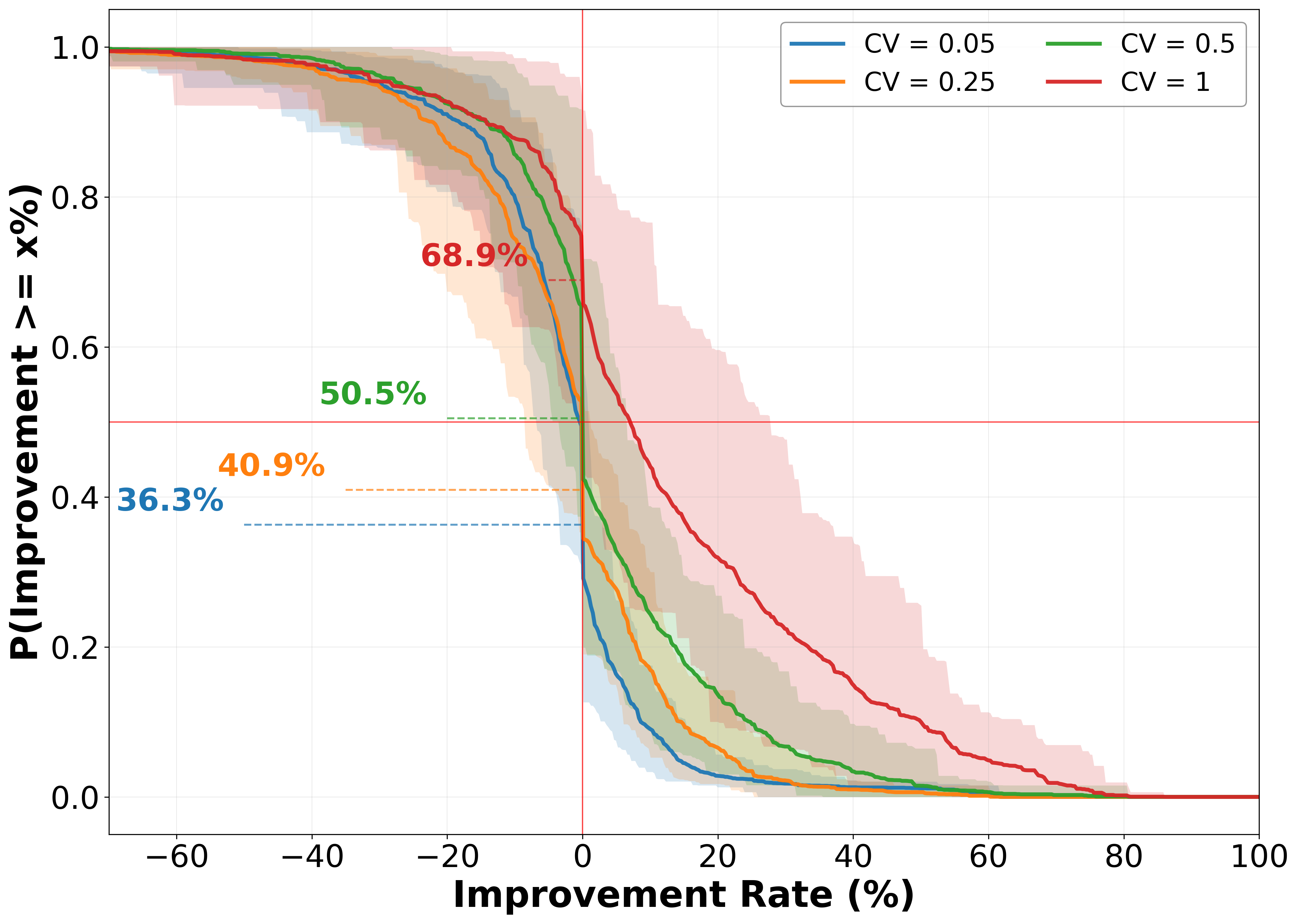}
        \caption{10 lines}
        \label{fig:appendix_improvement_exceedance_10routes}
    \end{subfigure}
    
    \medskip
    \begin{subfigure}[t]{0.47\textwidth}
        \centering
        \includegraphics[width=\textwidth]{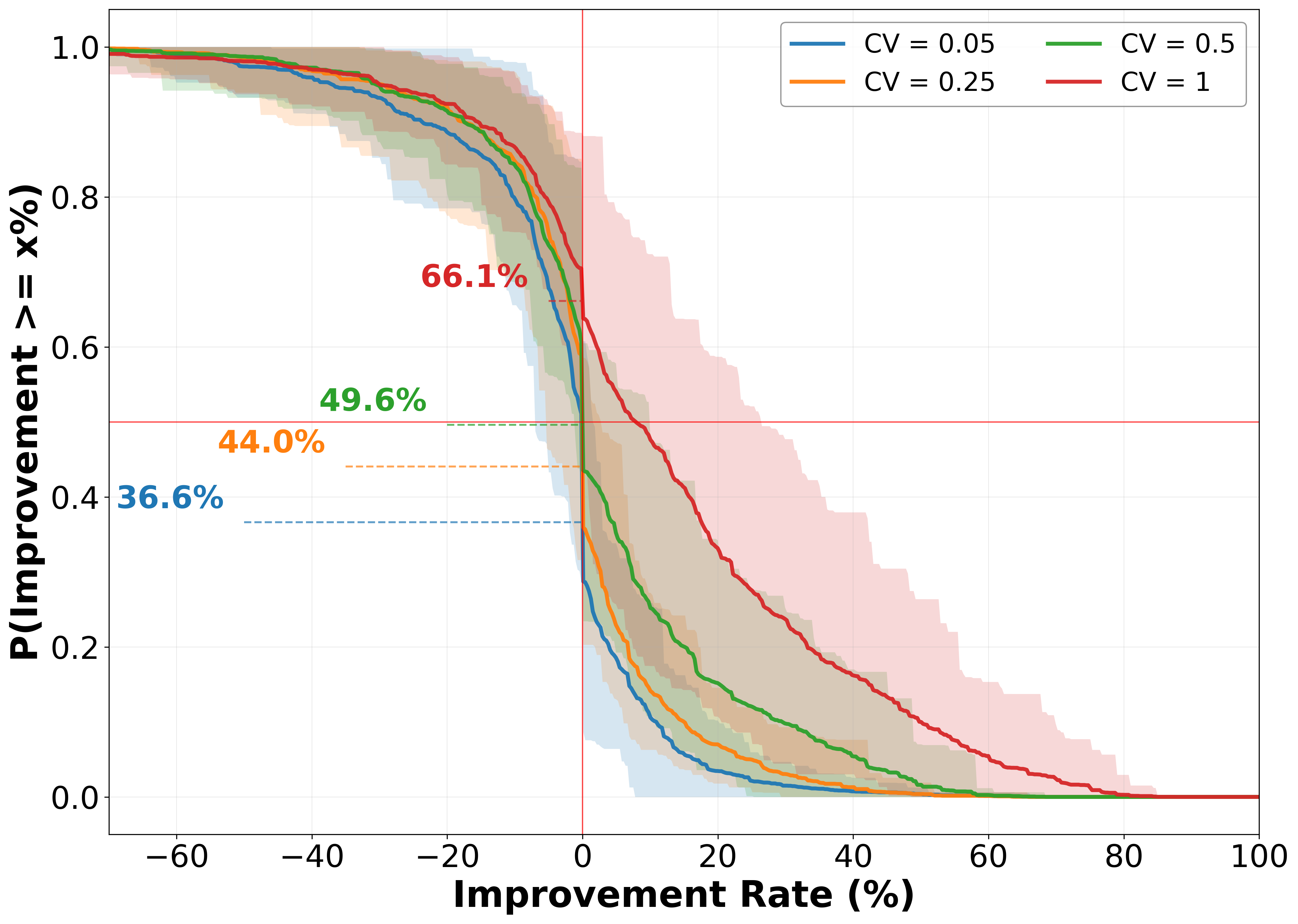}
        \caption{12 lines}
        \label{fig:appendix_improvement_exceedance_12routes}
    \end{subfigure}

    \caption{\rev{zg}{}{Trimmed-mean exceedance probability curves of OD-level improvement rates under different line quantities and variability levels. Each line is computed by first forming one per-seed quantile curve $Q_\omega(p)$ and then averaging the remaining curves after trimming extreme seed curves.}}
    \label{fig:appendix_improvement_exceedance_all_routes}
\end{center}

\section{Additional Beijing Network Visualisations}\label{app:beijing_networks}

\rev{zg}{}{Fig.~\ref{fig:app:beijingAllNetworks} complements the Beijing case-study visualisation in Fig.~\ref{fig:beijingRouteNetworks} by reporting the original network and the adapted line layouts at intermediate variability levels. These additional panels show that the redesign remains spatially targeted as $cv$ increases: lines are adjusted around affected corridors while preserving the broad structure of the original service network.}

\begin{center}
\captionsetup{type=figure}
\centering
\begin{subfigure}[t]{0.48\textwidth}
    \centering
    \includegraphics[width=\linewidth]{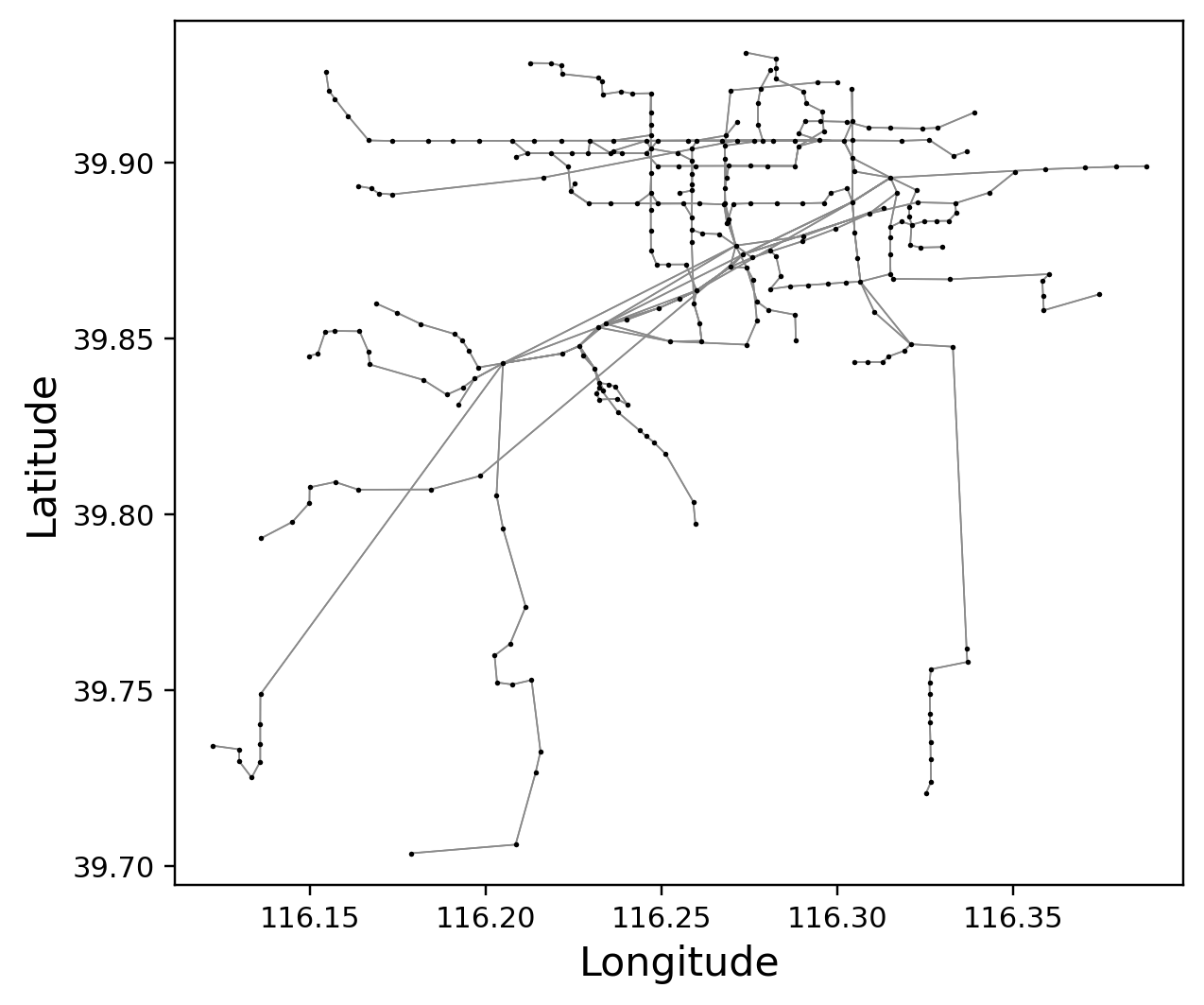}
    \caption{Original network}
    \label{fig:app:beijingOriginal}
\end{subfigure}
\hfill
\begin{subfigure}[t]{0.48\textwidth}
    \centering
    \includegraphics[width=\linewidth]{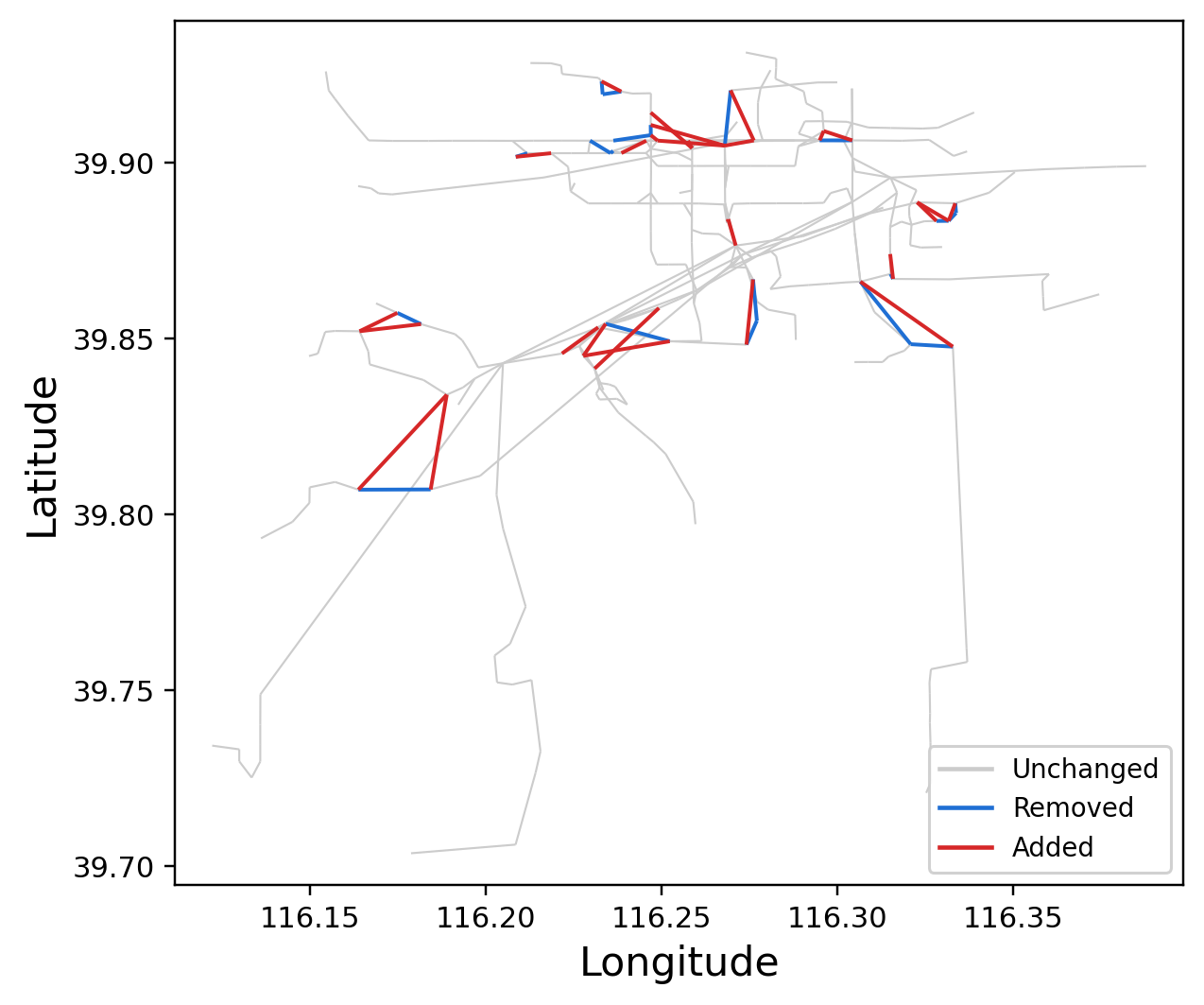}
    \caption{$cv = 0.05$}
    \label{fig:app:beijingCV0p05}
\end{subfigure}

\vspace{1em}

\begin{subfigure}[t]{0.48\textwidth}
    \centering
    \includegraphics[width=\linewidth]{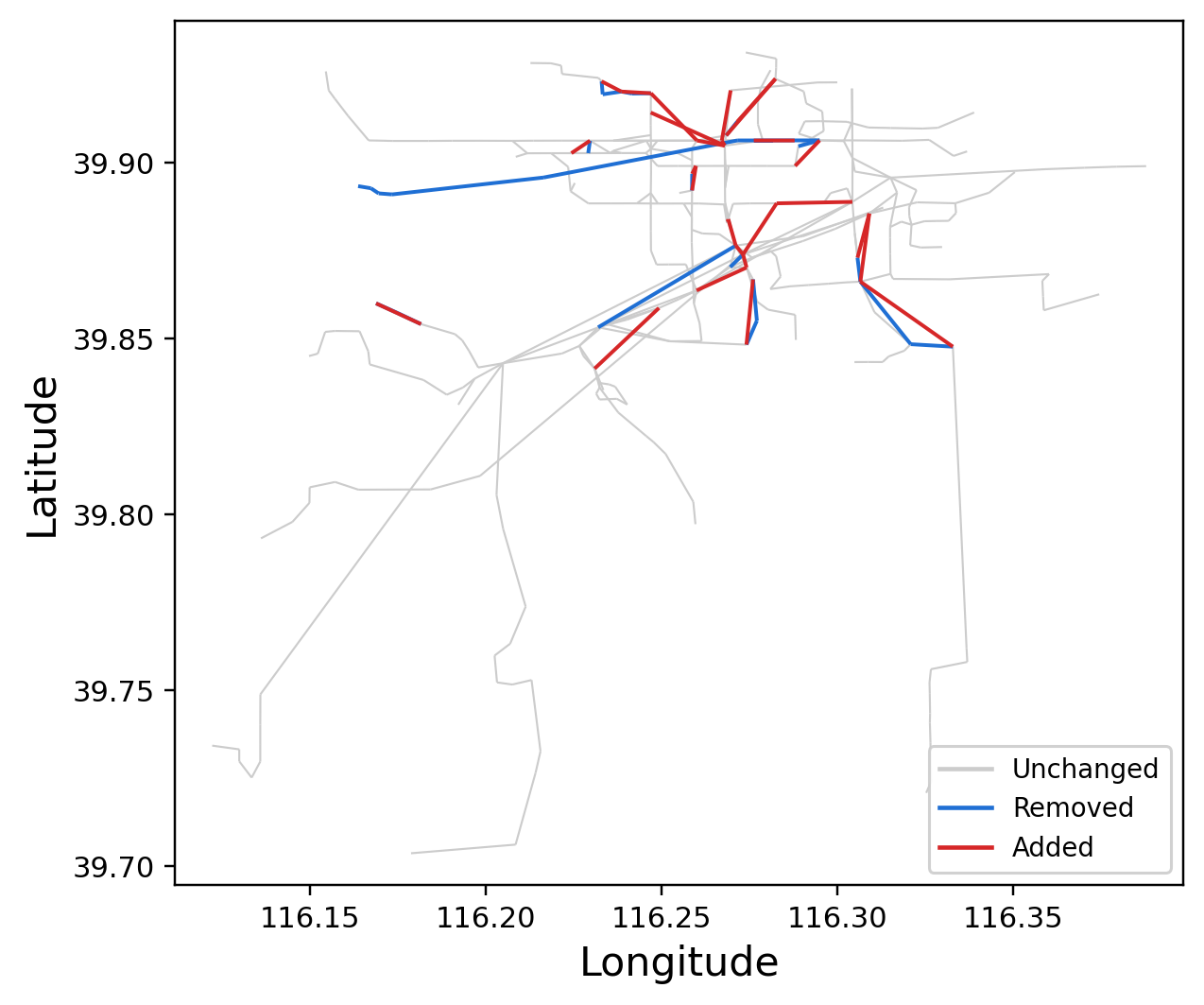}
    \caption{$cv = 0.25$}
    \label{fig:app:beijingCV0p25}
\end{subfigure}
\hfill
\begin{subfigure}[t]{0.48\textwidth}
    \centering
    \includegraphics[width=\linewidth]{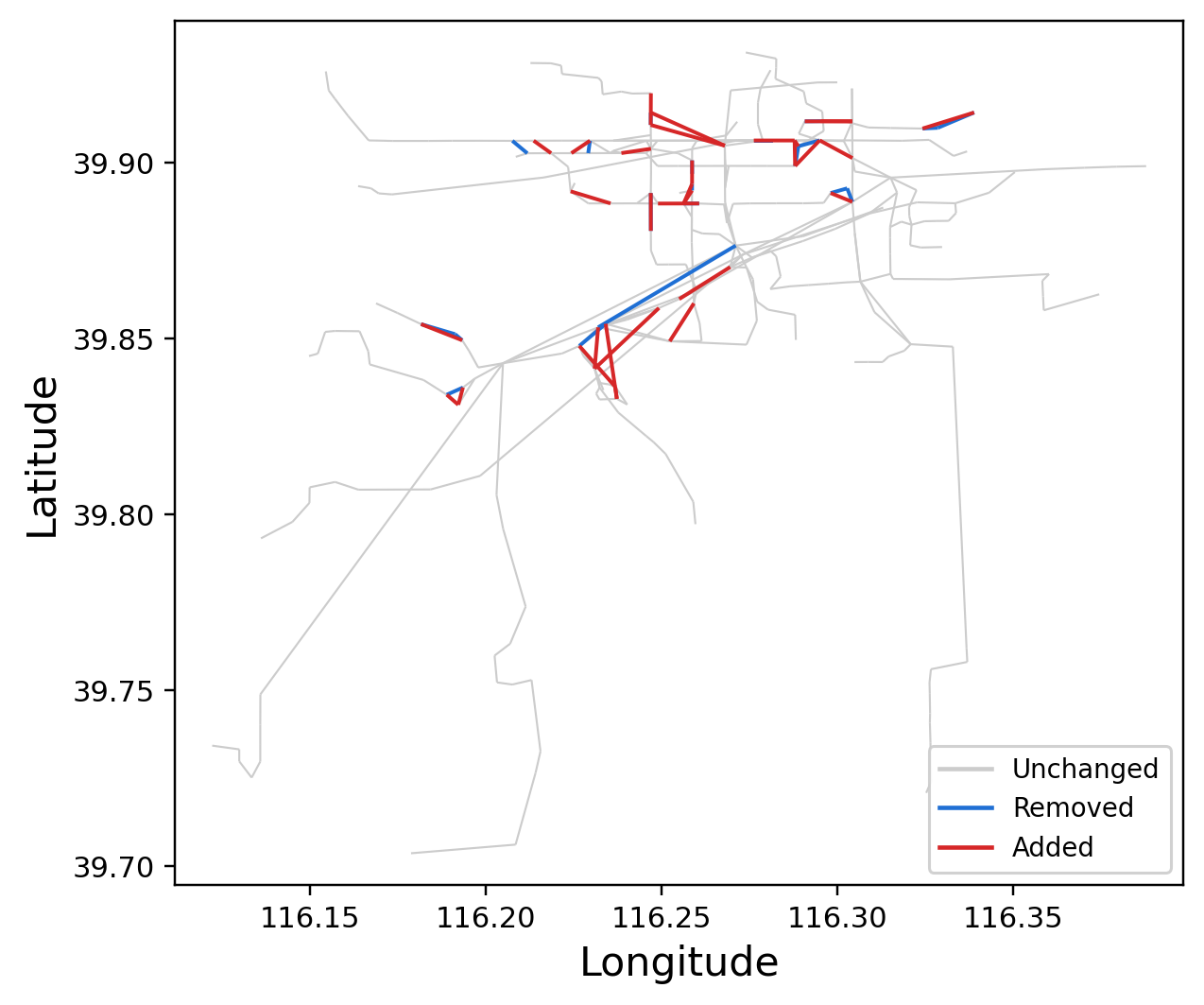}
    \caption{$cv = 0.5$}
    \label{fig:app:beijingCV0p5}
\end{subfigure}
\caption{\rev{zg}{}{Beijing PT networks: original design and adapted networks under $cv \in \{0.05, 0.25, 0.5\}$. The $cv=0$ and $cv=1$ cases are shown in the main text (Fig.~\ref{fig:beijingRouteNetworks}).}}
\label{fig:app:beijingAllNetworks}
\end{center}

%%%%%%%%%%%%%%%%%
\end{document}